\documentclass{amsart}
\usepackage{righttag}
\usepackage{epsf}
\def\hepsffile{\leavevmode\epsffile}

\theoremstyle{plain}
\newtheorem{thm}{Theorem}[subsection]

\newtheorem{lem}[thm]{Lemma}

\newtheorem{prop}[thm]{Proposition}

\theoremstyle{definition}
\newtheorem{defin}[thm]{Definition}

\newtheorem{emf}[thm]{}
\newtheorem{rem}[thm]{Remark}

\def\id{\protect\operatorname{id}}
\def\Im{\protect\operatorname{Im}}

\def\sign{\protect\operatorname{sign}}

\def\St{\protect\operatorname{St}}
\def\Sp{\protect\operatorname{Sp}}

\def\pr{\protect\operatorname{pr}}

\def\Z{{\mathbb Z}}

\def\R{{\mathbb R}}
\def\N{{\mathbb N}}

\def\1{\hbox{\rm\rlap {1}\hskip.03in{\rom I}}}
\def\Bbbone{{\rm1\mathchoice{\kern-0.25em}{\kern-0.25em}
	{\kern-0.2em}{\kern-0.2em}I}}

\def\p{\partial}

\begin{document} 
\title[Invariants of wave fronts on surfaces]
{Arnold-type invariants of wave fronts on surfaces}
\author[V.~Tchernov]{Vladimir Tchernov}
\address{D-MATH, HG G66.4, ETH Zentrum, CH-8092 Z\"urich, Switzerland}  
\email{chernov@math.ethz.ch}
\begin{abstract}
Recently Arnold's $\St$ and $J^{\pm}$ invariants of 
generic planar curves have been generalized to the
case of generic planar wave fronts. 
We generalize these 
invariants 
to the case of wave fronts on an arbitrary surface $F$. 
All invariants satisfying the axioms which  
naturally generalize the axioms used by Arnold 
are explicitly described. 
We also give an explicit formula for the finest order one 
$J^+$-type invariant of fronts on an orientable surface $F\neq S^2$.
We obtain necessary and
sufficient conditions for an invariant of nongeneric fronts with one 
nongeneric singular point to be the Vassiliev-type 
derivative of an invariant of generic fronts.
As a byproduct, we 
calculate all homotopy groups of the space of Legendrian immersions of 
$S^1$ into the spherical cotangent bundle of a surface.
\end{abstract}

\maketitle

By a surface we mean a smooth two-dimensional Riemannian manifold possibly 
with a boundary. The codimensions of all strata are calculated with respect
to the space of all fronts.

\section{Introduction}
\begin{figure}[htbp]
 \begin{center}
  \epsfxsize 4cm
  \hepsffile{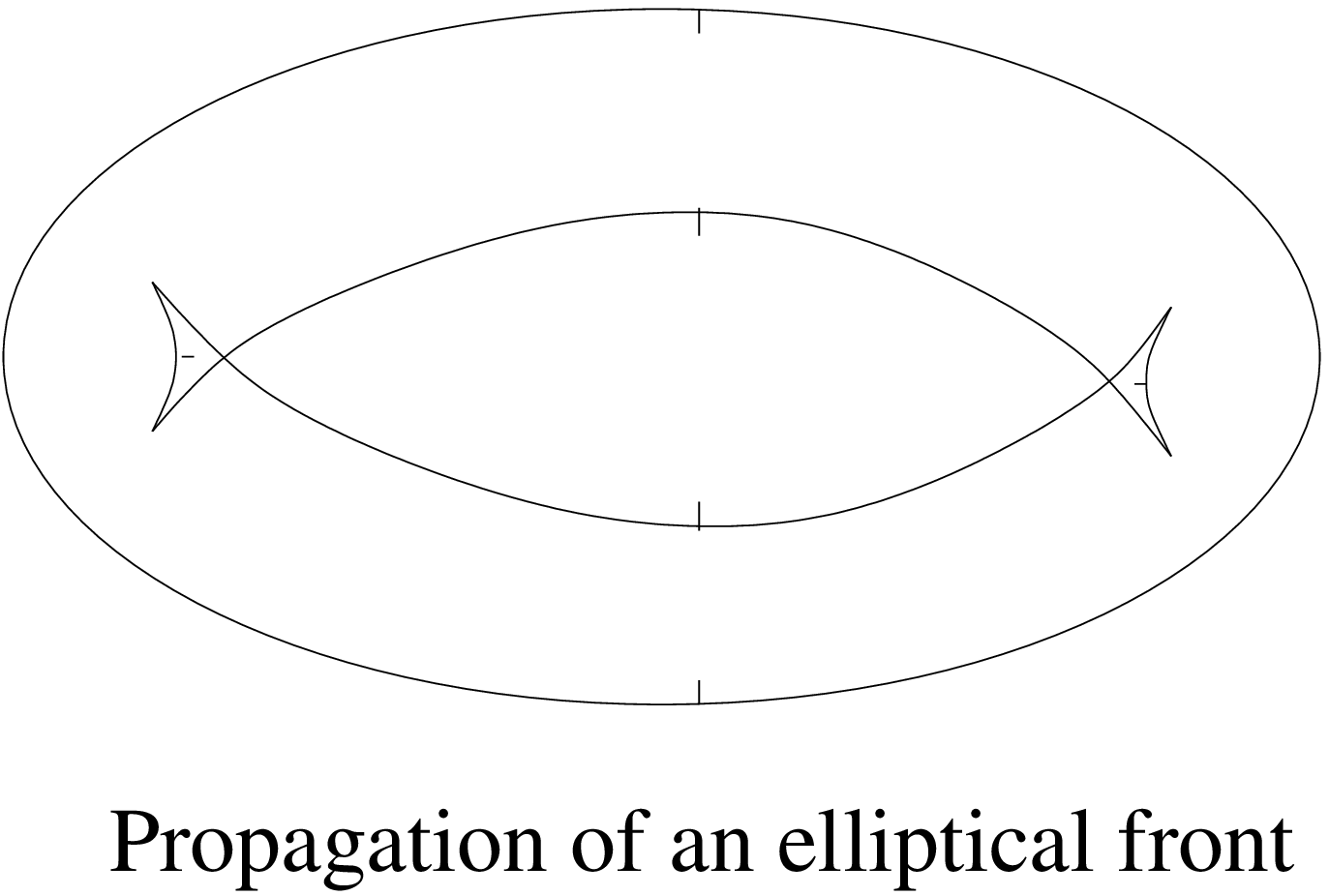}
 \end{center}
\caption{}\label{equidistant.fig}
\end{figure}

Consider an initially smooth closed generic 
wave front $L$ propagating on a surface $F$. 
If the matter forming the surface is uniform and isotropic, then at each
moment of time the propagating front forms an equidistant of $L$, i.e. 
the family of points of geodesics normal to $L$ located at the same distance 
from $L$. The propagating wave front is in general not smooth and has
semicubical cusp points as its singularities, see
Figure~\ref{equidistant.fig}. During the propagation
there are instances when the front has nongeneric singularities. 
The singularities 
arising during the propagation of a generic wave front are a 
triple point, a cusp crossing a
branch, a point of degree $4/3$ (this is a moment of birth of two cusp
points), and a self-tangency point at which the two
tangent branches propagate in opposite directions, see
Figure~\ref{codimone.fig}.
The self-tangency point
at which the branches are propagating in the same direction does not occur
for the simple physical reason that 
if it appears on the front, then the two 
branches stay tangent during the whole propagation process and hence were
tangent on the initial front.  

A line tangent to a front has a natural coorientation (transversal
orientation) given by the direction of propagation of the front.
A front $L$ can be naturally lifted to a curve $l$ in the spherical
cotangent bundle $ST^*F$ of $F$. (A point of the front is mapped to the point 
in $ST^*F$ corresponding to a functional which is zero on the line tangent
to the front at the point and positive in the coorienting 
half-plane of the tangent plane.) 

The space $ST^*F$ has a natural contact structure given by the distribution
of hyperplanes. The Huygens principle implies that at each moment of time 
the lifting of an initially smooth wave front is
a smooth knot which is everywhere tangent to the distribution
of contact hyperplanes, i.e. it is a {\em Legendrian\/} knot in $ST^*F$.
Since it never happens that two points of the wave front coincide and have
the same direction of propagation, we see that the 
propagation of the front induces an isotopy of the Legendrian knot. In
particular wave fronts corresponding to nonisotopic Legendrian knots 
can not be obtained from each other under the propagation.  
These observations  
make the study of Legendrian curves in $ST^*F$ started by V. Arnold 
a very attractive subject.

We consider the space of Legendrian immersed curves in $ST^*F$. A front 
of a Legendrian immersed curve in $ST^*F$ is its projection to $F$ equipped
with the natural coorientation and orientation.  
(The Legendrian immersed curve is 
uniquely determined by its front.)
A front $L$ is called generic if its only 
singularities are transversal double points and semicubical cusp points.
Nongeneric fronts form a discriminant in 
the space $\mathcal L$ of all fronts. We consider four 
codimension one strata of the discriminant. They are formed by 
fronts with a triple point (triple point stratum); 
fronts with a cusp point crossing a 
branch of the front (cusp crossing stratum); 
fronts with a self-tangency point at which the
coorienting normals of the two branches are pointing to the same direction
(dangerous self-tangency stratum); and fronts with a self-tangency point at which the 
coorienting normals of the two branches are pointing to the opposite
directions (safe self-tangency stratum).

A sign is associated to a generic crossing of 
each of these strata. In~\cite{Arnoldtopology} V.~Arnold constructed 
$J^+$ and $J^-$ invariants of generic fronts on $\R^2$.  
They increase 
by a constant under a positive 
crossing of respectively dangerous and safe self-tangency strata 
and do not change under crossings of all other codimension one 
strata of the discriminant. F.~Aicardi~\cite{Aicardidiscr} and 
M.~Polyak~\cite{Polyak} independently constructed an invariant that
increases 
respectively by a constant and by one half of the constant under a
positive crossing of triple point and cusp crossing strata and does
not change under crossings of all other codimension one strata of
the discriminant. Aicardi denoted this invariant by $\Sp$ and Polyak by 
$\St'$. The normalizations they used for this invariant are different, 
namely $\Sp=4\St'$. (In this paper I use Polyak's definition for the
invariant.) These invariants are natural generalizations of
V.~Arnold's~\cite{Arnold} $J^+$,
$J^-$, and $\St$ invariants of generic immersions of $S^1$ to $\R^2$. 
They give a lower bound for the number of crossings of the corresponding 
parts of the discriminant that are necessary to transform one generic 
front on $\R^2$ to another. $J^+$ seems to be the most
interesting of the three invariants because, as it was explained above, 
it corresponds to an isotopy invariant of Legendrian knots.

In this paper we construct generalizations of these invariants 
to the case where $F$ is any
surface (not necessarily $\R^2$). We follow an approach similar to the one
that was used in~\cite{Tchernov1} to generalize Arnold's 
invariants of generic 
immersions of $S^1$ to $\R^2$
to the case of generic immersions of $S^1$ to an arbitrary surface $F$. 
(Results similar to those of~\cite{Tchernov1} 
were independently obtained by A.~Inshakov~\cite{Inshakov},~\cite{Inshakov2}.)

The fact that for most surfaces
the fundamental group is nontrivial allows us to decompose in a natural way 
each of the four strata of the discriminant into pieces. 
To generalize the $J^+$ invariant   
we take an integer-valued function $\psi$ on the set of pieces 
obtained from the dangerous self-tangency stratum 
and try to construct an invariant that increases by $\psi(P)$ 
under a positive
crossing of a piece $P$ of this stratum
and does not change under crossings of the other
codimension one strata. In an obvious sense, $\psi$ is a derivative of such
an invariant and the invariant is an integral of $\psi$.
Any integrable (in the sense above)
function $\psi$ defines this kind of an invariant up to an additive
constant. 
Similar constructions generalize $J^-$ and $\St'$.
%invariants. 

We introduce a condition on $\psi$ which is necessary and sufficient 
for the existence of such invariants. 
If the surface is orientable, then the conditions that correspond to the 
generalizations of $J^+$ and $J^-$ are automatically satisfied, 
and such an invariant
exists for any function $\psi$. For the generalization of $\St'$ the
condition is not trivial. We reduce it to a
simple condition on $\psi$ which is sufficient for the existence of such an
invariant. All these conditions are automatically 
satisfied in the case of orientation-reversing fronts. 
We also get a very general statement giving necessary and sufficient 
conditions for
an invariant of nongeneric fronts with one  
singular point (of codimension one) to be a derivative of an
invariant of generic fronts.

The pieces into which the dangerous self-tangency stratum is decomposed  
are in a natural one-to-one
correspondence with the connected components arising under the normalization
of the part of the discriminant containing dangerous self-tangency points.
This means that any order one invariant of the $J^+$-type 
can be obtained as an integral of some $\psi$. Analogous facts are true for 
our decompositions of triple point and safe self-tangency strata, provided
that the surface $F$ is orientable.
We introduce a finer way to subdivide the strata 
into pieces to obtain a similar result for nonorientable surfaces.

We give an explicit formula for an order one $J^+$-type invariant of fronts 
on orientable surfaces. For $F\neq S^2$ this invariant is the finest order 
one $J^+$-type invariant and it distinguishes every two fronts 
that one can distinguish using order one $J^+$-type invariants with values
in an Abelian group.

The proofs of the main theorems are based on certain properties 
of $\pi_1(\mathcal L)$. As a byproduct result we explicitly calculate all the
homotopy groups of $\mathcal L$ or, which is the same, of the space of
Legendrian immersions of $S^1$ into $ST^*F$.

\section{Invariants of planar fronts}
\subsection{Basic facts and definitions}\label{basicdefin}
A {\em coorientation\/} of a smooth hypersurface in a functional space is
a local choice of one of the two parts 
separated by this hypersurface in a neighborhood of any of its points.
This part is called {\em positive\/}.

A {\em contact element\/} on the manifold is a hyperplane 
in the tangent space to the manifold at a point.
For a surface $F$ we denote by $ST^*F$ the space of all
cooriented (transversally oriented) contact elements of $F$. This space is a
spherical cotangent bundle of $F$. Its natural contact structure is a
distribution of tangent hyperplanes given by a condition 
that a velocity vector of an incidence point of a contact element 
belongs to the element. 
The Riemannian structure on $F$ allows us to
identify $ST^*F$ with the spherical tangent bundle  $STF$ of $F$.  
We denote by $CSTF$ the space  
of directions in the planes of the contact structure of the manifold $STF$. 
(One can show that if $F$ is orientable, then $CSTF=STF\times S^1$. 
For nonorientable $F$ the fiberwise projectivization $PCSTF$ of $CSTF$ is
isomorphic to $STF\times S^1$.)

A {\em Legendrian curve\/} $l$ in $STF$ is a smooth immersion 
of an oriented circle to $STF$  
such that the velocity vector of $l$ at every point lies 
in the plane of the contact structure. We denote by $\mathcal M$ the space of 
all Legendrian curves in $STF$.

The $h$-principle proved
for the Legendrian curves by M.~Gromov~\cite{Gromov} says that 
$\mathcal M$ is weak homotopy equivalent to the
space $\Omega CSTF$ of all free loops in $CSTF$. The equivalence is
given by $h:\mathcal M\rightarrow \Omega CSTF$ that sends a point
on $S^1$ (parameterizing a Legendrian curve) to the direction of the velocity
vector of the curve at the point. 
The connected components 
of $\mathcal M$ admit a rather simple description. They are naturally identified 
with the connected components of $\Omega CSTF$ or, which is the same, with
the conjugacy classes of $\pi_1(CSTF)$.

We denote by $f_1\in \pi_1(CSTF)$ and by 
$f_2\in \pi_1(STF)$ the classes of oriented fibers of the natural 
$S^1$-fibrations $\pr^1:CSTF\rightarrow STF$ and $\pr^2:STF\rightarrow F$.
We denote by $\pr^1_*$ and $\pr^2_*$ the homomorphisms of the fundamental
groups induced by the fibrations.

For a Legendrian curve $l$ we denote by $L$ the corresponding {\em wave
front\/}, which is the naturally cooriented and oriented 
projection of $l$ to $F$.  
A wave front on $F$ can be
naturally lifted to the Legendrian curve in $STF$ 
by mapping a point of it to the direction of the coorienting normal at the 
point. 
We denote by $\mathcal L$ the space 
of all fronts on $F$. (Note that the spaces $\mathcal L$ and $\mathcal M$
are naturally homeomorphic.)

A front $L$ on a surface $F$ is said to be {\em orientation-preserving\/} 
if it represents an orientation-preserving loop in $F$, and it is said to be
{\em orientation-reversing\/} otherwise.

A {\em generic wave front\/} has only transversal double points and
semicubical cusp
points as its singularities.
The nongeneric fronts 
form a {\em discriminant\/} in $\mathcal L$.

\begin{figure}[htbp]
 \begin{center}
  \epsfxsize 12cm
  \hepsffile{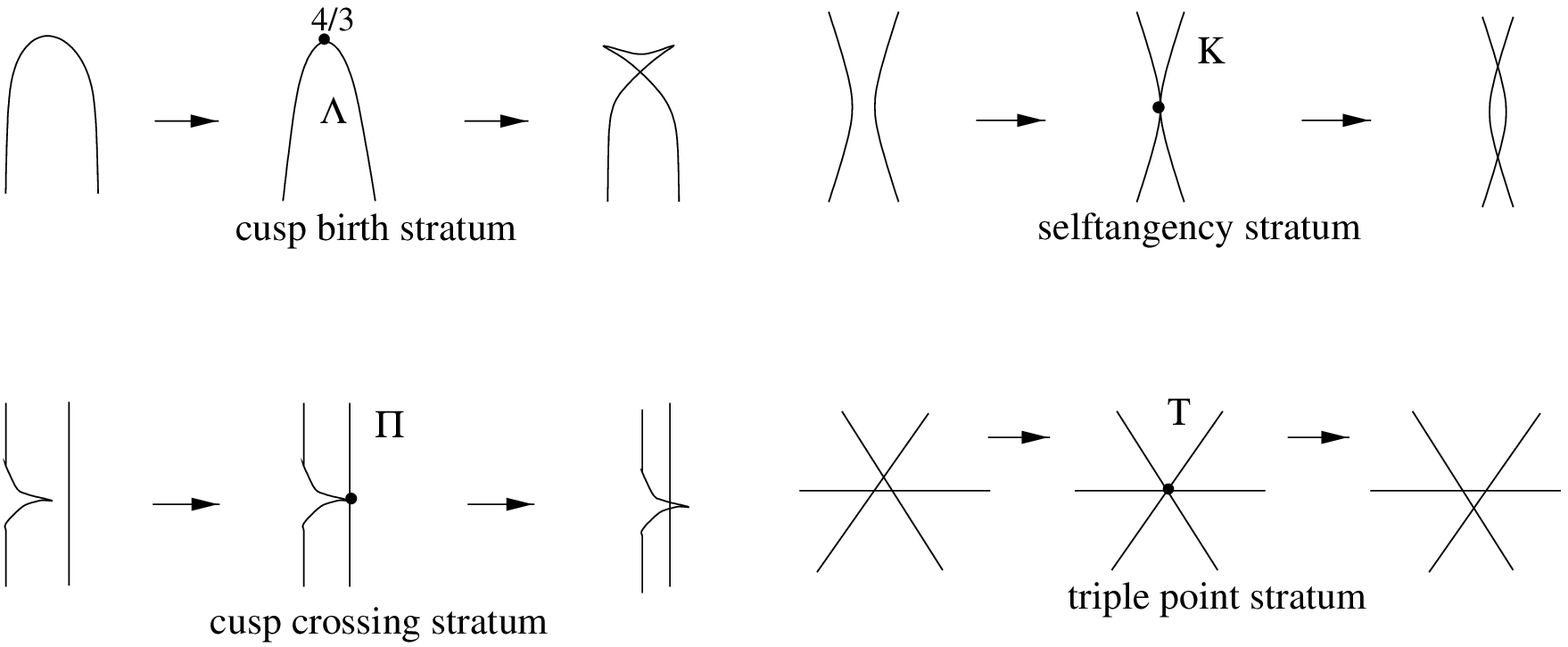}
 \end{center}
\caption{}\label{codimone.fig}
\end{figure}

\begin{thm}[Arnold~\cite{Arnoldsplit}] 
The codimension one strata of the discriminant of $\mathcal L$ are 
formed by fronts with one nongeneric singular point which is 
of one of the following types (see Figure~\ref{codimone.fig}):

1) A singular point of degree $\frac{4}{3}$. (This is
a moment of birth of two cusps.)
This stratum is denoted by $\Lambda$ and is called the cusp birth stratum.

2) A self-tangency point of order one of a front. This stratum is denoted
by $K$ and is called the self-tangency stratum.

3) A cusp point of a front passing through a branch.
(Here it is assumed that the line tangent to the front at
the cusp point is transverse to the branch.)
This stratum is denoted by $\Pi$ and is called the cusp crossing stratum.

4) A triple point of a front with pairwise transverse tangent lines at it. 
This stratum is denoted by $T$ and is called the triple point stratum.

\end{thm}

\begin{figure}[htbp]
 \begin{center}
  \epsfxsize 8cm
  \hepsffile{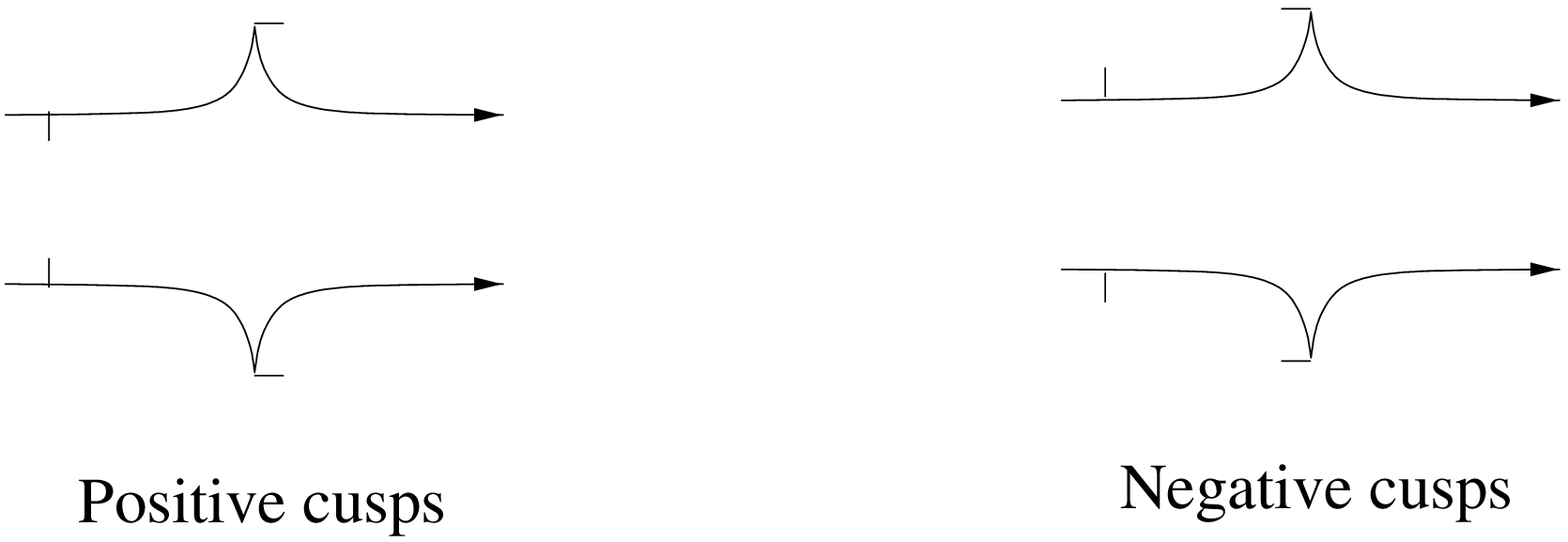}
 \end{center}
\caption{}\label{signcusp.fig}
\end{figure}

\begin{emf}\label{MaslovWhitney}{\em Whitney and Maslov indices.\/}
The {\em Whitney index\/} of a planar wave front is 
the total rotation number of the coorienting normal vector of the front.
The {\em Maslov index\/} of a generic planar 
wave front is the difference between the
number of positive and negative cusps. A cusp is said to be {\em
positive\/}
if the branch of the front going away from the cusp belongs to the
coorienting half-plane, see Figure~\ref{signcusp.fig}. 
A cusp is said to be {\em negative\/} otherwise. 

Whitney and Maslov indices of a front $L$ are
denoted by $\omega(L)$ and $\mu(L)$ respectively. Both these indices do not
change under a {\em regular homotopy\/} of a front, which is the projection 
of a homotopy in the class of Legendrian immersed curves in $ST\R^2$.
(Note that the Maslov index of a front on an arbitrary
surface $F$ is well defined.)
\end{emf}

The following theorem can be found in~\cite{Arnoldtopology}.
\begin{thm}
Two planar wave fronts $L_1$ and $L_2$ 
can be transformed to
each other by a regular homotopy if and only if $\omega(L_1)=\omega(L_2)$
and $\mu(L_1)=\mu(L_2)$.
\end{thm}

\subsection{Invariants $J^+$, $J^-$, and $\St$}

\begin{defin}[of the sign of a crossing of the $K$-stratum, 
Arnold~\cite{Arnoldsplit}]\label{signdirect}
%(Arnold~\cite{Arnoldsplit}.)
A self-tangency of a front is called {\em direct\/} 
if the velocity vectors of the two tangent branches have the 
same direction. A self-tangency is called {\em inverse\/} otherwise.
The $K$-stratum is decomposed into direct self-tangency and inverse
self-tangency parts.

A transversal crossing of the direct self-tangency part of the $K$-stratum 
is said to be {\em positive\/} 
if it increases (by two) the number of the
double points of the front. It is called {\em negative\/} otherwise.
A transversal crossing of the inverse self-tangency part of the $K$-stratum 
is said to be {\em positive\/} 
if it decreases (by two) the number of the
double points of the front. It is called {\em negative\/} otherwise.
\end{defin}

\begin{defin}[of $K^+$- and $K^-$-strata, Arnold~\cite{Arnoldsplit}]
%(Arnold~\cite{Arnoldsplit})
A self-tangency of a wave front is said to be {\em dangerous\/}
if the the coorientations of the tangent branches coincide, 
and it is said to be {\em safe\/} otherwise. 
(A front with a point of dangerous self-tangency
lifts to a Legendrian knot with a double point.) 
This relation induces a decomposition
of the $K$-stratum into the strata $K^+$  and $K^-$ of respectively
dangerous and safe self-tangencies.  
\end{defin}

\begin{defin}[of the sign of the
$T$-stratum crossing, Arnold~\cite{Arnold}]\label{signtriple}
A {\em vanishing triangle\/} is the triangle formed by the three branches
of the front corresponding to a subcritical or to a supercritical value of
the parameter near the triple point of the critical front.
The orientation of the front defines the cyclic
order on the sides of the vanishing triangle. (It is the order of the
visits of the triple point by the three branches.) Hence the sides of the
triangle acquire orientations induced by the ordering. But each side has
also its own orientation which may coincide or not with the
orientation defined by the ordering.
For a vanishing triangle we put $q\in\{0,1,2,3\}$ to be the
number of sides of it that are equally oriented by the ordering
and their direction. The {\em sign of the vanishing triangle\/} is $(-1)^q$.
The {\em sign\/} of a transversal crossing of the $T$-stratum 
is put to be the sign of the new born vanishing triangle. 
\end{defin}

\begin{figure}[htbp]
 \begin{center}
  \epsfxsize 9cm
  \hepsffile{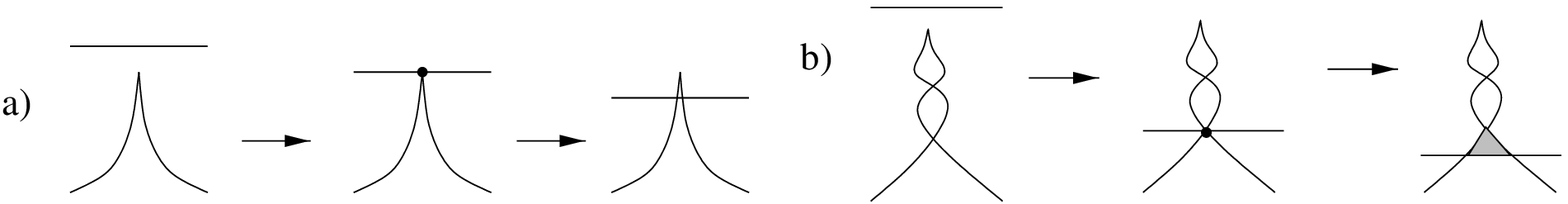}
 \end{center}
\caption{}\label{subst1.fig}
\end{figure}

\begin{defin}[of the sign of a crossing of the $\Pi$-stratum]\label{signcusp}
To define a sign of a transversal crossing of the $\Pi$-stratum 
shown in Figure~\ref{subst1.fig}.a,  we substitute a
cusp by a figure eight shape with a cusp on it. 
The sign of the crossing of
the $\Pi$-stratum is put to be the sign of the new-born
vanishing triangle shown in Figure~\ref{subst1.fig}.b.,
cf. Aicardi~\cite{Aicardidiscr} and Polyak~\cite{Polyak}.
(A similar way of defining the sign was suggested by
M.~Polyak~\cite{Polyakpriv}.)
\end{defin}

\begin{figure}[htbp]
 \begin{center}
  \epsfxsize 9cm
  \hepsffile{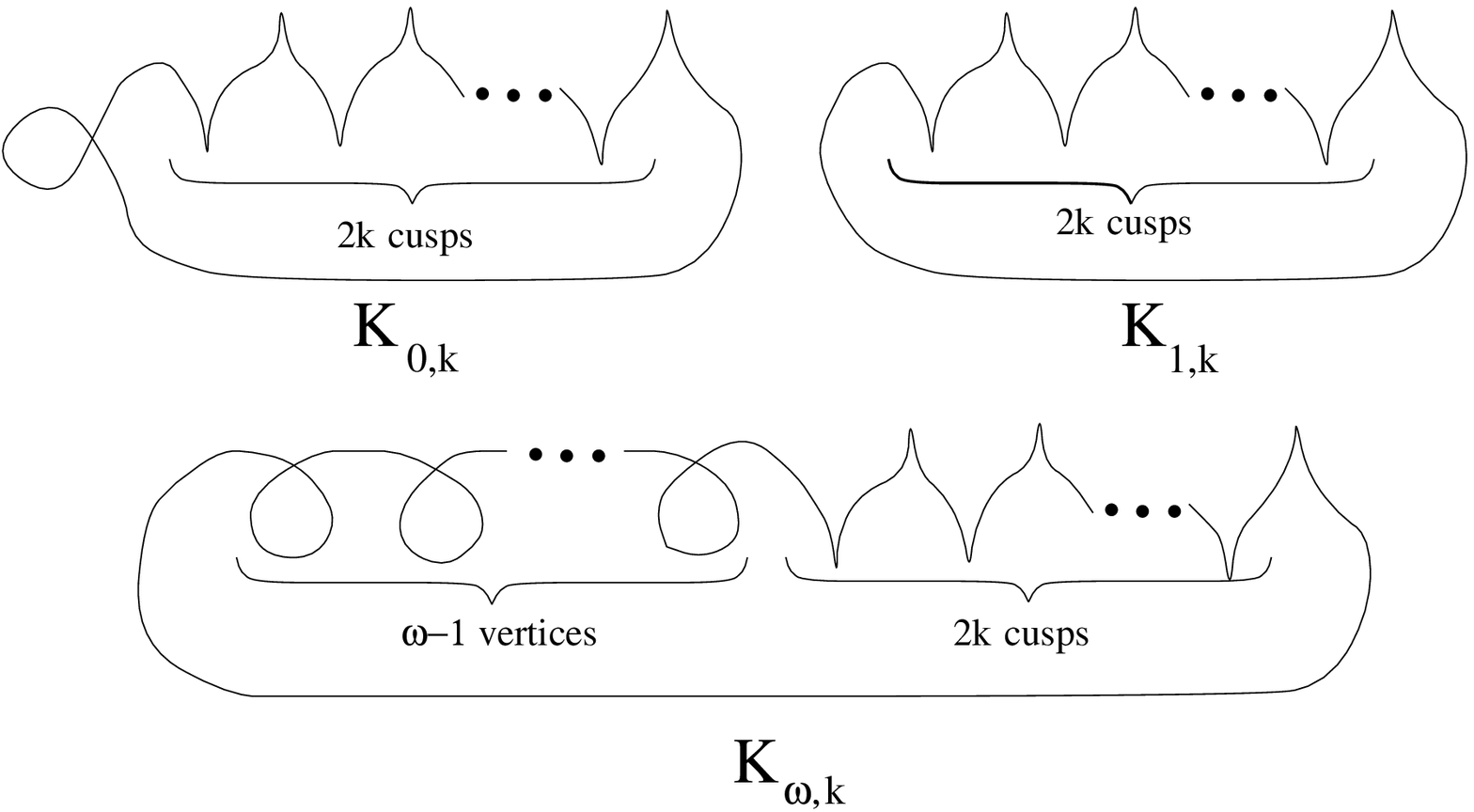}
 \end{center}
\caption{}\label{basic.fig}
\end{figure}

\begin{thm}[Aicardi~\cite{Aicardidiscr}, Arnold~\cite{Arnoldtopology},
Polyak~\cite{Polyak}]$\text{   }$

There exist three numbers $\St'(L)$, $J^+(L)$, and $J^-(L)$ assigned to a
generic planar front $L$ which are uniquely defined by the following   
properties:

1. $\St'(L)$, $J^+(L)$, and $J^-(L)$ are invariant under a regular homotopy  
in the class of generic fronts.

2. $\St'(L)$ does not change under crossings of $K^{\pm}$- and $\Lambda$-strata. 
It increases by one and by $\frac{1}{2}$ under positive crossings of
respectively $T$- and $\Pi$-strata.

3.  $J^+(L)$ does not change under crossings of $T$-, $\Pi$-, $\Lambda$-,
and $K^-$-strata, and it increases by two 
under a positive crossing of the $K^+$-stratum.

4.  $J^-(L)$ does not change under crossings of $T$-, $\Pi$-, $\Lambda$-,
and $K^+$-strata, and it increases by two
under a positive crossing of the $K^-$-stratum.

5. On the standard fronts $K_{\omega, k}$ (see
Figure~\ref{basic.fig}) $\St'(L)$, $J^+(L)$, and $J^-(L)$ take the
following values (independent of the choice of orientation and
coorientation of the standard fronts):

\begin{equation}
\St'(K_{0,k})=\frac{k}{2},\text{   } \St'(K_{\omega +1, k})=\omega+\frac{k}{2} 
\text{   }(\omega=0,1,2,\dots);
\end{equation}
\begin{equation}
J^+(K_{0,k})=-k,\text{   } J^+(K_{\omega +1, k})=-2\omega-k     
\text{   }(\omega=0,1,2,\dots);
\end{equation}
\begin{equation}
J^-(K_{0,k})=-1,\text{   }J^-(K_{\omega +1, k})=-3\omega       
\text{   }(\omega=0,1,2,\dots);
\end{equation}
where $k=0,1,2,\dots$.
\end{thm}

These values of the invariants on the standard fronts are chosen to make the
invariants additive under a certain connected summation of fronts and
independent of the choice of orientation and coorientation of the fronts.  

\section{Invariants of Fronts on Orientable Surfaces}

\subsection{Natural decomposition of the $K^+$-stratum}

\begin{emf}
Let $B_2$ be a bouquet of two oriented circles and $b$ its base point. 
Let $L\in K^+$ be a 
%generic 
front with a dangerous self-tangency point $q$. 
It can be lifted to the mapping $\bar L:S^1 \rightarrow STF$
that sends $p\in S^1$ to the point in $STF$ corresponding to the
direction of the
coorienting normal of $L$ at $L(p)$. (Note that $q$ lifts to a double point
$\bar q$ of
$\bar L$.) 

Let $\alpha :S^1 \rightarrow B_2$ be a continuous mapping such that:
\begin{description}
\item[a] $\alpha (\bar L^{-1}(\bar q))=b$;
\item[b] $\alpha$ is injective on the complement of $\bar L^{-1}(\bar q)$;
\item[c] The orientation of $B_2 \setminus b$ induced by $\alpha$ coincides 
with the orientation of the circles of $B_2$.
\end{description}

The mapping $\phi:B_2\rightarrow STF$ such that $\bar L=\phi \circ \alpha$
is called an {\em associated\/} with $L$ mapping of $B_2$.
(Note that the free homotopy class of a mapping of $B_2$ to $STF$ realized by 
$\phi$ is well defined modulo the orientation-preserving automorphism of
$B_2$ interchanging the circles.)
\end{emf}

\begin{defin}[of $K^+$-equivalence]
We say that $L_1\in K^+$ and $L_2\in K^+$  are
{\em $K^+$-equivalent\/} if there exist mappings of $B_2$ 
associated with the two of them that are free homotopic.
This equivalence relation induces a decomposition of the $K^+$-stratum  
into parts corresponding to different $K^+$-equivalence classes. 
We denote by $[L^+]$ the $K^+$-equivalence class corresponding to 
$L\in K^+$ and by $\mathcal K^+$ the set of all $K^+$-equivalence classes.
\end{defin}

\subsection{Axiomatic description of $\overline{J^+}$}

A natural way to introduce $J^+$-type invariant of generic wave fronts
on a surface $F$ is to take a function $\psi:\mathcal K^+\rightarrow \Z$ and to
try to construct an invariant of generic fronts from a fixed 
connected component $\mathcal C$ of
the space $\mathcal L$ (of all fronts on $F$) such that:

1. It increases by $\psi ([L^+])$
under the positive crossing of the part of the $K^+$-stratum 
that corresponds to a $K^+$-equivalence class $[L^+]$.

2. It does not change under crossings of $K^-$-, $T$-, $\Pi$-, and
$\Lambda$-strata of the discriminant.

If for a given function $\psi:\mathcal K^+\rightarrow \Z$ 
there exists such an invariant of wave fronts
from $\mathcal C$, then
we say that there exists a $\overline {J^+}$ invariant of
fronts from $\mathcal C$ that {\em integrates\/} $\psi$. 
Such $\psi$ is said to be $\overline{J^+}$-{\em integrable\/} in $\mathcal C$.  

\begin{thm}\label{barJ+orient}
Let $F$ be an orientable surface, 
$\mathcal C$ a connected component of $\mathcal L$, and 
$\psi:\mathcal K^+\rightarrow \Z$ a function. 
Then there exists a unique (up to an additive constant) 
invariant $\overline{J^+}$ of generic fronts from $\mathcal C$
which integrates $\psi$.
\end{thm}

The Proof of this Theorem (see Section~\ref{pfbarJ+orient}) is based on 
Theorem~\ref{integrability}.

%\begin{rem}
Thus for orientable $F$ every $\psi:\mathcal K^+\rightarrow \Z$
is integrable in all connected components of $\mathcal L$.
However if $F$ is nonorientable, then such an invariant 
exists not for all functions $\psi$. 
In Theorem~\ref{integrability} we present a condition on $\psi$
which is necessary and sufficient for it to be $\overline{J^+}$-integrable 
in a fixed connected component of $\mathcal L$.
%\end{rem}

\begin{rem}\label{kalf}
Most likely the proof of Theorem~\ref{barJ+orient} in the case of 
$l\neq 1\in\pi_1(STF)$
can be obtained as a consequence of a version of 
E. Kalfagianni's~\cite{Kalfagianni} 
Theorem $3.7$ for framed knots (if one formulates and
proves this version). Similar remark holds for the 
$\overline {J^-}$ invariant of fronts on orientable surfaces,
see Theorem~\ref{barJ-orient}, and for the part of
statement~{\rm I} of Theorem~\ref{integrability} which is related to 
$\overline{J^{\pm}}$ invariants. Other statements of
Theorem~\ref{integrability} about $\overline{J^{\pm}}$ invariants can not be
obtained in this way. Statements of Theorems~\ref{barSTorient} 
and~\ref{integrability} about the existence of $\overline{\St'}$ invariants 
also can not be obtained in this way. 
\end{rem}

\begin{emf}\label{splitJ+}{\em Connection with the standard
$J^+$-invariant.\/}
Since $\pi_1(ST\R^2)=\Z$, there are countably many $K^+$-equivalence 
classes of
nongeneric planar fronts of fixed Whitney 
and Maslov indices.
(The Whitney and Maslov indices of a planar front $L$
define the connected component
of the space of planar fronts containing $L$.)
Thus the construction of $\overline{J^+}$ gives rise to a splitting of the
standard $J^+$ invariant of V.~Arnold. This is the splitting introduced by
V.~Arnold~\cite{Arnoldsplit} in the case of planar fronts of the 
zero Whitney index and generalized to the case of arbitrary 
planar fronts by F.~Aicardi~\cite{Aicardi}.
\end{emf}

\subsection{Natural decomposition of the $K^-$-stratum}

\begin{emf}
Let $B_2$ be a bouquet of two oriented circles and $b$ its base point. 
Let $L\in K^-$ be a front with a safe self-tangency point 
$q$. It can be lifted to the mapping $\bar L$
from the oriented circle to $PTF$ (the projectivized tangent bundle of F)
which sends $p\in S^1$ to the point in $PTF$ corresponding to the 
line normal to $L$ at $L(p)$. 
(Note that $q$ lifts to a double point
$\bar q$ of
$\bar L$.)

Let $\alpha :S^1 \rightarrow B_2$ be a continuous mapping such that:
\begin{description}
\item[a] $\alpha (\bar L^{-1}(\bar q))=b$;
\item[b] $\alpha$ is injective on the complement of $\bar L^{-1}(\bar q)$;
\item[c] The orientation of $B_2 \setminus b$ induced by $\alpha$ coincides 
with the orientation of the circles of $B_2$.
\end{description}

The mapping $\phi:B_2\rightarrow PTF$ such that $\bar L=\phi \circ \alpha$
is called an {\em associated\/} with $L$ mapping of $B_2$.
(Note that the free homotopy class of a mapping of $B_2$ to $PTF$ realized by 
$\phi$ is well defined modulo the orientation-preserving automorphism of
$B_2$ interchanging the circles.)
\end{emf}

\begin{defin}[of $K^-$-equivalence]
We say that $L_1\in K^-$ and $L_2\in K^-$ are
{\em $K^-$-equivalent\/} if there exist associated with the two of them
mappings of $B_2$ that are free homotopic. 
The $K^-$-stratum is naturally
decomposed into parts corresponding to different 
$K^-$-equivalence classes. 
We denote by $[L^-]$ the $K^-$-equivalence class corresponding to 
$L\in K^-$ and by $\mathcal K^-$ the set of all $K^-$-equivalence classes.
\end{defin} 

\subsection{Axiomatic description of $\overline{J^-}$}

A natural way to introduce $J^-$-type invariant of generic fronts
on a surface $F$ is to take a function $\psi:\mathcal K^-\rightarrow \Z$ and to
try to construct an invariant of generic wave fronts from a fixed 
connected component $\mathcal C$ of $\mathcal L$ such that:

1. It increases by $\psi ([L^-])$
under a positive crossing of the part of the $K^-$-
stratum 
corresponding to a $K^-$-equivalence class $[L^-]$.

2. It does not change under crossings of $K^+$-, $T$-, $\Lambda$-, and
$\Pi$-strata of the discriminant.

If for a given function $\psi:\mathcal K^-\rightarrow \Z$ 
there exists such an invariant of wave fronts
from $\mathcal C$, then
we say that there exists a $\overline{J^-}$ invariant of
fronts in $\mathcal C$ which {\em integrates\/} $\psi$. 
Such $\psi$ is said to be $\overline{J^-}$-{\em integrable\/} in $\mathcal C$.

\begin{thm}\label{barJ-orient}
Let $F$ be an orientable surface, 
$\mathcal C$ a connected component of $\mathcal L$, and 
$\psi:\mathcal K^-\rightarrow \Z$ a function. 
Then there exists a unique (up to an additive constant) 
invariant $\overline{J^-}$ of generic fronts from $\mathcal C$
which integrates $\psi$.
\end{thm}

The Proof of this Theorem is analogous to the Proof of
Theorem~\ref{barJ+orient} (see Section~\ref{pfbarJ+orient}) 
and is based on Theorem~\ref{integrability} (see also~\ref{kalf}).

%\begin{rem}
Thus for orientable $F$ every $\psi:\mathcal K^-\rightarrow \Z$
is integrable in all connected components of $\mathcal L$.
However, if $F$ is nonorientable, then such an invariant 
exists not for all functions $\psi$. 
In Theorem~\ref{integrability} we present a condition on $\psi$
which is necessary and sufficient for it to be $\overline{J^-}$-integrable 
in a fixed connected component of $\mathcal L$.
%\end{rem}

\begin{emf}\label{splitJ-}
{\em Connection with the standard $J^-$-invariant.\/}
Since $\pi_1(PT\R^2)=\Z$, there are countably many $K^-$-equivalence 
classes of
nongeneric wave fronts on $\R^2$
of the fixed Whitney and Maslov indices.
(Whitney and Maslov indices of a planar front $L$
define the connected component
of the space of planar fronts containing $L$.)
Thus the construction of $\overline{J^-}$ gives rise to a splitting of the
standard $J^-$ invariant of V.~Arnold. This splitting is analogous to the 
splitting of $J^+$ introduced by
V.~Arnold~\cite{Arnoldsplit} in the case of planar fronts of the zero Whitney 
index and generalized to the case of arbitrary planar wave fronts 
by F.~Aicardi~\cite{Aicardi}.
\end{emf}

\centerline{}

\subsection{Natural decomposition of $T$- and $\Pi$-strata}
\begin{emf}
Let $B_3$ be a bouquet of three oriented circles with a fixed
cyclic order on the set of them, and let $b$ be the base point of $B_3$. 
Let $L\in T$ be a front on $F$ with a triple point $q$. 

Let $\alpha:S^1\rightarrow B_3$ be a continuous
mapping such that:

a) $\alpha (L^{-1}(q))=b$;

b) $\alpha$ is injective on the complement of $L^{-1}(q)$;

c) The orientation induced by $\alpha$ on $B_3\setminus b$
coincides with the orientation of the circles of $B_3$;

d) The cyclic order induced on the set of circles of $B_3$ by traversing
$\alpha(S^1)$ according to the orientation of $S^1$
coincides with the fixed one.

The mapping $\phi:B_3\rightarrow F$ such that $L=\phi\circ\alpha$ is called
an {\em associated\/} with $L$ mapping of $B_3$.
(Note that the free homotopy class of the mapping of $B_3$ to $F$ realized
by $\phi$ is well
defined modulo an automorphism of $B_3$ that preserves the orientation and
the cyclic order of the circles.) 
\end{emf}

\begin{defin}[of $T$-equivalence]
We say that $L_1\in T$ and $L_2\in T$ are
{\em $T$-equivalent\/} if 
there exist associated with them mappings of $B_3$ which are
free homotopic. The $T$-stratum is naturally
decomposed into parts corresponding to different 
$T$-equivalence classes. 
Amazingly enough, the $T$-equivalence relation 
induces also a subdivision of the $\Pi$-stratum of the discriminant. 
To see it, one substitutes the cusp on $L\in \Pi$ by a small figure 
eight shape with a cusp on it (see
Figure~\ref{subst1.fig}). 
As a result of this operation $L\in \Pi$ changes to
a front $L '\in T$. (Note that $L$ and $L'$ belong to
the same component of $\mathcal L$.)
We take the $T$-equivalence class
of the front $L$ to be the $T$-equivalence class
of the front $L'$.
We denote by $[L]$ the $T$-equivalence class corresponding to $L\in T$ or to
$L\in \Pi$ and 
by $\mathcal T$ the set of all $T$-equivalence classes.
\end{defin}

\subsection{Axiomatic description of $\overline{\St'}$}
A natural way to introduce $\St'$-type invariants of generic wave fronts
on $F$ is to take $\psi:\mathcal T\rightarrow \Z$ and to
try to construct an invariant of generic wave fronts from a fixed 
connected component $\mathcal C$ of
the space $\mathcal L$ (of all fronts on $F$) such that:

1. It does not change under crossings of $K^{\pm}$- and $\Lambda$-
strata.

2. It increases by $\psi ([L])$
under a positive crossing of the part of the $T$-stratum that
corresponds to a $T$-equivalence class $[L]$.

3. It increases by $\frac{1}{2}\psi ([L])$
under a positive crossing of the part of the $\Pi$-stratum that
corresponds to a $T$-equivalence class $[L]$. (As it is explained 
in~\ref{halfchange}, one can not substitute $\frac{1}{2}$ by another 
constant and construct an invariant of this sort, unless $\psi$ is put to be
identically zero on all $T$-equivalence classes 
appearing on the $\Pi$-stratum.) 

If for a given function $\psi:\mathcal T \rightarrow \Z$ 
there exists such an invariant of wave fronts
from $\mathcal C$, then
we say that there exists a $\overline{\St'}$ invariant of
fronts from $\mathcal C$ which {\em integrates\/} $\psi$. 
Such $\psi$ is said to be $\overline{\St'}$-{\em integrable\/} in $\mathcal C$.  

Not all functions $\psi$ are integrable.
In Theorem~\ref{integrability} we present a condition on $\psi$
which is necessary and sufficient for it to be integrable.
In the case of orientable $F$ there is a  simple condition 
which is
sufficient for the integrability of $\psi$.

\begin{thm}\label{barSTorient}
Let $F$ be an orientable surface, 
$\mathcal C$ a connected component of $\mathcal L$, and
$\psi:\mathcal T\rightarrow \Z$ a function taking equal
values on any two $T$-equivalence classes such that:

a) The free homotopy classes of the
mappings of $B_3$ representing them 
are different by an orientation-preserving automorphism of
$B_3$ which changes the cyclic order of the circles.

b) The restrictions of the mappings representing these classes to
one of the circles of $B_3$ are homotopic to a trivial loop.

Then there exists a unique (up to an additive constant) 
invariant $\overline{\St'}$ of generic wave fronts from $\mathcal C$
which integrates $\psi$.
\end{thm}

The Proof of this Theorem is analogous to the Proof of
Theorem~\ref{barJ+orient} (see Section~\ref{pfbarJ+orient}) 
and is based on Theorem~\ref{integrability} (cf.~\ref{kalf}).

%\begin{rem}
Note, that if $F$ is orientable, then a function $\psi:\mathcal T\rightarrow \Z$
described in Theorem~\ref{barSTorient}
is integrable in all connected components of $\mathcal L$.

In Section~\ref{integrability} we present a condition on $\psi$
which is necessary and sufficient for it to be $\overline{\St'}$-integrable 
in a fixed connected component of $\mathcal L$.
%\end{rem}

\begin{emf}\label{splitST}
{\em Connection with the standard $\St'$-invariant.\/}
Since $\R^2$ is simply connected, 
there is just one $T$-equivalence class of
nongeneric wave fronts on $\R^2$.
Thus the construction of $\overline{\St'}$ does not give anything new in the
classical case of planar wave fronts. 
\end{emf}

\section{Necessary and sufficient conditions for the integrability of
functions}

\subsection{Obstructions for the integrability}\label{changealong}
Let $\psi:\mathcal T\rightarrow \Z$ be a function, and let 
$\gamma$ be a generic loop in a connected component $\mathcal C$ of $\mathcal L$.
Let $I_1$ and $I_2$ be the sets of moments when $\gamma$ crosses
respectively $T$- and $\Pi$-strata.  
Let $\{\sigma_i\}_{i\in I_1}$ and $\{\sigma_i\}_{i\in I_2}$ 
be the signs of the corresponding crossings of the strata, 
and let $\{s_i\}_{i\in I_1}$ and $\{s_i\}_{i\in I_2}$
be the $T$-equivalence classes corresponding to the parts of the strata
where the crossings occurred. 
We call 

\begin{equation}
\Delta_{\overline{\St'}}(\gamma)=\sum_{i\in I_1}\sigma_i\psi(s_i)+
\sum_{i\in I_2}\sigma_i\frac{1}{2}\psi(s_i)
\end{equation} 
the {\em change of $\overline{\St'}$ along $\gamma$\/}. 
If $\Delta_{\overline{\St'}}(\gamma)=0$, then $\psi$ is said to be 
{\em integrable along\/} $\gamma$. 
In a similar way we introduce
the notion of {\em integrability along\/} $\gamma$ for integer valued functions on 
$\mathcal K^+$ and on $\mathcal K^-$. 
(For this purpose we use the intersections of $\gamma$
with $K^+$- and $K^-$-strata respectively.) The {\em
changes\/} of $\overline{J^-}$ and of $\overline{J^+}$ along $\gamma$ 
are also defined in a similar way. 
(Using Lemma~\ref{onlyhomotop} and the versions of it for the 
$\overline{J^\pm}$ invariants one can verify that the change along 
$\gamma$ depends only on the homology class
realized by $\gamma$ in $H_1(\mathcal L, \Z)$.)

Clearly if a function $\psi$ is integrable in $\mathcal C$, then it is
integrable along any generic loop $\gamma\subset \mathcal C$.
In this section we describe two loops $\gamma_1$ and 
$\gamma_2$ in $\mathcal C$ such that integrability 
along them implies integrability in $\mathcal C$. In a sense, the changes along 
them are the only obstructions for the integrability.
The loop $\gamma_1$ is going to be well defined (and needed) only in the
case of $\mathcal C$ consisting of 
orientation-preserving fronts on $F$.
The loop $\gamma_2$ is going to be well defined (and needed) only in the
case of $F$ being a Klein bottle and $\mathcal C$ consisting of 
orientation-preserving fronts on it.

\begin{figure}[htbp]
 \begin{center}
  \epsfxsize 9cm
  \hepsffile{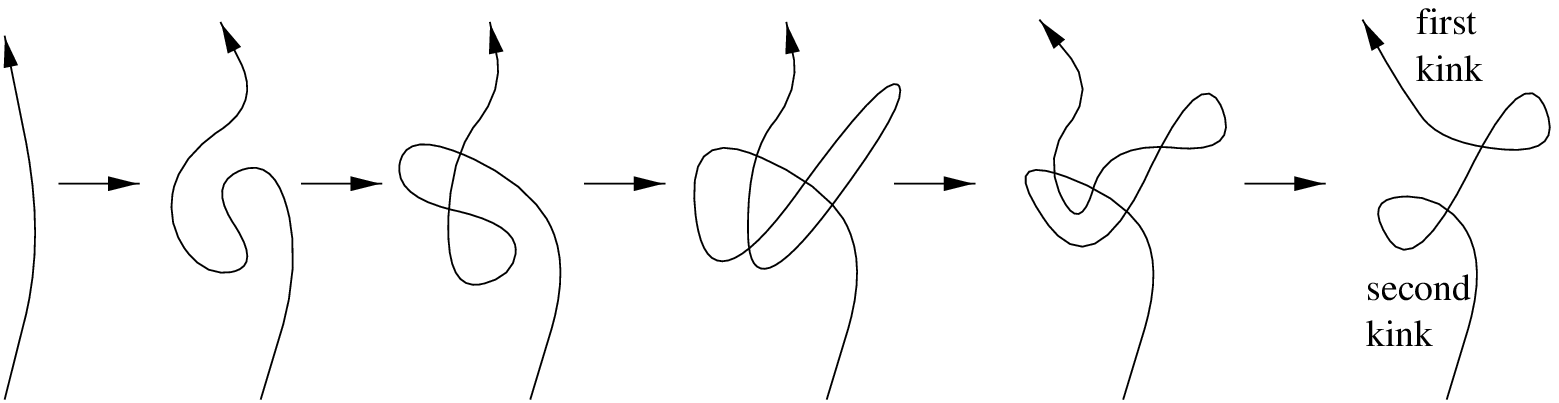}
 \end{center}
\caption{}\label{twokink.fig}
\end{figure}

\begin{figure}[htbp]
 \begin{center}
  \epsfxsize 8cm
  \hepsffile{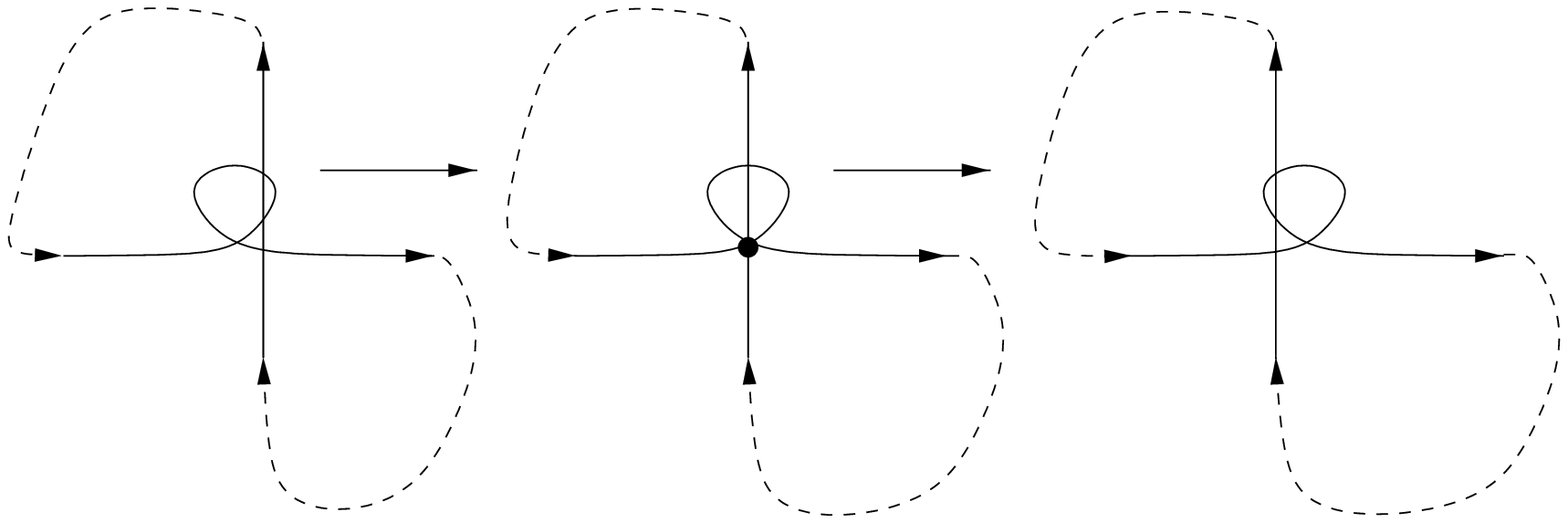}
 \end{center}
\caption{}\label{factorst1.fig}
\end{figure}

\begin{figure}[htbp]
 \begin{center}
  \epsfxsize 9cm
  \hepsffile{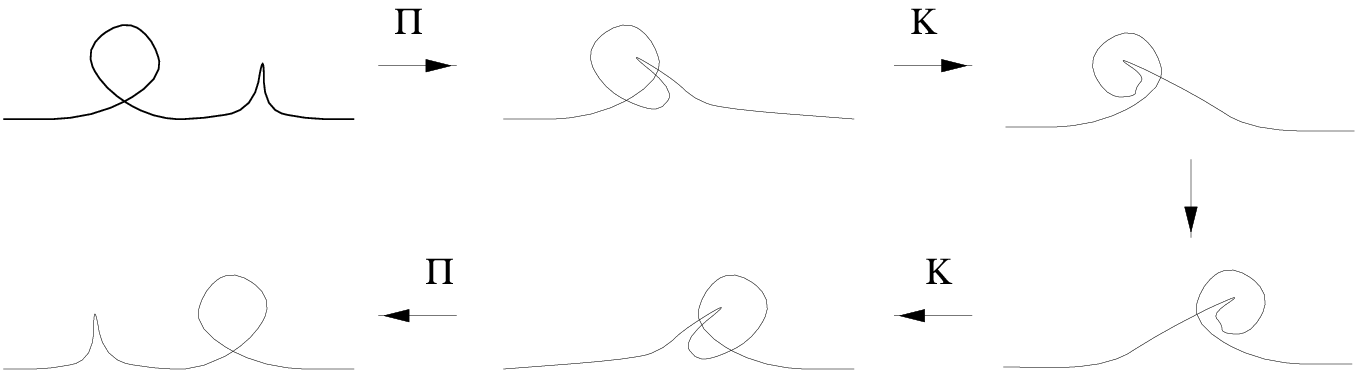}
 \end{center}
\caption{}\label{kinkcusp.fig}
\end{figure}

\begin{emf}\label{gamma1}{\em Loop $\gamma_1$.\/}
Let $\mathcal C$ be a connected component of $\mathcal L$ consisting of
orientation-preserving fronts, and let $L\in\mathcal C$ be a
generic front. Let $\gamma_1$ be the loop starting at $L$ which is 
described below.

Deform $L$ along a generic path $t$ in $\mathcal C$ to
get two opposite kinks,  as it shown in Figure~\ref{twokink.fig}. Make
the first kink small and slide it along the front till it comes back. 
We require the deformation to be such that at
each moment of time points of $L$ located outside of a small neighborhood
of the kink do not move.
(In Figure~\ref{factorst1.fig}
and Figure~\ref{kinkcusp.fig} it is shown how the kink passes through 
a neighborhood
of a double and of a
cusp point.)
Finally deform $L$ to its original shape along $t^{-1}$.
\end{emf}

\begin{emf}\label{gamma2}{\em Loop $\gamma_2$.\/}
Let $L$ be a generic orientation-preserving
front on the Klein bottle $K$. Let $\gamma_2\subset \mathcal C$ be the loop
starting at $L$ that is constructed below.

Consider $K$ as a quotient of a rectangle modulo the identification on its
sides shown in Figure~\ref{klein1.fig}. Let $p$ be the orientation
covering $T^2\rightarrow K$. There is a loop $\alpha$ in the space of all
autodiffeomorphisms of $T^2$
which is the sliding of $T^2$ along the unit
vector field parallel to the lifting
of the curve $c\subset K$
(see Figure~\ref{klein1.fig}). Since $L$ is an orientation-preserving front
it can
be lifted to a front $L '$ on $T^2$.
The loop $\gamma_2$ is the composition of $p$
and of the sliding of $L '$ induced by $\alpha$. 

\begin{figure}[htbp]
 \begin{center}
  \epsfxsize 3cm
  \hepsffile{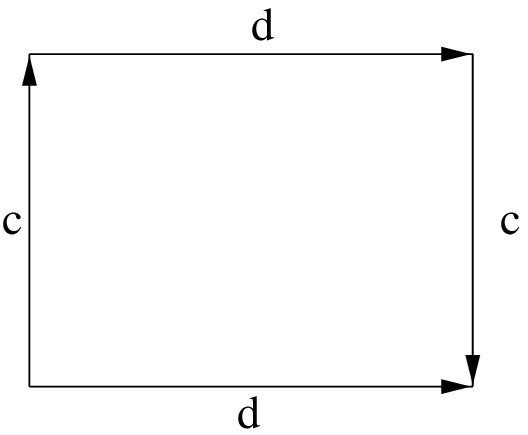}
 \end{center}
\caption{}\label{klein1.fig}
\end{figure}

\end{emf}

\subsection{Main Integrability Theorem}
Now we are ready to formulate the main integrability theorem.

\begin{thm}\label{integrability}
Let $F$ be a surface (not necessarily compact or orientable), $\mathcal C$ 
a connected component of $\mathcal L$, and $L\in \mathcal C$ a generic front.  
Let $\psi_1:\mathcal T\rightarrow \Z$, $\psi_2:\mathcal K^+\rightarrow \Z$ and
$\psi_3:\mathcal K^-\rightarrow \Z$ be functions. 

\begin{description}
\item[\textrm{I}] If $\mathcal C$ consists of orientation-reversing 
fronts on $F$, then there exists a 
$\overline{\St'}$ (resp.~$\overline{J^+}$, resp.~$\overline{J^-}$) 
invariant which integrates $\psi_1$ (resp.~$\psi_2$, resp.~$\psi_3$) 
in $\mathcal C$.

\item[\textrm{II}] If $\mathcal C$ consists of orientation-preserving fronts 
and $F\neq K$, then 
the condition that $\psi_1$ (resp.~$\psi_2$, resp.~$\psi_3$) 
is integrable along the loop $\gamma_1$ starting at $L$ 
is necessary and sufficient for the existence of a
$\overline{\St'}$ (resp.~$\overline{J^+}$, resp.~$\overline{J^-}$) 
invariant which integrates $\psi_1$ (resp.~$\psi_2$, resp.~$\psi_3$) 
in $\mathcal C$.

\item[\textrm{III}] If $\mathcal C$ consists of orientation-preserving fronts
and $F=K$, then the condition that $\psi_1$ (resp.~$\psi_2$,
resp.~$\psi_3$) is integrable along loops $\gamma_1$ and 
$\gamma_2$ starting at $L$ 
is necessary and sufficient for the existence of a
$\overline{\St'}$ (resp.~$\overline{J^+}$, resp.~$\overline{J^-}$)
invariant which integrates $\psi_1$ (resp.~$\psi_2$, resp.~$\psi_3$) 
in $\mathcal C$.
\end{description}
\end{thm}

For the proof of Theorem~\ref{integrability} see
Section~\ref{pfintegrability} (cf. also~\ref{kalf}).

\begin{emf}\label{importrem}{\em Remarks to Theorem~\ref{integrability}.\/}
If an invariant from the statement of the Theorem exists, then it is unique up 
to an additive constant.

The choice of $L\in \mathcal C$ does not matter and to check integrability of
a given function it is easier to take $L\in \mathcal C$ that
already has a small kink. Clearly all the crossings of the discriminant
under the deformation $\gamma_1$ of $L$ 
occur when the kink passes through a neighborhood of a double point 
or of a cusp of $L$. It easy to verify that inputs into $\Delta(\gamma_1)$
of the crossings of the discriminant corresponding to the kink passing
through a neighborhood of a cusp occur in pairs and cancel out. A kink
passes through a neighborhood of a double point $x$ twice (once along each 
of the two intersecting branches). One can verify that if $x$ separates the
front into two orientation-preserving loops, then the inputs into 
$\Delta_{\overline{J^{\pm}}}(\gamma_1)$ corresponding to $x$ also cancel out.

A straightforward modification of the proof of Theorem~\ref{integrability} 
shows that it holds for $\psi_1$, $\psi_2$, and $\psi_3$
taking values in any torsion free Abelian group.
\end{emf}

\subsection{A very general integrability Theorem}\label{generaltheorem} 
From the proof of Theorem~\ref{integrability} one can see that for any 
$\alpha \in H_1({\mathcal C}, \Z)$ a certain (nonzero) multiple of $\alpha$
is equal to a sum of a homology class realizable by a loop not 
crossing the discriminant and of a linear combination of homology classes 
realized by $\gamma_i$, $i\in\{ 1,2,3\}$ (see~\ref{gamma3} for the
definition of $\gamma_3$). (If $\gamma_i$ for some $i\in\{
1,2,3\}$ is not defined in $\mathcal C$, then it does not participate in the 
linear combination.)
A very important consequence of this is the following very general statement. 

Let $\Sigma$ be the
discriminant in $\mathcal C$, $G$ an Abelian group, and $\chi$
a $G$-valued invariant of generic fronts in $\mathcal C$, that is a mapping
from the set of the connected components of ${\mathcal C}\setminus \Sigma$
to $G$. (The condition that $\chi$ takes values in an Abelian group is not very
restrictive, since $\chi$ with values in an abstract set $S$ can be 
viewed as an invariant
taking values in $\Z[S]$ an Abelian group of abstract 
finite integer linear combinations of the elements of $S$.)
 
The invariant $\chi$ gives rise to the invariant $\chi'$ 
of nongeneric fronts with one singular point (which is singular of 
codimension one), that is a mapping to
$G$ from the set of connected components of $\Lambda$-, $K^{\pm}$-, $\Pi$-, 
and $T$-strata. The value of $\chi'$ on a component of a stratum is set to be 
the difference 
between the values of $\chi$ on the positive and negative sides of it. 
(The positive side of the $\Lambda$-stratum is the one with more cusps.)
In an obvious sense $\chi'$ is a derivative of $\chi$.

A very natural question is the following (cf. Kalfagianni~\cite{Kalfagianni}):
does a given $G$-valued 
function $\chi'$ on the set of connected components of $\Lambda$-,
$K^{\pm}$-, $\Pi$-, and $T$-strata correspond to some $\chi$ under the
construction above? (Is the integral of $\chi '$ well defined?) 

For a generic loop $\gamma$ in $\mathcal C$ put 
$\Delta_{\chi}(\gamma)=\sum_{x\in \gamma\cap\Sigma}\sign(x)\chi'(x)$.
Clearly the necessary condition for the existence of $\chi$ 
is that $\Delta_{\chi}(\gamma)=0$ for 
any small generic loop $\gamma$ going around a
codimension two stratum of $\mathcal C$ (see~\ref{strata} for
the list of the strata). We call this condition 
{\em a local integrability condition\/}.
If the local integrability 
condition is satisfied, then the change of $\Delta_{\chi}(\gamma)$ 
depends only on the homology class of the generic loop $\gamma$. 
Using the observation above one can easily modify the 
proof of Theorem~\ref{integrability} 
to get the following very general Theorem.

\begin{thm} Let $G$ be a torsion free Abelian group, $\mathcal C$ a
connected component of $\mathcal L$, and $\chi'$ a
$G$-valued invariant of nongeneric fronts from $\mathcal C$ 
whose only nongeneric singularity is one codimension one singular point.  
Then $\chi'$ is a derivative of an invariant $\chi$ of generic fronts 
from $\mathcal C$
if and and only if 
$\chi'$ satisfies the local integrability condition and 
$\Delta_{\chi}(\gamma_i)=0$ for those $i\in \{1,2,3\}$ for which 
$\gamma_i$ is well defined in $\mathcal C$. 
\end{thm}

Potential interesting application of this Theorem lie in the theory
of Legendrian knots in $ST^*F$. The invariants of Legendrian 
knots correspond to
$\chi'$ being identically zero on all the components of $\Lambda$-, $K^-$-, 
$\Pi$-,
and $T$-strata. In this case the local integrability condition appears to
be very simple. Using this Theorem one can easily obtain the
generalizations of other local invariants of planar fronts studied by
F.~Aicardi~\cite{Aicardidiscr}.

\section{Singularity theory interpretation of the invariants}
The invariants $\overline{\St'}$, $\overline{J^+}$, and $\overline{J^-}$ admit 
a rather simple singularity theory interpretation. Namely, the set of all 
$K^+$-equivalence classes appearing on the discriminant in $\mathcal C$
enumerates the components of the normalized dangerous self-tangency part of the
discriminant in $\mathcal C$. If $F$ is
orientable then similar facts are true for 
the sets of $K^-$- and $T$-equivalence classes. (See
Proposition~\ref{orient}.)
 
In this section we introduce finer versions of equivalence relations 
to obtain the complete classification of the components of the four parts of
the discriminant arising under the normalizations described below.
(This is done for all $F$, not necessarily orientable, see
Theorem~\ref{interpretation}.) We also formulate the corresponding versions 
of Theorem~\ref{integrability}

\subsection{Normalizations}
Let $S^1(2)$ be the configuration space of unordered pairs of distinct
points on
$S^1$. 
Let $\mathcal N^+$ be the subspace of
$S^1(2)\times\mathcal L$ consisting of $t\times L\in S^1(2)\times\mathcal L$ such
that $L$ maps the two points from $t$ to one point in $F$ and the
coorienting normals of $L$ at these two points have the same direction.
(This is a normalization of the part of the discriminant 
containing points of dangerous self-tangency.)

We say that $n_1^+, n_2^+\in \mathcal N^+$ are
$\overline{K^+}$-equivalent
if they belong to the same path connected component of $\mathcal N^+$.
Clearly for $L\in K^+$ 
there is a
unique $\overline{K^+}$-equivalence class associated with it. 
Thus the
$\overline{K^+}$-equivalence relation induces a decomposition of
the $K^+$-stratum.

Similarly we normalize the part of the discriminant containing 
safe self-tangency (resp. triple) points
and introduce the notion of $\overline {K^-}$- (resp. $\overline T$-) 
equivalence relation of fronts in $K^-$ (resp. in $T$).

We consider the following normalization of the closure 
of the $\Pi$-stratum.
Let $N$ be the closure of the subspace 
of $S^1\times S^1\times \mathcal L$
consisting of $t_1\times t_2\times L$ such that $t_1\neq t_2$, $L$ has
a cusp point at $t_1$, and $L(t_1)=L(t_2)$. We say that
$n_1, n_2\in \mathcal N$ are
$\overline{\Pi}$-equivalent
if they belong to the same path connected component of $\mathcal N$. Clearly,
for $L\in \Pi$ 
there is a                                                       
unique $\overline{\Pi}$-equivalence class associated with it. 
Thus the
$\overline{\Pi}$-equivalence relation induces a decomposition of the
$\Pi$-stratum. We denote by $\overline {\mathcal P}$ the set of all
$\overline {\Pi}$-equivalence classes. 

\begin{rem}
The reason why we treat 
cusp crossings differently from other codimension one 
singularities is that in the neighborhood of two of the codimension two 
strata of the discriminant the $\Pi$-stratum is not connected.
(These strata are $\Lambda \Lambda$ and 
$\Pi \Lambda$ in the notation of Theorem~\ref{strata}, see 
Figure~\ref{bifur5.fig}.) One can verify that
the two branches on the bifurcation diagrams 
of these singularities corresponding
to the $\Pi$-stratum belong to the same irreducible real algebraic 
curve, and hence it is natural to glue together the components of 
the $\Pi$-stratum around these codimension two strata. 
\end{rem}

\subsection{Description of the sets of 
$\overline {K^+}$-, $\overline {K^-}$-, $\overline T$-, and $\overline
\Pi$-equivalent fronts}
In this subsection we state Theorem~\ref{interpretation} 
which gives an explicit description of the sets of components 
of the four parts of the discriminant arising 
under the normalizations described above. These components are enumerated 
by the sets of $K^+_i$-, $K^-_i$-, $T_i$-, and $\Pi_i$-equivalence classes
introduced below.

We fix $d\in STF$ and denote by $\pi_1(STF)$ and $\pi_1(PTF)$ the groups
$\pi_1(STF, d)$ and $\pi_1(PTF, p(d))$. (Here $p:STF\rightarrow PTF$ is the
natural double covering.) The definitions introduced below are similar to the
ones introduced by Inshakov~\cite{Inshakov} in the case of immersed 
curves on a surface. For this reason we use a subscript $i$ in the notation 
of the arising equivalence classes.

\begin{defin}[of $K_i^+$-equivalence]
Let $R_{K^+}$ be the set consisting of all triples
$(\delta_1, \delta_2, i)\in \pi_1(STF)\oplus\pi_1(STF)\oplus\Z$
such that 
$i$ is even provided that $\pr^2(\delta_1\delta_2)$ is an
orientation-preserving loop in $F$ and odd otherwise.

The Legendrian curve $l$ corresponding to $L\in K^+$
has a double point separating it into two oriented loops. Deform $l$
preserving the double point so that the double point is located at $d$.
Choosing one of the two loops of $l$ we obtain
an ordered set of two elements $\delta_1, \delta_2\in \pi_1(STF)$. 
We also correspond to the
front its Maslov index $\mu(L)\in \Z$ that is even if and only if
$l=\delta_1\delta_2$ projects to an orientation-preserving loop in $F$.
Thus we obtain an element of $R_{K^+}$ corresponding to the deformed $l$.
There is a unique $K_i^+$-equivalence class of elements of $R_{K^+}$
corresponding to the undeformed $l$, where two elements of $R_{K^+}$ are 
$K_i^+$-equivalent if one can be transformed to the other by the consequent 
actions of the following groups (which all act trivially
on the last summand in $R_{K^+}$):

1. $\pi_1(STF)$ whose elements act via conjugation of the first two
summands in $R_{K^+}$. (This corresponds to the ambiguity in deforming $l$, so
that the double point is located at $d$.)  

2. $\Z_2$ which acts via the cyclic permutation of the first two summands.
(This corresponds to the ambiguity in the choice of one of the two 
loops of $l$.) 

The set of all $K_i^+$-equivalence classes of elements of $R_{K^+}$ is
denoted by $\mathcal K^+_i$.
\end{defin}

\begin{defin}[of $K^-_i$-equivalence] 
Let $p:STF\rightarrow PTF$ be the natural covering, 
$\pi_1^+(PTF)=p_*(\pi_1(STF))$, and let $\pi_1^-(PTF)$ 
be a set $\pi_1(PTF)\setminus \pi_1^+(PTF)$. 
Let $R_{K^-}$ be the set of all
$(\delta_1, \delta_2, i)\in \pi_1^-(PTF)\oplus\pi_1^-(PTF)\oplus\Z$
such that $i$ is even provided that $\pr^2(\delta_1\delta_2)$
is an orientation-preserving loop in $F$ and odd otherwise.

Let $L\in K^-$ be a front and $a\in S^1$  
one of the two points projecting to the self-tangency point. 
Deform $l$ keeping the two preimages of the tangency point 
opposite to each other in the fiber of $STF$, so
that $l(a)$ is located at $d\in STF$. The point $p(d)$ separates
$p(l)$ into two closed loops. Thus we obtain an ordered set of two elements 
$\delta_1, \delta_2\in\pi_1(PTF)$ corresponding to the deformed $l$. (The first element
is the projection of the arc of $l$ that has $a$ as its starting point.)
Clearly both elements belong to $\pi_1^-(PTF)$. 
We also correspond to the
front its Maslov index $\mu(L)$, which is even if and only if
$l=\delta_1\delta_2$ projects to an orientation-preserving loop in $F$.
Thus we obtain an element of $R_{K^-}$ corresponding to the deformed $l$.

There is a unique $K_i^-$-equivalence class of elements of $R_{K^-}$
corresponding to the undeformed $l$, where two elements of $R_{K^-}$ are
$K_i^-$-equivalent if one can be transformed to the other by an
action of the following group (which acts trivially
on the last summand in $R_{K^-}$):

The group is the index two subgroup of $\pi_1(PTF)\oplus\Z_2$ which is 
$\bigl(\pi_1^+(PTF)\oplus 0\bigr)\cup \bigl(\pi_1^-(PTF)\oplus 1\bigr)$. 
The action of the first summand of the group is via conjugation of the first
two summands in $R_{K^-}$ and the action of the second summand of the group
is via
permutation of the first two summands in $R_{K^-}$. (The factorization by
the action of this group corresponds to the ambiguity in the choices of one
of the two points of $S^1$ that project to the selftangency point and in 
the deformation of $l$ so that the chosen point is located at $d\in STF$.)

The set of all $K_i^-$-equivalence classes of elements of $R_{K^-}$ is
denoted by $\mathcal K^-_i$.

\end{defin}

\begin{defin}[of $T_i$-equivalence] 
Let $R_T$ be the set consisting of all quadruples
$(\delta_1, \delta_2, \delta_3, i)\in \pi_1(STF)\oplus\pi_1(STF)
\oplus\pi_1(STF)\oplus\Z$ such that $i$ is even provided that 
$\pr^2(\delta_1\delta_2\delta_3)$ is an orientation-preserving loop in 
$F$ and odd otherwise.

Let $L\in T$ be a front. 
Deform the lifting $l$ of $L$ in the neighborhood of the fiber of $\pr^2$ 
over the triple point, so that it maps the three preimages of the
triple point to one point in 
$STF$. (This triple point in $STF$ separates $l$ into three cyclicly
ordered oriented loops based at this point.) 
Then deform the singular knot $l$ (preserving the triple point), so that the
triple point is located at $d\in STF$. 
Choosing which one of the
three closed arcs of $l$ is first we obtain an ordered set of three elements
of $\pi_1(STF)$. We also correspond to the front its Maslov index
$\mu(L)\in \Z$, which is even if and only if $l=\delta_1\delta_2\delta_3$ 
projects to an orientation preserving loop in $F$.
Thus we obtain an element of $R_T$ corresponding to the deformed $l$.
There is a unique $T_i$-equivalence class of elements of $R_T$ corresponding
to the undeformed $l$, where two elements of $R_T$ are $T_i$-equivalent 
if one of them can be transformed to the other by the consequent action 
of the following groups (which all act trivially on the last
summand in $R_T$):

1. $\Z ^3$ whose element $(i_1, i_2, i_3)$ acts on 
$(\delta_1, \delta_2, \delta_3, i)\in R_T$ by mapping it to
$(f^{i_1}_2\delta_1 f^{-i_2}_2, f^{i_2}_2\delta_2 f^{-i_3}_2,
f^{i_3}_2\delta_3f^{-i_1}_2, i)$.
(This corresponds to the ambiguity in deforming $l$ 
to a singular knot with a triple point.)

2. $\pi_1(STF)$ whose elements act via conjugation of the first three
summands. 
(This corresponds 
to the ambiguity in deforming a singular knot with a triple point, so that
the triple point is at $d$.)

3 $\Z_3$ which acts via the cyclic permutation of the first three summands. 
(This corresponds to the ambiguity in the choice of the first closed arc of 
$l$.)  

The set of all $T_i$-equivalence classes of elements of $R_T$ is denoted by 
$\mathcal T_i$.
\end{defin}

\begin{defin}[of $\Pi_i$-equivalence]\label{Pii}
Let $R_{\Pi}$ be the set consisting of all quadruples
$(\delta_1, \delta_2, j, i)
\in \pi_1(STF)\oplus\pi_1(STF)\oplus\Z_2\oplus \Z$ such that 
$i$ is even provided that $\pr^2(\delta_1\delta_2)$ 
is an orientation-preserving loop in $F$ and odd otherwise.

Fix a local orientation at $\pr^2(d)\in F$.
Let $L\in \Pi$ be a front.
The direction of rotation of the coorienting normal at the cusp point $x$ 
under
traversing of a small neighborhood (in $L$) of it along the orientation of
$L$ defines a local orientation at $x\in F$. 
Deform the lifting $l$ of $L$, so that $l$ maps the two preimages of the
cusp crossing point to one point in $STF$.  
The double point of $l$ separates it into two 
oriented loops $\delta_1, \delta_2$ based at this point. 
The set of these two loops is
ordered by taking the loop corresponding to the branch of the cusp 
going away from it as the first one. Deform $l$ in $STF$ (preserving the
double point), so that the
double point is located at $d\in STF$. Transfer the local
orientation at the original location of the cusp along
the projection of the path of the double point of $l$ under the 
deformation to get the local orientation at $\pr^2(d)$. 
It is an element of $\Z_2$. 
We also correspond to $L$ its Maslov
index $\mu(L)\in \Z$. Clearly $\mu(L)$ is even if and only if
$\pr^2(l)=\pr^2(\delta_1\delta_2)$ is an orientation-preserving loop in $F$. 

Hence we obtain an element of $R_{\Pi}$ corresponding to the deformed $l$. One
verifies that there is a unique $\Pi_i$-equivalence class of elements of
$R_{\Pi}$
associated with the undeformed $l$, where two elements of $R_{\Pi}$ are said to be 
$\Pi_i$-equivalent if one of them can be transformed to the other by a
consequent action of the following groups (which all act trivially on the
last summand in $R_{\Pi}$):

1. $\Z ^2$ whose element $(i_1, i_2)$ acts on 
$(\delta_1, \delta_2, j, i)\in R _{\Pi}$
by mapping it to   
$(f_2^{i_1}\delta_1 f_2^{-i_2},$ $f_2^{i_2}\delta_2 f_2^{-i_1}, j,i)$. 
(Factorization by this action corresponds to the ambiguity in deforming $l$ 
to a singular knot with a double point.)

2. $\pi_1(STF)$ whose element $\alpha$ acts on the first two summands by
conjugation and acts trivially on the $\Z_2$-summand 
if and only if $\pr^2(\alpha)$ is an orientation-preserving loop in
$F$. (Factorization by this action corresponds 
to the ambiguity in deforming a singular knot with a double point, so that
the double point is located at $d$.)

3. $\Z_2$ whose elements act nontrivially only  
if $\pr^2(\delta_1)$ or $\pr^2(\delta_2)$ is $1\in\pi_1(F)$, 
in which case the action is by permutation of
the first two summands. (The reason for this factorization is that we want
to be able to identify the sets of 
$\Pi_i$- and $\overline \Pi$-equivalence classes, and one
can verify that under the passage
through a point of degree $5/4$ along the $\Pi$-stratum 
the order of the two loops is changed, see Figure~\ref{bifur5.fig}.)

The set of all $\Pi_i$-equivalence classes of
elements of $R_{\Pi}$ is denoted by $\mathcal P_i$.
\end{defin}

One verifies that: 

1. if $L_1, L_2\in K^+$ are $\overline {K^+}$-equivalent, then
they are $K^+_i$-equivalent; 

2. if $L_1, L_2\in K^-$ are $\overline {K^-}$-equivalent, then
they are $K^-_i$-equivalent; 

3. if  $L_1, L_2\in T$ are $\overline T$-equivalent, then
they are $T_i$-equivalent;

4. if $L_1, L_2\in \Pi$ are $\overline \Pi$-equivalent, then they are 
$\Pi_i$-equivalent.

And we get mappings 
$\psi ^+:\overline {\mathcal K^+}\rightarrow {\mathcal K^+_i}$, $\psi ^-:\overline
{\mathcal K^-}\rightarrow {\mathcal K^-_i}$, $\psi :\overline {\mathcal
T}\rightarrow {\mathcal T_i}$, $\psi ^{\pi}:\overline  
{\mathcal P}\rightarrow {\mathcal P _i}$.

\begin{thm}\label{interpretation}
The mappings $\psi^+$, $\psi^-$, $\psi$, and $\psi ^\pi$ are bijective.
\end{thm}

For the proof of Theorem~\ref{interpretation} see
Section~\ref{pfinterpretation}.

The connected components of the normalization of the 
dangerous self-tangency part of the discriminant in   
$\mathcal C$ are in a natural one-to-one 
correspondence with the set of $K^+$-equivalence classes of fronts from
$\mathcal C$. 
Analogous statement is true for the safe self-tangency and triple point 
parts of the discriminant, 
provided that $F$ is orientable.

Let $\mathcal C$ be a connected component of $\mathcal L$.
Let $\overline{\mathcal T_{\mathcal C}}$,  $\overline{\mathcal K^+_{\mathcal C}}$,
$\overline{\mathcal K^-_{\mathcal C}}$, $\mathcal T_{\mathcal C}$, $\mathcal
K^+_{\mathcal
C}$, and $\mathcal K^-_{\mathcal C}$ be the sets of equivalence classes 
corresponding to fronts from $\mathcal C$. 
Similarly to the 
case above we get the mappings $\phi_{\mathcal C}:\overline{\mathcal T_{\mathcal
C}}\rightarrow \mathcal T_{\mathcal C}$, $\phi^+_{\mathcal C}:
\overline{\mathcal K^+_{\mathcal  
C}}\rightarrow \mathcal K^+_{\mathcal C}$, and 
$\phi^-_{\mathcal C}:\overline{\mathcal K^-_{\mathcal  
C}}\rightarrow \mathcal K^-_{\mathcal C}$. 
\begin{prop}\label{orient}
The mapping $\phi^+_{\mathcal C}$ is bijective. The mappings
$\phi^-_{\mathcal C}$ and  
$\phi_{\mathcal C}$ are bijective provided that $F$ is orientable.
\end{prop}

For the proof of Proposition~\ref{orient} see
Section~\ref{pforient}.

One can construct examples showing that for nonorientable surfaces
the mappings $\phi ^-_{\mathcal C}$ and $\phi _{\mathcal C}$ fail to be
injective.

\subsection{An important generalization of
Theorem~\ref{integrability}}\label{importantgeneralization}
We say that an invariant of generic fronts is of $J^+$-type if it changes 
under crossings of the $K^+$-stratum and does not change
under crossings of the other codimension one strata of the discriminant.
Similarly to~\ref{generaltheorem} 
for such an invariant with values in an Abelian group
$G$ we define its $n$-th derivative, which assigns an element of $G$ 
to a front whose only nongeneric singularities are 
$n$ points of dangerous order one self-tangency. 
The invariant is said to be of order 
$k$ if $k$ is the minimal number such that the 
$k+1$-st derivative of the invariant is identically zero.

The invariants whose change under the crossing of a 
part of the $K^+$-stratum depends only on the $\overline
{K^+}$-equivalence class corresponding to the part are exactly the 
$J^+$-type invariants of order one. Proposition~\ref{orient} implies that
every such invariant of fronts from $\mathcal C$ can be obtained as a
$\overline {J^+}$ invariant for some $\psi:\mathcal K^+\rightarrow \Z$.

Such interpretation does not hold if we do not restrict ourselves 
to one component $\mathcal C$ of $\mathcal L$
and consider invariants defined in all components.
The reason is that nongeneric fronts belonging to 
different components of $\mathcal L$, and hence realizing different
$\overline {K^+}$-equivalence classes, can realize the same $K^+$-equivalence
class. 
(This happens only if the Maslov indices of the two fronts are different.) 
Theorem~\ref{interpretation} says that there is a natural bijection
between the sets $\mathcal K^+_i$ and $\overline {\mathcal K^+}$. 
Hence every order one $J^+$-type invariant of fronts on $F$
can be obtained by integration in all components of $\mathcal
L$ of some $\psi:\mathcal K^+_i\rightarrow \Z$. 
One can easily verify that statements of Theorem~\ref{integrability} 
hold if one substitutes 
$K^+$- by $K^+_i$-equivalence classes in the
formulation of the Theorem and in the axiomatic description of the $\overline
{J^+}$ invariant. 
(The proof of this version of the theorem is the same as the
original one.)

The interpretation of the $\overline {J^-}$ invariant is similar to the one
of $\overline {J^+}$, with the
difference that in the case of nonorientable $F$ not every order one 
$J^-$-type invariant of fronts in $\mathcal C$ 
is an integral of some $\psi:\mathcal
K^-\rightarrow \Z$.
(For nonorientable $F$ two fronts from $\mathcal C$
realizing the same $K^-$-equivalence class can realize different 
$\overline {K^-}$-classes.) 
Similarly to the case of $J^+$, we see that every order one
$J^-$-type invariant of fronts in
$\mathcal L$ can be obtained by integration in all components of $\mathcal
L$ of some $\psi:\mathcal K^-_i\rightarrow \Z$. 
One verifies that statements of
Theorem~\ref{integrability} 
hold if one substitutes
$K^-$- by $K^-_i$-equivalence classes in the
formulation of the Theorem and in the axiomatic description of the $\overline
{J^-}$ invariant.

The operation of changing $L\in \Pi$ to $L'\in T$
shown in Figure~\ref{subst1.fig} induces a decomposition of the 
$\Pi$-stratum into parts corresponding to different $T_i$-equivalence
classes, and one obtains the corresponding version 
of Theorem~\ref{integrability}. 

One verifies that this operation induces a mapping $g:\mathcal 
P_i\rightarrow \mathcal T_i$
that sends a $\Pi_i$-equivalence class of $(\delta_1, \delta_2, i, j)\in
R_{\Pi}$ to a $T_i$-equivalence class of $(\delta_1, \delta_2, 1, j)\in
R_{T}$. (Here $1\in\pi_1(STF)$ is a class of a trivial loop.) One verifies
that the part of the $\Pi$-stratum corresponding to $\pi_i\in \mathcal
P_i$ is adjacent (along the $\Pi\Lambda$-stratum, see~\ref{strata})
to the part of the $T$-stratum corresponding to
$g(\pi_i)\in \mathcal T_i$. This means
(see~\ref{halfchange}) that if the magnitudes of the 
change of a ${\St'}$-type invariant 
under the crossings of $T$- and
$\Pi$-strata depend only on the component of the normalization
corresponding to the crossings, then the change under the positive crossing 
of 
the $\pi_i$-part of the $\Pi$-stratum should be half of the change
under the positive crossing of the $g(\pi_i)$-part of the $T$-stratum. 
Every ${\St'}$-type invariant of this sort 
is a $\overline {\St'}$ invariant for
some $\psi:\mathcal T_i\rightarrow \Z$.

An important observation is that Remarks~\ref{importrem} hold for the
versions of the three invariants described above. 

\begin{rem}
A front with a dangerous self-tangency point lifts to a singular 
Legendrian knot in $ST^*F$. (It has a double point.) One can verify 
that order one $J^+$-type invariants of fronts are exactly order one invariants 
of Legendrian knots in $ST^*F$.

\begin{figure}[htbp]
 \begin{center}
  \epsfxsize 11cm
  \hepsffile{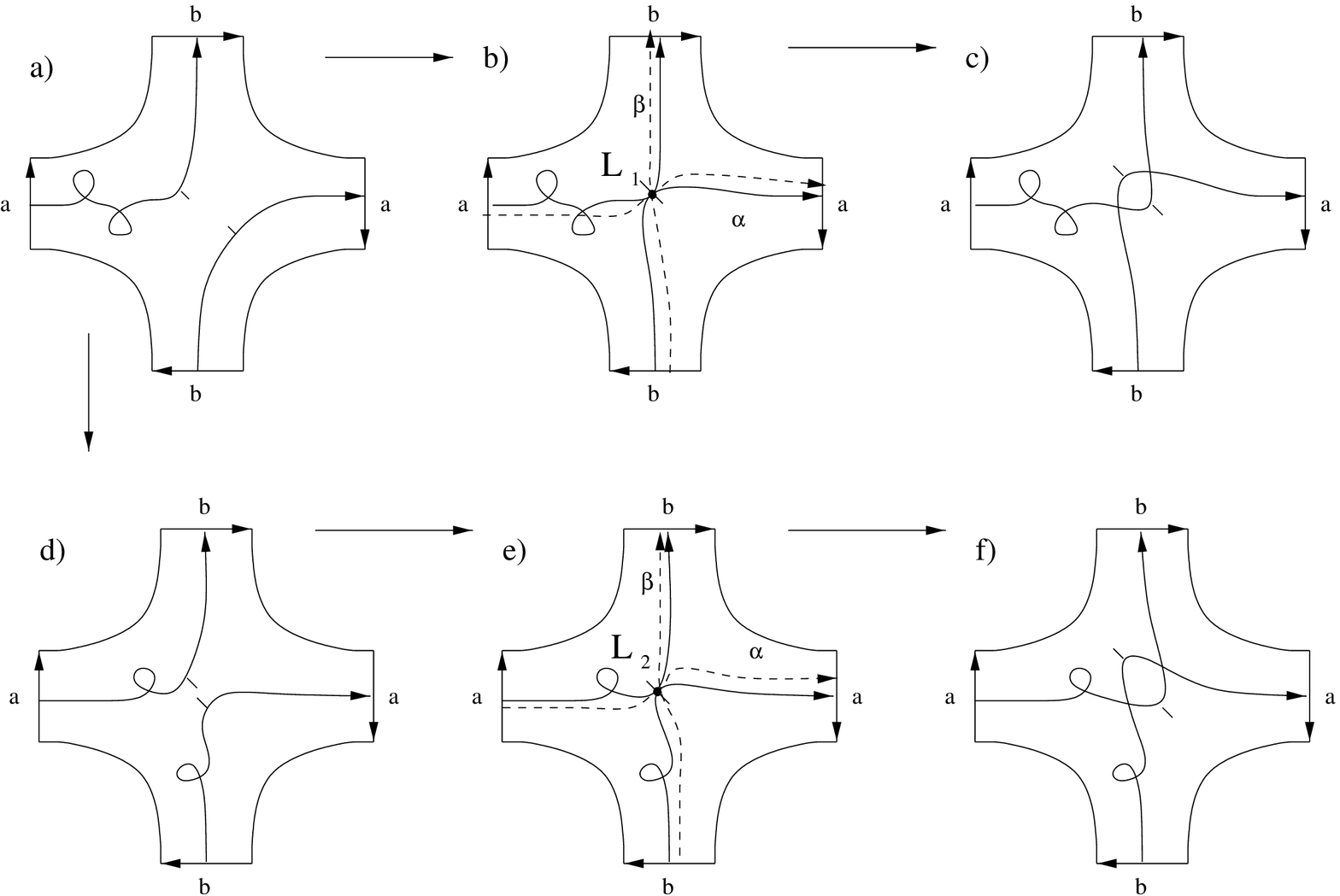}
 \end{center}
\caption{}\label{example1.fig}
\end{figure}

Similarly a front with a safe self-tangency point lifts to a
singular Legendrian knot in $PT^*F$. However for nonorientable $F$ it is 
not true that order one $J^-$-type invariants are order one invariants of 
Legendrian knots in $PT^*F$. To see this consider a nonorientable 
surface shown in Figure~\ref{example1.fig}. 
One can verify that any order one invariant of 
Legendrian knots in $PT^*F$ takes equal values on the Legendrian lifting 
of fronts shown in Figure~\ref{example1.fig}~c and~f, but there
exists an order one $J^-$-type invariant which takes different values on the
two fronts.
To construct such $J^-$-type order one invariant we observe that 
any function on $\mathcal K^-_i$ is integrable along the loop $\gamma_1$ and 
hence gives rise to an order one $J^-$-type invariant, see
Subsection~\ref{importantgeneralization}. Both fronts are
obtained from the front in Figure~\ref{example1.fig} by an isotopy and a
(positive) crossing of the $K^-$-stratum. Finally we observe that the fronts 
$L_1$ and $L_2$ in Figure~\ref{example1.fig}~b and~e 
realise different $K^-_i$-equivalence classes. If one of them is 
$(\alpha, \beta, 0)\in R_{K^-}$, then the second one is $(\alpha f^2, \beta
f^2, 0)\in R_{K^-}$. (Here $f\in\pi_1^-(PTF)$ is the class of the fiber of
$PT^*F\rightarrow F$). The fact that the two classes are different can be
easily obtained from Proposition~\ref{Preissman} and the identities similar
to 2 and 3 of Proposition~\ref{commute}. (Recall that $(f,0)$ does not
belong to the group by the action of which we quotient $R_{K^-}$ to define 
the $K^-_i$-equivalence relation.) Hence any function on $\mathcal K^-_i$
that takes different values on the $K^-_i$-equivalence classes of $L_1$ and
$L_2$ gives rise to the desired invariant. 

One can show that the Legendrian liftings of $L_1$ and $L_2$ can be
transformed to each other in the class of Legendrian knots with a double point 
(so that the points on the parameterizing circle corresponding to the double
point change continuously under the transformation). Hence the values of the
derivative of any order one invariant of Legendrian knots in $PT^*F$ on
$L_1$ and on $L_2$ are equal, and thus the values of the invariant on the 
Legendrian liftings to $PT^*F$ 
of fronts in Figure~\ref{example1.fig}~c and~f are equal.

\end{rem}

\section{An explicit formula for the 
finest order one $J^+$-type invariant on orientable $F\neq S^2$}

Below we give an explicit formula for $I^+$ an order one 
$J^+$-type invariant  of generic 
wave fronts on an orientable surface $F$.
If $F\neq S^2$ then the invariant 
distinguishes every two generic wave fronts that one can 
distinguish using order one $J^+$-type invariants with values in any Abelian
group (not necessarily torsion free), see Theorem~\ref{universal}.

We assign a positive (resp. negative) sign to a cusp point
if the coorienting vector turns in the positive (resp. negative)
direction while traversing a small neighborhood of the cusp point along the
orientation of the front. We denote half of the number of
positive and negative cusp of the front $L$ by
$C^+$ and $C^-$ respectively.

\begin{figure}[htbp]
 \begin{center}
  \epsfxsize 8cm
  \hepsffile{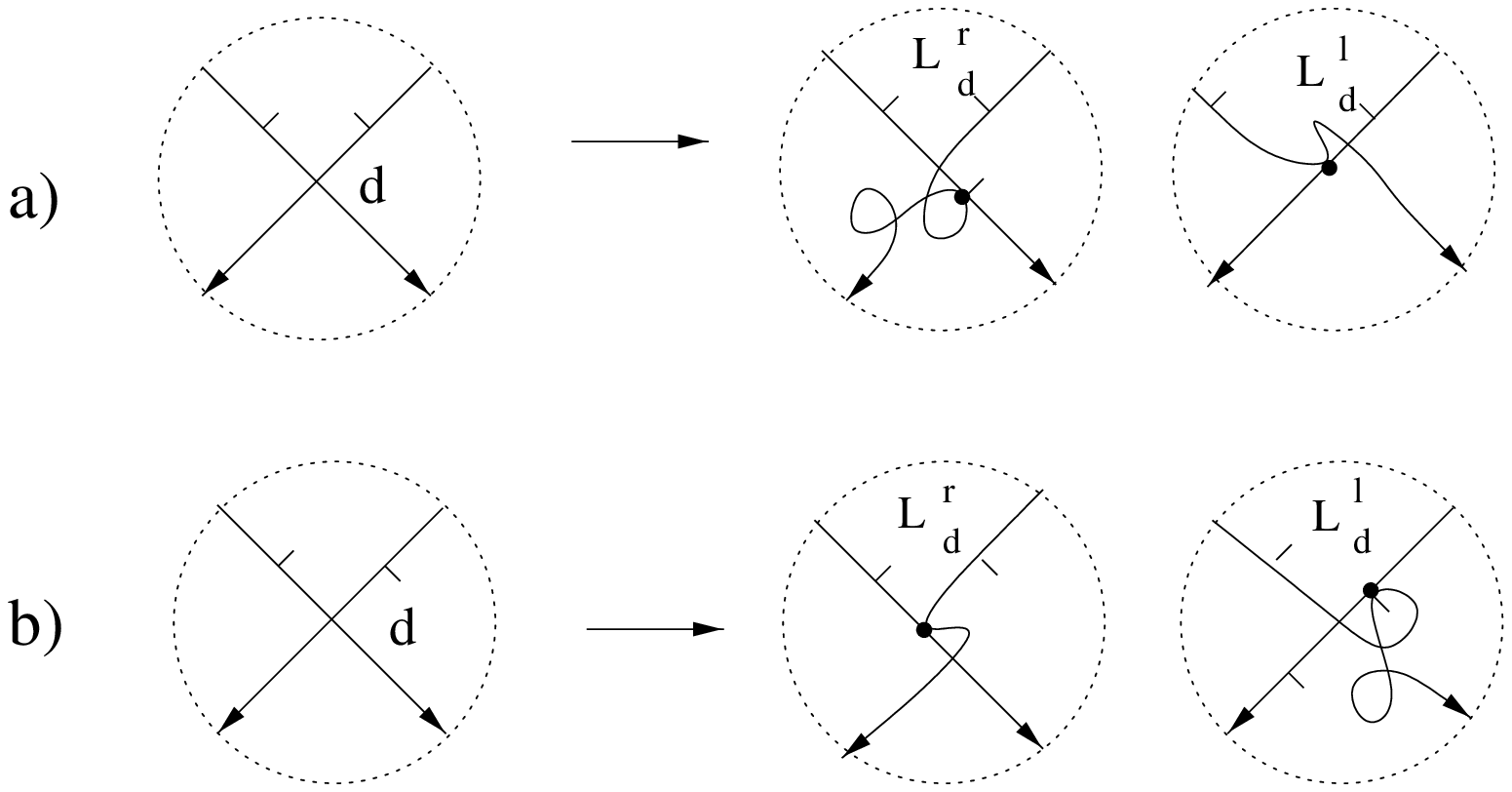}
 \end{center}
\caption{}\label{refine+.fig}
\end{figure}

Using the orientation of $F$ one gets that there are four types of double
points of a wave front. Two of them are shown in Figure~\ref{refine+.fig}, 
two more are obtained by a change of coorientation on both participating
branches. To a double point $d$ of $L$ we correspond two nongeneric 
fronts $L^r_d, L^l_d\in K^+$, as it is shown in Figure~\ref{refine+.fig}.
(The $L^r_d, L^l_d$ fronts for the double points of the types not shown in the
Figure are obtained by a change of coorientation.) We denote by 
$[L^r_d], [L^l_d]$ the $K^+$-equivalence classes of these fronts. 
The set $\mathcal K^+$ is naturally identified with
the set $R^+$ which is the factor of 
$\pi_1(STF)\oplus\pi_1(STF)$ modulo the action of $\pi_1(STF)$ by
conjugation of the first two summands and by the action of $\Z_2$ permuting
the summands.

We denote by $\Z[\mathcal K^+]$ the free $\Z$-module
of all formal finite integer combinations of elements of $\mathcal K^+$. 
For a wave front $L$ we denote by 
$[lf_2^{-1}, f_2]$, $[l, 1]$, $[lf_2,f_2^{-1}]$ the $K^+$-equivalence
classes described by the corresponding elements of $R^+$.

\begin{thm}\label{niceformula}
Put $I^+(L)\in \Z[\mathcal K^+]$ to be
\begin{equation}
\bigl(\sum_d([L^r_d]-[L^l_d])\bigr)-C^-([l,1]-[lf_2,
f_2^{-1}])-C^+([l,1]-[lf_2^{-1}, f_2]).
\end{equation}
Then $I^+(L)$ is an order one $J^+$-type invariant. Under the positive
crossing of the $K^+$-stratum
corresponding to $[s_1, s_2]\in R^+$
it
increases by 
$\bigl(2[s_1, s_2]-[s_1f_2, s_2f_2^{-1}]-[s_1f_2^{-1}, s_2f_2]\bigr)$.
\end{thm}

The proof of the Theorem is straightforward. One verifies that $I^+(L)$ does
not change under crossings of $\Lambda$-, $K^-$-, $\Pi$-, and $T$-strata.
To see that $I^+(L)$ has the described behaviour under the
crossings of the $K^+$-stratum, one uses the fact that $f_2$ is in the center
of $\pi_1(STF)$, see~\ref{commute}. 

(To construct a similar invariant of wave fronts on non-orientable surfaces
one symmetrizes the construction of $L^r_d$ and $L^l_d$ and gets four
elements of $\mathcal K^+$ corresponding to the double point $d$.)  

The following theorem says that for orientable $F\neq S^2$ the invariant 
$I^+$ is the finest order one $J^+$-type invariant.

\begin{thm}\label{universal}
Let $\mathcal C$ be a connected component of the space of fronts on an
orientable $F\neq S^2$. 
Let $L_1, L_2\in \mathcal C$ be generic fronts, and
$\overline{I^+}$ an order one $J^+$-type invariant of fronts with
values in some Abelian group $G$ (not necessarily torsion free).
Then $\overline{I^+}(L_1)=\overline{I^+}(L_2)$, provided that
$I^+(L_1)=I^+(L_2)$.
\end{thm}

For the proof of Theorem~\ref{universal} see Section~\ref{pfuniversal}.

It is well known that in the case of $F=\R^2$ the partial linking polynomial 
of F.~Aicardi~\cite{Aicardi} appears to be the finest order one invariant 
in the above sense.
Other order one $J^+$-type invariants of wave fronts on surfaces 
were constructed by M.~Polyak~\cite{Polyakbennequin} 
and by the author in~\cite{Tchernovamsvolume}. 
They are significantly easier for calculation but are not the finest in the
above sense.

\section{Homotopy groups of the space of fronts}
Fix $a\in S^1$, then $L\in \mathcal L$ represent an element of 
$\pi_1(F, L(a))$, the lifting $l$ of $L$ to
$STF$ represents an element of $\pi_1(STF, l(a))$, and the lifting $\vec l$
of $L$ to $CSTF$ represents an element of $\pi_1(CSTF, \vec l(a))$.

\subsection{Fundamental group of the space of fronts on an orientable
surface.}

For orientable surfaces the group $\pi_1(\mathcal L,
L)$ appears to be much simpler than for nonorientable surfaces.

\begin{thm}\label{pi1S2}
Let $F=S^2$ and $L$ a front on $S^2$. Then $\pi_1(\mathcal L,
L)=\Z \oplus \Z_2$.
\end{thm}

\begin{thm}\label{pi1T2}
Let $F=T^2$ (torus) and $L$ a front on $T^2$. Then $\pi_1(\mathcal L,
L)=\Z^4$.
\end{thm}

For the proofs of Theorems~\ref{pi1S2} and~\ref{pi1T2} 
see Subsection~\ref{pf4surfaces}.

\begin{thm}\label{pi1orient}
Let $F\neq S^2, T^2$ be an orientable surface (not necessarily compact),
and let $L$ be a front on $F$.
\begin{description}
\item[\textrm{I}] If $L\neq 1\in \pi_1(F)$, then
$\pi_1(\mathcal L, L)=\Z^3$.
\item[\textrm{II}] If $L=1\in\pi_1(F)$, then $\pi_1(\mathcal
L, L)=\Z\oplus\pi_1(STF)$.
\end{description}
\end{thm}

For the proof of Theorem~\ref{pi1orient} see
Subsection~\ref{pfpi1arbitrary}.

\subsection{Fundamental group of the space of fronts on a nonorientable
surface}

We denote by $\pi_1^{pres}(F)$ the subgroup of $\pi_1(F)$ consisting of all
orientation-preserving loops and by $\pi_1^{rev}(F)$ the subset of
$\pi_1(F)$ which is 
$\pi_1(F)\setminus\pi_1^{pres}(F)$. We denote by $\pi_1^{pres}(STF)$ the
subgroup of $\pi_1(STF)$ which is a preimage of $\pi_1^{pres}(F)$ under 
$\pr^2_*:\pi_1(STF)\rightarrow \pi_1(F)$ and by $\pi_1^{rev}(STF)$ the
subset of $\pi_1(STF)$ which is a preimage of $\pi_1^{rev}(F)$ under    
$\pr^2_*$.
We denote by $\Z^{ev}$  the subgroup of even numbers in $\Z$
and by $\Z^{odd}$ the subset of odd numbers in $\Z$.

\begin{thm}\label{pi1RP2}
Let $F=\R P^2$ and $L$ a front on $\R P^2$. Then $\pi_1(\mathcal L,
L)$ is isomorphic to $\Z\oplus\Z_2$.
\end{thm}

\begin{thm}\label{pi1K}
Let $F=K$ (Klein bottle), and let $L$ be a front on $K$.
\begin{description}
\item[\textrm{I}] If $L\in \pi_1^{pres}(K)$, then 

{\bf a:} $\pi_1(\mathcal L, L)$ is isomorphic to $\Z\oplus\pi_1(STK)$,
provided that $l=b^{2k}\in\pi_1(STK, l(a))$
for some $b\in\pi_1^{rev}(STK, l(a))$. 

{\bf b:} $\pi_1(\mathcal L, L)=\Z^4$
otherwise.
\item[\textrm{II}] If $L\in \pi_1^{rev}(K)$, then
$\pi_1(\mathcal L, L)$ is isomorphic to $\Z^2$.
\end{description}
\end{thm}

For the proofs of Theorems~\ref{pi1RP2} and~\ref{pi1K} 
see Subsection~\ref{pf4surfaces}.

Let $F\neq \R P^2, K$ be a nonorientable
surface (not necessarily compact), and
let $L$ be a front on $F$ such that $L\neq 1\in\pi_1(F,
L(a))$. 
One can show that there exists a unique maximal Abelian
subgroup  $G_{L}<\pi_1(F, L(a))$ containing
$L\in \pi_1(F, L(a))$, and that this $G_{L}$ is isomorphic to $\Z$
(see also Proposition~\ref{Preissman}).
Let $g$ be a generator of $G_L$ and $L_g$ a front 
such that $\vec l_g(a)=\vec l(a)$ and $L_g=g\in\pi_1(F, L(a))$.
One can show that
$l\in\pi_1(STF, l(a))$ can be presented in the unique way
as $l_g^k f_2^m\in\pi_1(STF, l(a))$ (see also the Proof of
Theorem~\ref{pi1nonorient}). 

\begin{thm}\label{pi1nonorient}
Let $F\neq \R P^2,K$
be a nonorientable surface (not necessarily compact), and
let $L$ be a front on $F$.

\begin{description}
\item[\textrm{I}] If $L\in \pi_1^{rev}(F)$,
then $\pi_1(\mathcal L,L)$ is isomorphic to $\Z^2$.

\item[\textrm{II}] If $L\in \pi_1^{pres}(F)$ and $L\neq 1\in \pi_1(F)$,
then:

{\bf a:}  $\pi_1(\mathcal L, L)$ is isomorphic to $\Z\oplus\pi_1(K)$,
provided that $L_g\in
\pi_1^{rev}(F)$ and that 
$l=l_g^{2k}\in\pi_1(STF, l(a))$, for some $k\in \Z$. 

{\bf b:} $\pi_1(\mathcal L, L)=\Z^3$
otherwise.

\item[\textrm{III}] If $L=1\in \pi_1(F)$, then:

{\bf a:} $\pi _1(\mathcal L, L)$ is isomorphic to 
$\Z\oplus \pi_1^{pres}(STF)$, provided that
$l\neq 1\in \pi_1(STF)$.

{\bf b:} $\pi_1(\mathcal L, L)$ is isomorphic to the subgroup of 
$\Z\oplus\pi_1(STF)$ which is $(\Z^{even}\oplus\pi_1^{pres}(STF))\cup
(\Z^{odd}\oplus\pi_1^{rev}(STF))$, provided that
$l=1\in \pi_1(STF)$. 
\end{description}
\end{thm}

For the proof of Theorem~\ref{pi1nonorient} see
Subsection~\ref{pfpi1arbitrary}.

Statement \textrm{III}.a (resp. \textrm{III}.b) corresponds to $L$ that is 
regular homotopic to one of the fronts of type $K_{i,k}$, $i>0$, 
(resp. of type $K_{0,k}$) see Figure~\ref{basic.fig}.

\subsection{Higher homotopy groups of the space of fronts.}
\begin{thm}\label{pin}
Let $F$ be a surface (not necessarily compact or orientable)
and let $L$ be a front on $F$.

\begin{description}
\item[\textrm{I}] If $F$ is $S^2$ or  $\R P^2$, 
then $\pi_2(\mathcal L, L)=\Z$,
and
$\pi_n(\mathcal L,
L)=\pi_n(S^2)\oplus\pi_{n+1}(S^2)$, $n\geq 3$.

\item[\textrm{II}] If $F\neq S^2, \R P^2$, then $\pi_n(\mathcal L, L)=0$, $n\geq 2$.
\end{description}
\end{thm}

For the Proof of Theorem~\ref{pin} see Section~\ref{pfpin}.

\section{Proof of Theorem~\ref{integrability}}\label{pfintegrability}
We prove only the statements of 
Theorem~\ref{integrability} related to the existence of the $\overline{\St '}$
invariant integrating $\psi_1$ in $\mathcal C$.
The proofs of statements related to the
existence of $\overline{J^+}$ and $\overline{J^-}$ invariants are 
obtained in a similar way.

In order for $\overline{\St'}$ to be well defined,
the change of it along any generic loop 
has to be zero. This proves the necessity of the conditions
described in Theorem~\ref{integrability}.
Let us prove that these conditions are sufficient for the existence of
$\overline{\St'}$ invariant integrating $\psi_1$ in $\mathcal C$.

Put $\overline{\St'}(L)$ to be any number.
Let $L'\in \mathcal C$ be a generic front and $p$ a generic path
connecting $L$ to $L'$. Similarly to the case of a closed loop we define
$\Delta_{\overline{\St'}}(p)$.
Put $\overline{\St'}(L')=\overline{\St'}(L)+\Delta_{\overline{\St'}}(p)$.
To prove the theorem it suffices to show that $\overline{\St'}(L')$ 
is independent of the generic path $p$ which was used to define it. 
The last statement follows from Lemmas~\ref{onlyhomotop} and~\ref{zerojump}.
Thus we proved Theorem~\ref{integrability} modulo these two lemmas.
\qed

\begin{lem}[Cf. Arnold~\cite{Arnoldsplit}]\label{onlyhomotop}
Let $p$ be a generic path in $\mathcal C$ connecting $L$ to itself.
Then $\Delta_{\overline{\St'}}(p)$ depends only on the element 
of $\pi_1(\mathcal C, L)$ realized by $p$.
\end{lem}

\begin{lem}\label{zerojump}
Let $F$ be a surface, $\mathcal C$ a connected component of $\mathcal L$,
and $\psi_1:\mathcal T\rightarrow \Z$ a function 
integrable along those of the 
loops $\gamma_1$ and $\gamma_2$ that participate in the 
statement of Theorem~\ref{integrability} corresponding to $F$ and $\mathcal C$. 
Then every $\alpha\in\pi_1(\mathcal C, L)$ can be realized by 
a generic loop $q_{\alpha}$
in $\mathcal C$ such that $\Delta_{\overline{\St'}}(q_{\alpha})=0$.
\end{lem}

\subsection{Proof of Lemma~\ref{onlyhomotop}}\label{pfonlyhomotop}
To prove the Lemma, 
it suffices to show that if we go around any codimension two 
stratum of the discriminant 
along a small generic loop $r$ (not necessarily starting at $L$), then 
$\Delta_{\overline{\St'}}(r)=0$. All the codimension two strata are described
in the following Theorem.

\begin{figure}[htbp]
 \begin{center}
  \epsfxsize 8cm
  \hepsffile{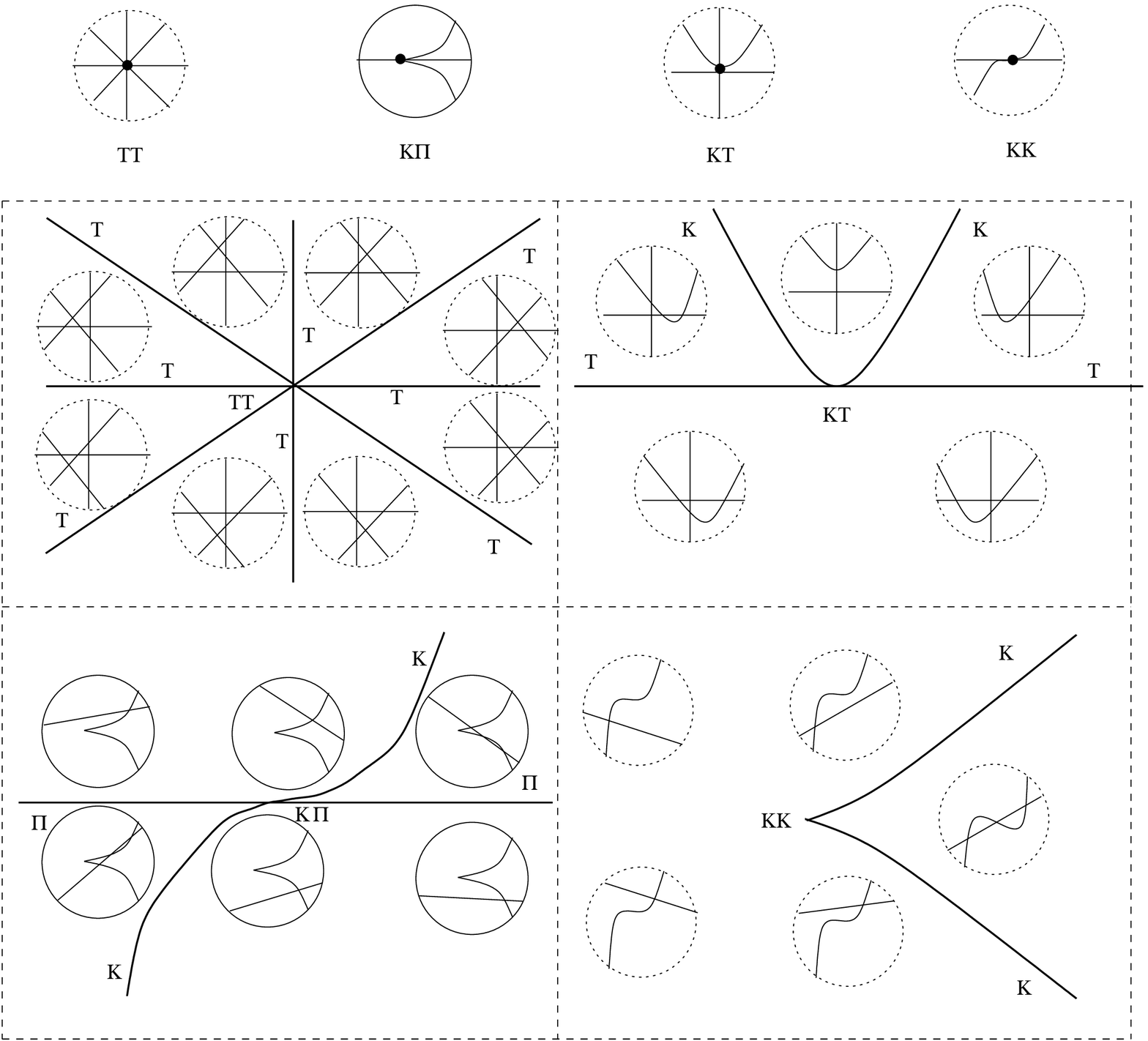}
 \end{center}
\caption{}\label{bifur4.fig}
\end{figure}

\begin{figure}[htbp]
 \begin{center}
  \epsfxsize 8cm
  \hepsffile{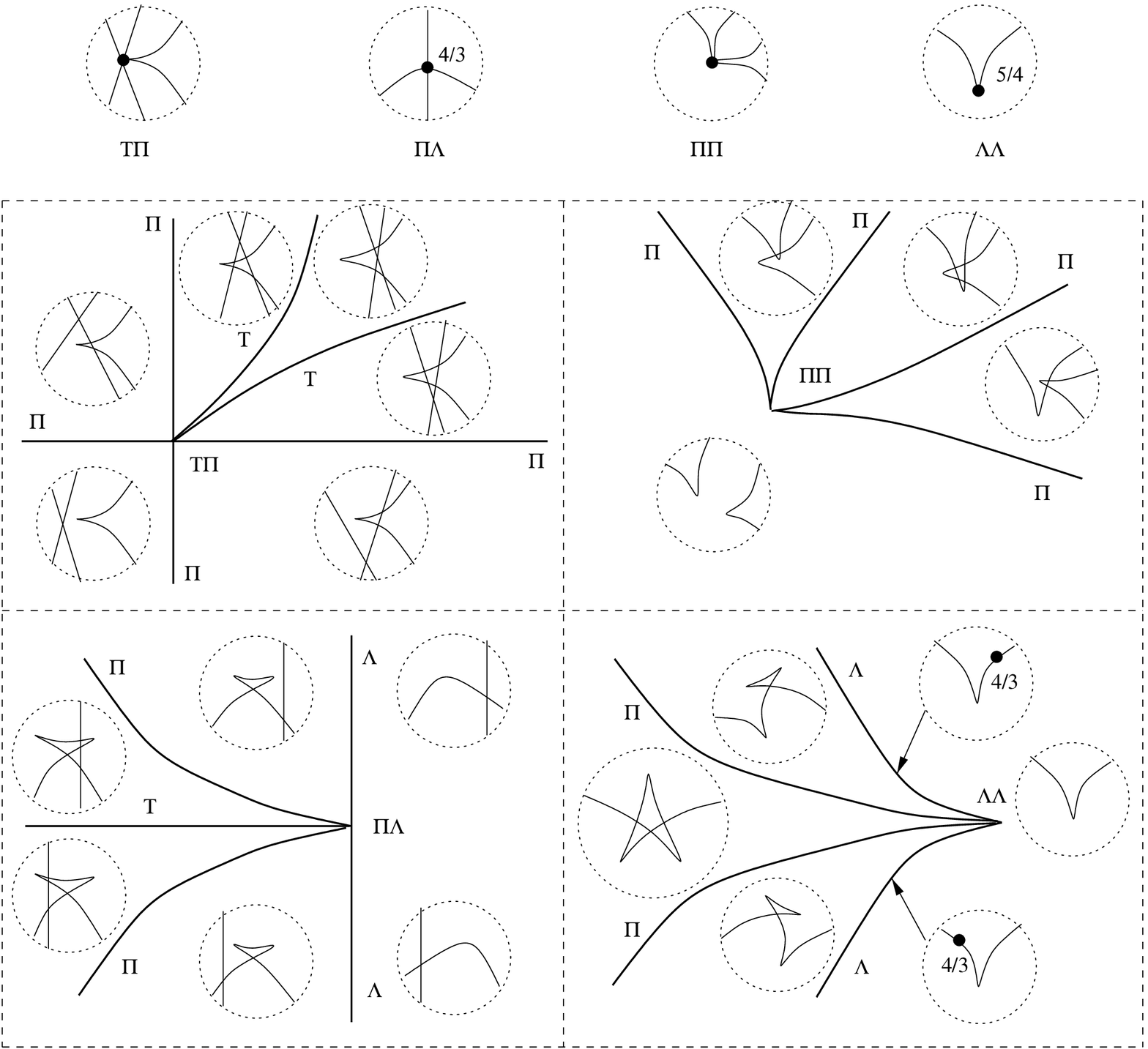}
 \end{center}
\caption{}\label{bifur5.fig}
\end{figure}

\begin{thm}[Arnold~\cite{Arnoldsplit}]\label{strata}
The strata of codimension two of the discriminant of $\mathcal L$ 
are formed
by fronts with two nongeneric singular points that are singular of codimension
one, and by fronts with one nongeneric singular point that is one of the
following 
(see Figure~\ref{bifur4.fig} and Figure~\ref{bifur5.fig}):

1) A quadruple point with pairwise transverse tangent lines.
This stratum is denoted by $TT$.

2) A cusp passing through a branch in such a way that they 
have the same tangent line. This stratum is denoted by $K\Pi$.

3) A degenerate triple point at which two branches are tangent of order one
and the third branch is transverse to them. This stratum is denoted by $KT$.

4) A point of a cubical self-tangency. This stratum is denoted by $KK$.

5) A cusp point passing simultaneously  through 
two branches. (Here it is assumed that the lines tangent to the three
participating branches are different.) This stratum is denoted by $T\Pi$.

6) A point of degree $\frac{4}{3}$ passing through a branch of the
front. (Here it is assumed that the two branches have transverse tangent
lines.)
This stratum is denoted by $\Pi\Lambda$.

7) Two coinciding cusp points. (Here it is assumed that the lines 
tangent to the two branches are different.)
This stratum is denoted by $\Pi \Pi$.

8) A point of degree $\frac{5}{4}$. This stratum is denoted by 
$\Lambda\Lambda$.

\end{thm}

\begin{emf}
The only strata of the discriminant of codimension two 
in whose bifurcation
diagram triple points or cusp crossings are present are:
two distinct codimension one singular points one of which is a triple
point; two distinct codimension one singular points one of which is a
cusp crossing point; 
and strata 
$TT$, $KT$, $K\Pi$, $T\Pi$, $\Pi \Pi$, $\Pi \Lambda$ and $\Lambda \Lambda$ 
(in the notation of Lemma~\ref{strata}).

If $r$ is a small loop going around a stratum of two distinct
codimension one singular points
one of which is a cusp or a triple point, then in 
$\Delta_{\overline{\St'}}(r)$ 
we have each of
the participating 
$T$-equivalence classes twice, once with the plus sign of the
newborn vanishing triangle, once with the minus. 
Hence $\Delta_{\overline{\St'}}(r)=0$.

To prove the statement for the other strata we use the bifurcation diagrams
shown in Figure~\ref{bifur4.fig} and Figure~\ref{bifur5.fig}.

Let $r$ be a small loop going around the $TT$-stratum. 
We can assume that it corresponds to a loop in 
Figure~\ref{vers1.fig} directed counterclockwise. 
There are eight terms in $\Delta_{\overline{\St'}}(r)$. We split them
into pairs I, II, III, IV, as it is shown in Figure~\ref{vers1.fig}. One can 
see that the wave fronts from the same pair are $T$-equivalent.
For each branch the sign of the colored triangle 
is equal to the sign of the triangle that died under the
$T$-stratum crossing shown on the next (in the counterclockwise direction)
branch. 
The sign of the dying vanishing triangle 
is minus the sign of the newborn vanishing triangle.
Finally, one can see that the signs
of the colored triangles inside each pair are opposite. Thus all these eight
terms cancel out, and $\Delta_{\overline{\St'}}(r)=0$.

\begin{figure}[htbp]
 \begin{center}
  \epsfxsize 6cm
  \hepsffile{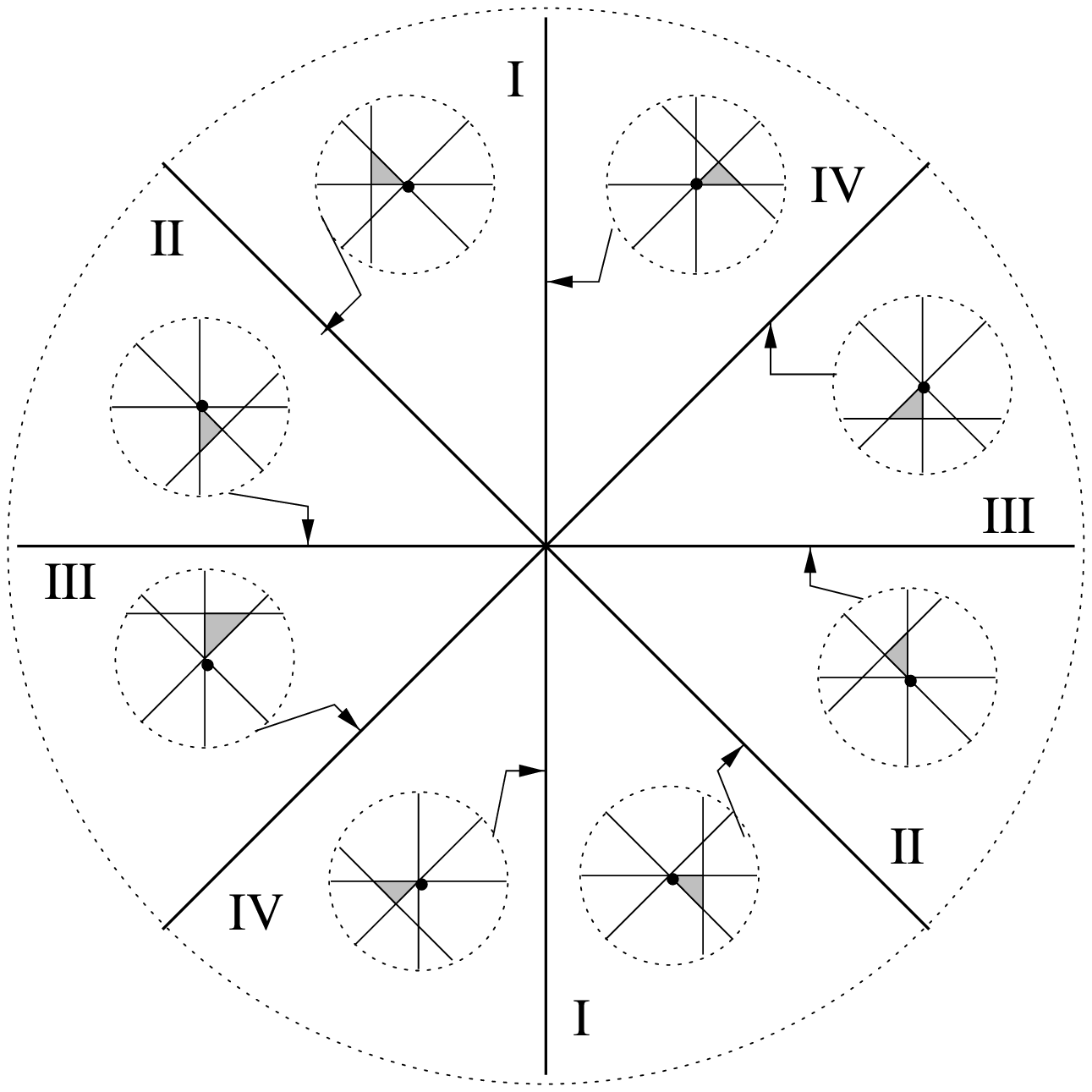}
 \end{center}
\caption{}\label{vers1.fig}
\end{figure}

There are six terms in $\Delta_{\overline{\St'}}(r)$ for $r$ a loop 
going around the $T\Pi$-stratum. Two of them correspond 
to $r$ crossing the 
$T$-stratum, and four to $r$ crossing the $\Pi$-stratum. 
One verifies that 
the two terms corresponding to $r$ crossing the  
$T$-stratum cancel
out. 
We split the other four terms into pairs such that the terms in
a pair correspond to the cusp crossing the same branch of the front. 
One verifies that the $T$-equivalence classes inside a pair coincide and the 
signs of the terms are opposite. Thus the terms inside each pair cancel 
out.

In a similar way one show that $\Delta_{\overline{\St'}}(r)=0$
for a small loop $r$ going around $K\Pi$-, $KT$-, $\Pi \Pi$-, and 
$\Lambda \Lambda$-strata.

Finally, let $r$ be a loop going around the $\Pi \Lambda$-stratum. There are
three terms in $\Delta_{\overline{\St'}}(r)$. One of them corresponds to $r$
crossing the $T$-stratum and two to $r$ crossing
the $\Pi$-stratum. One verifies that 
that the sign of the 
crossing of the $T$-stratum
is opposite from the signs of the crossings of the 
$\Pi$-stratum, and that all three $T$-equivalence classes are the same. 
We denote the class by $t$. 
Thus,
$\Delta_{\overline{\St'}}(r)=\frac{1}{2}\psi(t)+\frac{1}{2}\psi(t)-\psi(t)=0$.

\begin{rem}\label{halfchange}
Recall that the magnitude of
the change of $\St'$ under the crossing of the part of the 
$\Pi$-stratum corresponding to the $T$-equivalence class $t$ was put to be 
$\frac{1}{2}\psi(t)$. One can see 
that if we substitute $\frac{1}{2}$ by another
constant, then there is no hope for constructing
an invariant of this sort, unless $\psi$ is zero on all the 
$T$-equivalence classes appearing on the $\Pi$-stratum.
\end{rem}

This finishes the Proof of Lemma~\ref{onlyhomotop}.
\qed
\end{emf}

\subsection{Constructions and facts needed for the proof of
Lemma~\ref{zerojump}}\label{facts}

\begin{emf}\label{H-principle}{\em Parametric $h$-principle.\/}
The parametric $h$-principle proved
for the Legendrian curves by M.~Gromov~\cite{Gromov} says that 
the space of wave fronts $\mathcal L$ is weak homotopy equivalent to the
space $\Omega CSTF$ of all free loops in $CSTF$, see~\ref{basicdefin}.
The mapping $h:\mathcal L\rightarrow \Omega CSTF$ that gives the equivalence 
sends a wave front $L$ (corresponding to the Legendrian curve $l$) 
to the loop $\vec l\in \Omega CSTF$.

Fix $a\in S^1$. 
Let $q$ be a loop in $\mathcal L$ starting at $L$. At any moment of time $q(t)$
is a wave front that can be lifted to a loop in $CSTF$. Thus $q$ gives
rise to the mapping $q_h:S^1\times S^1\rightarrow CSTF$. 
(In the product $S^1\times S^1$ the first copy of $S^1$ corresponds
to the parameterization of a front and the second to the parameterization of 
the loop $q$.) The mapping $q_h$ restricted to $a\times S^1$ is a loop 
$t_a(q)$ in $CSTF$. 
One can verify that the mapping $t_a:\pi_1(\mathcal L, L)\rightarrow \pi_1(CSTF,
\vec l(a))$ is a homomorphism. 
\end{emf}

\begin{prop}\label{centralizer} 
$t_a:\pi_1(\mathcal L, L)\rightarrow
\pi_1(CSTF, \vec l(a))$ is an isomorphism of $\pi_1(\mathcal L, L)$
onto  the centralizer $Z(\vec l)$ of $\vec l \in \pi_1(CSTF, \vec l(a))$.
\end{prop}
\begin{emf}{\em Proof of Proposition~\ref{centralizer}.\/}
Let $p:\Omega CSTF\rightarrow
CSTF$ be the mapping that sends $\omega\in\Omega CSTF$ to $\omega(a)\in
CSTF$. (One can verify that this $p$ is a Serre fibration, with 
the fiber of it over a point isomorphic to the space of loops 
based at the point.) The $h$-principle
(see~\ref{H-principle}) implies
that to prove the Proposition it suffices to show that 
$p_*:\pi_1(\Omega CSTF, \vec l)\rightarrow \pi_1(CSTF, \vec l(a))$
is an isomorphism of $\pi_1(\Omega CSTF, \vec l)$ onto  $Z(\vec l)$.

A Proposition proved by V.L.~Hansen~\cite{Hansen} says that if $X$ is a
topological space with $\pi_2(X)=0$, then 
$\pi_1(\Omega X, \omega)=Z(\omega)<\pi_1(X, \omega (a))$. (Here $\Omega X$
is the space of free loops in $X$ and $\omega$ is an element of $\Omega X$.)
One can verify that $\pi_2(CSTF)=0$ for any surface $F$. Thus, 
$\pi_1(\Omega CSTF, \vec l)$ is isomorphic to $Z(\vec l)<\pi_1(CSTF,
\vec l(a))$. From the proof of the Hansen's Proposition it follows that the
isomorphism is given by $p_*$. \qed
\end{emf}

\begin{prop}\label{split}
The fiberwise projectivization $PCSTF$ of $CSTF$ is isomorphic to 
$S^1\times STF$. 
\end{prop}

\begin{emf}\label{pfsplit}{\em Proof of Proposition~\ref{split}.\/}
A local orientation at $x\in F$ induces an orientation of the $S^1$-fiber of 
$STF$ over $x$, which changes if we change the local orientation at $x$.
Hence $STF$ is canonically oriented. The planes of the contact structure 
of $STF$ are canonically cooriented. Thus they are also canonically
oriented. The orientations of them induce a coherent orientation of the $\R
P^1$-fibers of $PCSTF\rightarrow STF$. Hence to prove the proposition it
suffices to construct a section of $PCSTF$ over $STF$. 

A point $x\in STF$ is described by $\pr_2(x)\in F$ 
and a cooriented contact element at  $\pr_2(x)$. Consider an arc $L_x$ of the
geodesic passing through $\pr_2(x)$ that is tangent to the contact element. 
Equip the arc with the coorientation coherent with the one of
the contact element. Choose an orientation of $L_x$ and lift 
it to an immersed arc in $STF$. The direction of $l_x$ at the lifting
of the preimage of $\pr_2(x)$ defines a point in the $S^1$-fiber of 
$CSTF$ over $x$ and, consequently, a point $\bar x$ in the $\R P^1$-fiber of 
$PCSTF$ over $x$. Clearly the point $\bar x$ 
in the $\R P^1$-fiber is independent of 
the choice of the orientation of $L_x$.
The desired section is given by mapping $x\in STF$ to the 
point $\bar x \in PCSTF$.
\qed
\end{emf}

\begin{prop}\label{commute}
Let $f_1\in \pi_1(CSTF,\vec l(a))$ and  $f_2\in \pi_1(STF,l(a))$ be
the classes of oriented (in some way) 
fibers of the $S^1$-fibrations $\pr^1: CSTF\rightarrow STF$ and 
$\pr^2: STF\rightarrow F$ respectively. Then: 

1. $f_1\alpha=\alpha f_1\in\pi_1(CSTF, \vec l(a)$ 
for any $\alpha\in \pi_1(CSTF,\vec l(a))$.

2. $f_2\alpha=\alpha f_2\in\pi_1(STF, l(a))$ for any $\alpha\in
\pi_1(STF, l(a))$ projecting to an orientation-preserving loop in $F$.

3. $f_2\alpha=\alpha f_2^{-1}\in \pi_1(STF, l(a))$ for any $\alpha\in
\pi_1(STF, l(a))$ projecting to an orientation-reversing loop in $F$.
\end{prop}

\begin{emf}{\em Proof of Proposition~\ref{commute}.\/}
Consider the double covering $p:CSTF\rightarrow PCSTF$. The homomorphism 
$p_*:\pi_1(CSTF)\rightarrow \pi_1(PCSTF)$ is injective, and it maps 
$f_1$ to $f^2\in \pi_1(PCSTF)$. (Here $f$ is 
the class of an oriented $S^1$-fiber of 
$PCSTF\rightarrow STF$.) 
Proposition~\ref{split} implies that $f$ 
is in the center of $\pi_1(PCSTF)$, and we have proved the first statement
of the proposition.

If we move an oriented fiber along the loop 
$\alpha\subset STF$, then in the end it comes to itself either 
with the same or
with the opposite orientation. Thus for any $\alpha\in \pi_1(STF, 
l(a))$ either $f_2\alpha=\alpha f_2$ or $f_2\alpha=\alpha f_2^{-1}$. 
A local orientation of the neighborhood of a point in $F$ 
induces an orientation of the fiber of $\pr^2$ over the point.
Combining these facts we get
the proof of the other two statements of the proposition.
\qed
\end{emf} 

\begin{prop}\label{reprfront}
For any $\alpha\in \pi_1(CSTF, d)$ there exists a Legendrian curve $l$ 
such that $\vec l(a)=d$ and $\vec l=\alpha\in \pi_1(CSTF, d)$.
\end{prop}

\begin{emf}{\em Proof of Proposition~\ref{reprfront}.\/}
Let $L$ be a front such that $\vec l(a)=d$ and 
$L=\pr^2_*\pr^1_*(\alpha)\in\pi_1(F, \pr^2\pr^1(d))$.
A small extra kink on $L$ corresponds to the multiplication of
$l\in\pi_1(STF)$ by $f_2^{\pm 1}$ depending on the side of
front the kink points to. 
Thus, adding extra kinks we can modify $L$, so that 
$l=\pr^1_*(\alpha)\in \pi_1(STF)$. In~\cite{Arnoldsplit} it is shown that an 
extra pair of adjacent cusps pointing to
opposite sides of a planar front $L_1$ 
corresponds to the multiplication of $\vec l_1$ by 
$f_1^{\pm 1}$ depending on the sign of the cusps. 
Since the addition of the
pair of cusps is done locally and $f_2$ is in the center of $\pi_1(CSTF)$
for any $F$ (see~\ref{commute}),
this fact is true for any surface $F$.
Thus we can modify $L$ so
that $\vec l=\alpha\in \pi_1(CSTF,d )$, and we have proved the
proposition.\qed
\end{emf}
\begin{prop}\label{pi1CSTF}
The group $\pi_1(CSTF)$ is isomorphic to the index two subgroup of 
$\pi_1(PCSTF)=\pi_1(S^1)\oplus\pi_1(STF)=\Z\oplus \pi_1(STF)$ consisting of
elements of the form (odd number, element 
projecting to an orientation-reversing loop in $F$) and 
(even number, element projecting to an
orientation-preserving loop in $F$).
\end{prop}   
  
\begin{emf}\label{pfpi1CSTF}{\em Proof of Proposition~\ref{pi1CSTF}.\/}
By proposition~\ref{reprfront} 
an element of $\pi_1(CSTF)$ can be realized as a lifting $\vec l$ of some
front $L$. Orientation-preserving fronts 
have even number of cusps, and orientation-reversing fronts
have odd number of cusps. 
There are only two (opposite) 
points in the $S^1$-fiber of
$\pr^1:CSTF\rightarrow STF$ over $x$ corresponding to a cusp point of a
front at $x$. (These points are identified under $p:CSTF\rightarrow PCSTF$.)
Thus for $\alpha\in\pi_1(CSTF)$ the projection of 
$p_*(\alpha)\in \pi_1(PCSTF)=\Z\oplus\pi_1(STF)$ to the $\Z$-summand is even 
provided that $\pr^2\pr^1(\alpha)$ is an orientation-preserving loop in $F$
and odd otherwise. The difference between the projections of 
$p_*(\alpha)$ and $p_*(\alpha f_1)$ to the $\Z$-summand in $\pi_1(PCSTF)$ is
two. Since $\pi_1(CSTF)$ is isomorphic to $p_*(\pi_1(CSTF))<\pi_1(PCSTF)$ we
get the statement of the proposition.\qed
\end{emf} 

\begin{emf}\label{gamma3}{\em Loop $\gamma_3$.\/}
Let $\mathcal C$ be a connected component of $\mathcal L$ and $L\in\mathcal C$ a
generic front. Let $\gamma_3\in\pi_1(\mathcal L, L)$ 
be the loop constructed below.

Deform $L$ along a generic path $t$ in $\mathcal C$, so that all cusps are
concentrated on a small piece $P$ of $L$ and the side of $L$ they point
to alternates. (The notion of side is locally well defined.) 
This is possible because
we can cancel a pair of adjacent cusps pointing to the same side of $L$, 
see Figure~\ref{cancel.fig}.

\begin{figure}[htbp]
 \begin{center}
  \epsfxsize 6cm
  \hepsffile{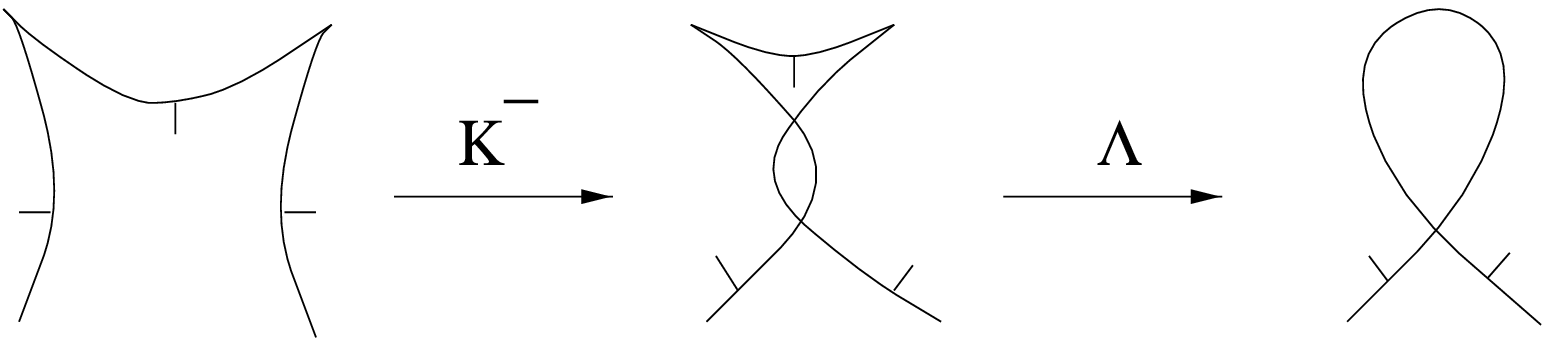}
 \end{center}
\caption{}\label{cancel.fig}
\end{figure}

If after this deformation the number of cusps is nonzero, then we take 
the last pair of cusps in $P$ and slide them along $L$ till they come to the
beginning of $P$. Then we shift all the cusps by two positions, so that $L$ gets
the shape it had before the sliding.
(If $L$ is an orientation-reversing front, then 
it can happen that there is only one cusp on $P$. Then to obtain $\gamma_3$
we slide the cusp twice around $L$ till it comes to the original position.) 
We require the deformation to be such that at
each moment of time points of $L$ located outside of a small neighborhood
of participating cusps do not move.

If after the deformation the number of cusps is zero (this happens 
if $\mu(L)=0$), then we perform a regular homotopy shown in
Figure~\ref{loop0.fig}.

\begin{figure}[htbp]
 \begin{center}
  \epsfxsize 8cm
  \hepsffile{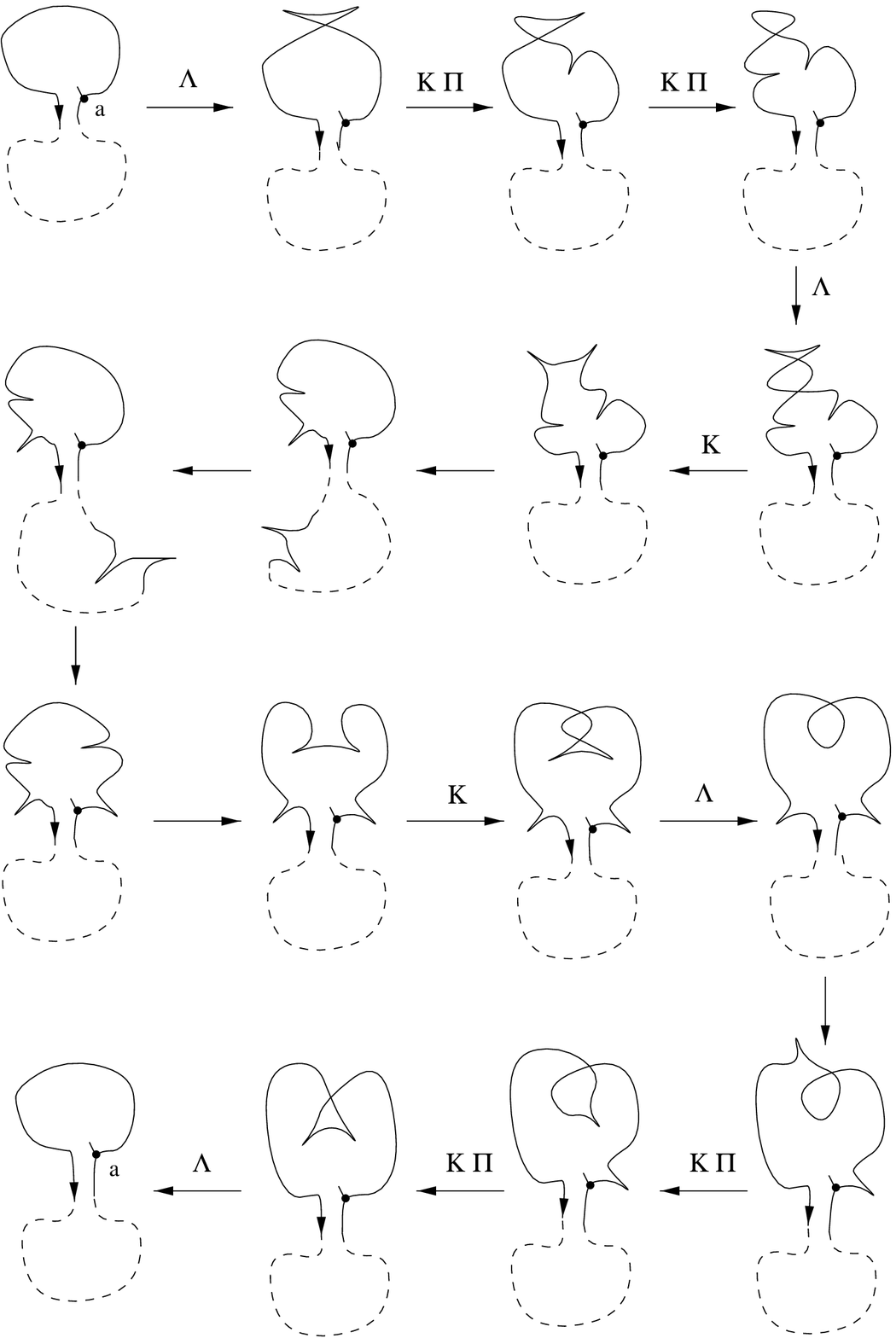}
 \end{center}
\caption{}\label{loop0.fig}
\end{figure}

Finally we deform $L$ to its original shape along $t^{-1}$.
\end{emf}

\begin{prop}\label{gammatraject}
1. Let $L$ be an orientation preserving front on $F$, then  
$\pr_*^1(t_a(\gamma_1))=f_2^{\pm 1}\in \pi_1(STF, l(a))$. 
(Here the sign of the power of $f_2$ depends on the orientation 
of the fiber we choose to define $f_2$.)

2. Let $L$ be a front on $F$, then $t_a(\gamma_3)=f_1^{\pm 1}\in
\pi_1(CSTF, \vec l(a))$. 
(Here the sign of the power of $f_1$ depends on the orientation 
of the fiber we choose to define $f_1$.)

\end{prop}
\begin{emf}{\em Proof of Proposition~\ref{gammatraject}.\/}
Under the deformation of $L$ described by $\gamma_1$ 
the point $L(a)$ never leaves 
a small neighborhood of its original position,
and the coorienting normal to $L$ at $L(a)$ is rotated
by $2\pi$. (The rotation happens when the kink passes through $L(a)$.) 
Thus the trajectory of $a$ under the lifting of
$\gamma_1$ to a loop in the space of Legendrian curves represents a class of
the fiber of the fibration $\pr^2$. Clearly, 
this trajectory coincides with $\pr^1_*(t_a(\gamma_1))$ and we have proved the
first statement of the Proposition.

From~\ref{commute} we know that $f_1$ is in the center of $\pi_1(CTSF)$.
Hence to prove the second statement it suffices
to prove the corresponding fact for a loop that is free
homotopic to $\gamma_3$. 
Thus we can assume that under the deformation
$\gamma_3$ the point $L(a)$ never leaves a small neighborhood of its
original position, and that $L$ has the property that 
all cusps of it are close to each other (in $L$), and the side of $L$ they
point to alternates. (The notion of side is locally well defined.)

One verifies that under $\gamma_3$ 
the total rotation angle of the coorienting normal at $L(a)$ is zero.
Thus the trajectory of $a$ under the lifting of 
$\gamma_3$ to a loop in the space of Legendrian curves in $STF$ represents
$1\in \pi_1(STF, l(a))$. Hence $t_a(\gamma_3)=f_1^k\in \pi_1(CSTF, \vec
l(a))$, for some $k\in \Z$. We have to show that
$k=\pm 1$.

One verifies that there are only two points in the $S^1$-fiber of $CSTF$ over 
$l(a)$ which correspond to the front having a cusp at $L(a)$. 
A lemma proved by Arnold~\cite{Arnoldsplit} says that 
under the deformations of the wave front shown in
Figure~\ref{rotate.fig} the velocity vector at the point $l(a)$ is turning
in the direction dependent only on, whether it is true or not, 
that after the deformation the
coorienting normal is pointing to the same direction as the curvature vector.
(In Figure~\ref{rotate.fig} the marked point is frozen into the surface
together with the coorienting normal at it.)

\begin{figure}[htbp]
 \begin{center}
  \epsfxsize 11cm
  \hepsffile{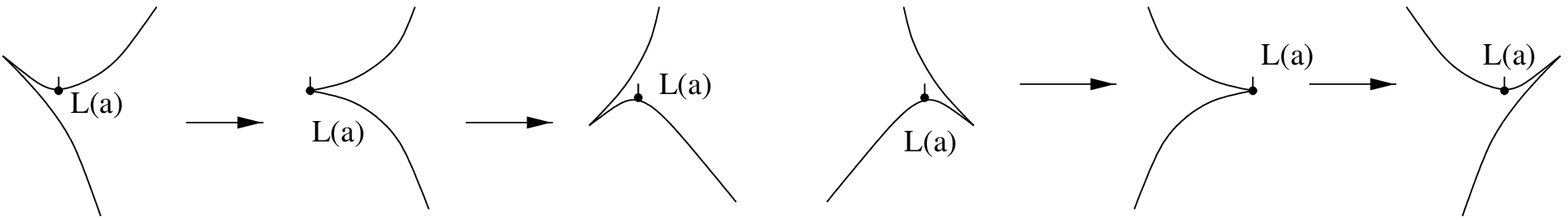}
 \end{center}
\caption{}\label{rotate.fig}
\end{figure}

It is clear that the loop $\gamma_3$ is free homotopic to $\gamma_3 '$,
in which the point $L(a)$ is frozen into $F$ together with the coorienting 
normal at it. Under the deformation described by $\gamma_3'$ 
the pair of cusps
passes through $L(a)$ in the way shown in Figure~\ref{cusppoint.fig}.
Using Arnold's lemma one verifies 
that under $\gamma_3'$  
the direction of the 
velocity vector of $l$ at $l(a)$ is rotated by the fiber of $\pr^1$.

\begin{figure}[htbp]
 \begin{center}
  \epsfxsize 8cm
  \hepsffile{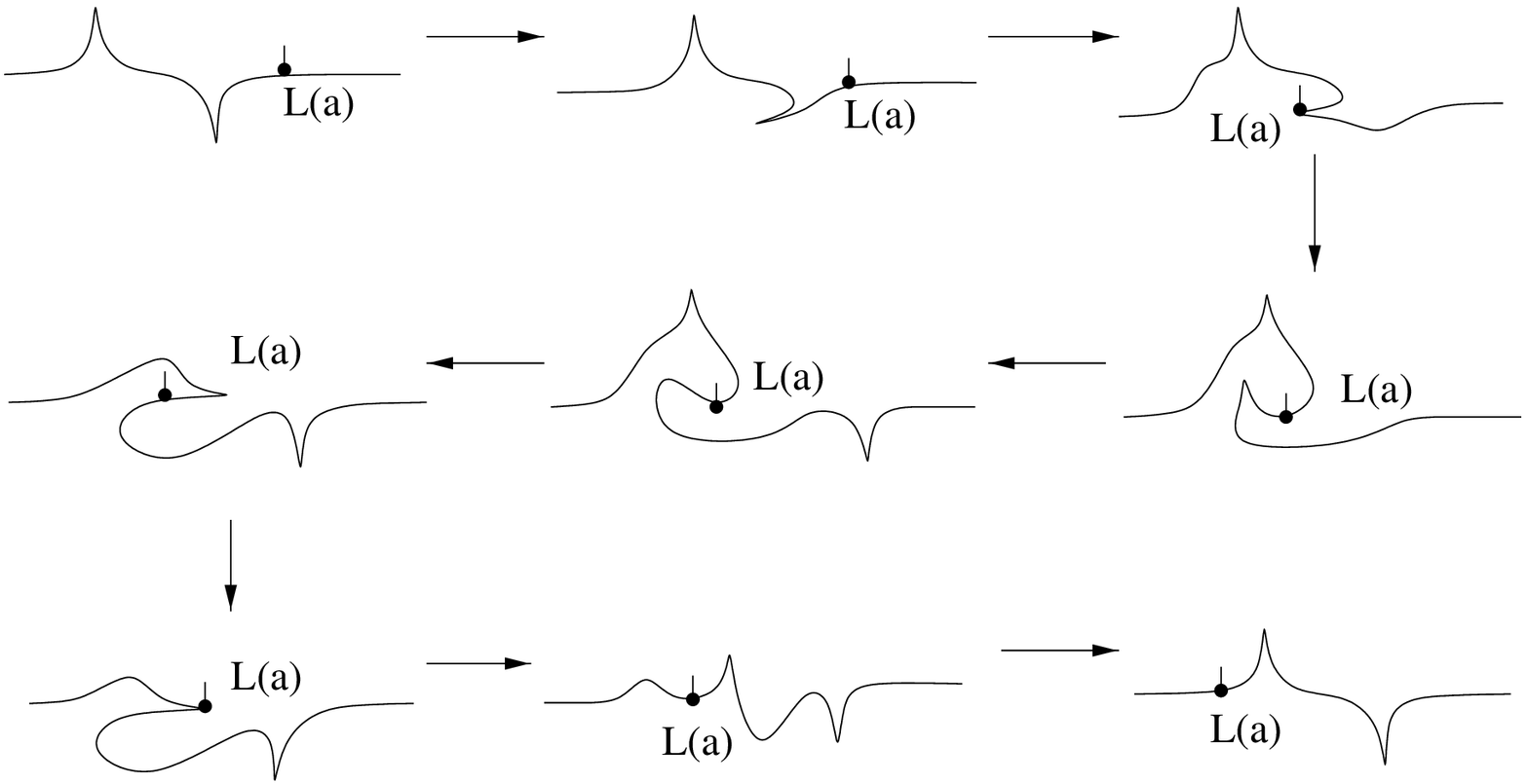}
 \end{center}
\caption{}\label{cusppoint.fig}
\end{figure}

This finishes the proof of the Proposition.
\qed
\end{emf}

\begin{prop}\label{zerodeltagamma3} $\Delta_{\overline{\St'}}(\gamma_3)=0$. 
\end{prop}

\begin{emf}{\em Proof of Proposition~\ref{zerodeltagamma3}.\/}
Clearly, the input into $\Delta_{\overline{\St'}}$ of the
deformation $r$ of $L$ to a front with cusps pointing to
alternating (locally well defined) sides of $L$ 
cancels 
out with the input into $\Delta_{\overline{\St'}}$ of the deformation along 
$r^{-1}$. 

Consider the case when $\mu(L)\not\in\{0, \pm 1\}$. 
No $T$-stratum crossings occur under the sliding of two cusps and 
the only inputs into $\Delta_{\overline{\St'}}(\gamma_3)$ 
come from the crossings of the $\Pi$-stratum. They occur only
when a cusp passes through a neighborhood of a double point $x$ of $L$, see
Figure~\ref{cuspst2.fig}. 
One verifies that for each double point $x$ the input corresponding to the 
first cusp passing through the neighborhood of $x$ cancels with the input 
corresponding to the second cusp passing through it.

If $\mu(L)=0$, then there are extra crossings of the $\Pi$-stratum, which occur
when we create (and later cancel) two pairs of cusps, one of which slides
along $L$. One verifies that the inputs of these extra crossings cancel out. 

If $\mu(L)=\pm 1$, then $L$ is orientation-reversing, and 
the only cusp present on $L$ slides twice 
along $L$ under $\gamma_3$. The input corresponding to the 
crossings of the $\Pi$-stratum  that occur in the neighborhood of a double
point $x$ under the first round of sliding cancels with the input under the 
second round of sliding. This finishes the proof of
Proposition~\ref{zerodeltagamma3}.\qed

\begin{figure}[htbp]
 \begin{center}
  \epsfxsize 7cm
  \hepsffile{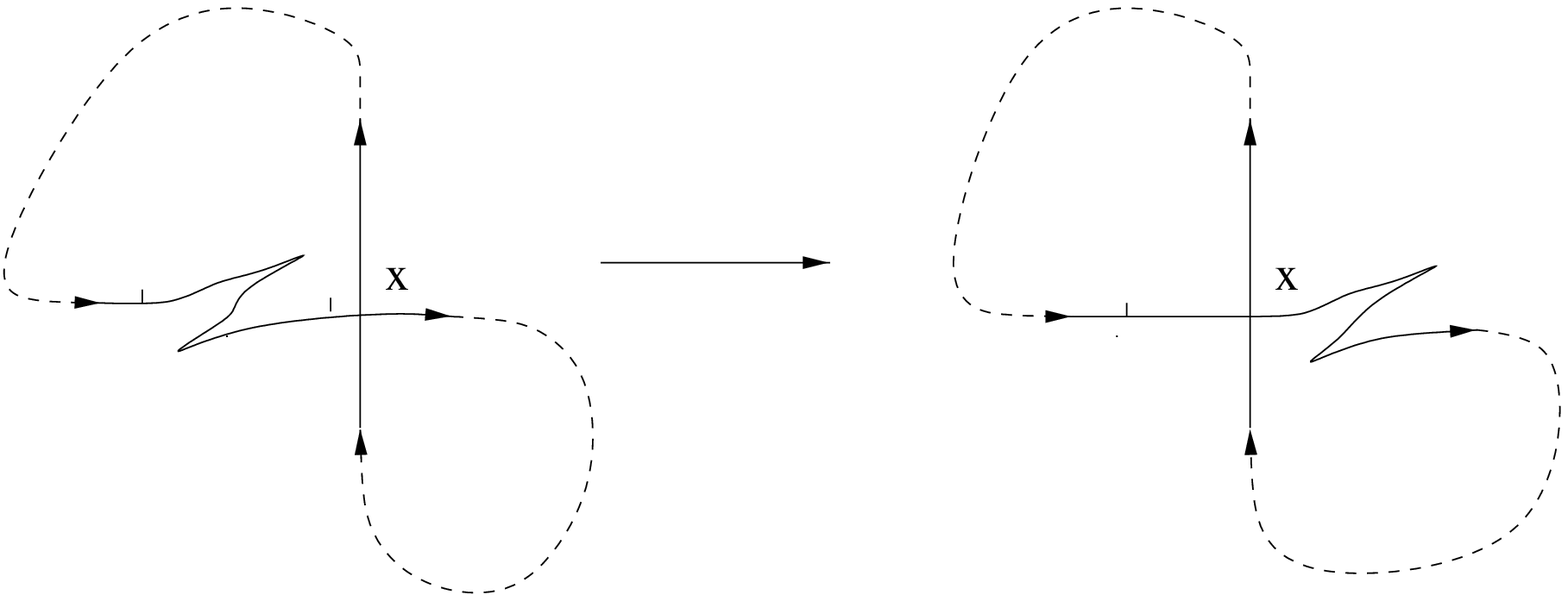}
 \end{center}
\caption{}\label{cuspst2.fig}
\end{figure}

\end{emf}

\begin{prop}\label{Preissman}
Let $F\neq S^2, T^2\text{ (torus), } \R P^2, K\text{ (Klein bottle)}$ 
be a surface (not necessarily compact or orientable),
and let $G$ be a nontrivial commutative subgroup of $\pi_1(F)$. 
Then $G$ is infinite cyclic 
and there exists a
unique maximal infinite cyclic $G'<\pi_1(F)$ containing $G$.
\end{prop}

\begin{emf}{\em Proof of Proposition~\ref{Preissman}.\/}
It is well known that any closed $F$, other than $S^2, T^2, \R P^2, K$,
admits a hyperbolic metric of a constant negative curvature. 
(It is induced from the universal covering of $F$ by the hyperbolic plane $H$.) 
The Theorem by A.~Preissman (see~\cite{Docarmo} pp. 258-265) 
says that if $M$ is a compact Riemannian manifold with a negative curvature,
then any nontrivial Abelian subgroup $G<\pi_1(M)$ is isomorphic to $\Z$.
Thus if $F\neq S^2, T^2, \R P^2, K$ is closed, then any nontrivial 
commutative $G<\pi_1(F)$ is infinite cyclic.

The proof of the Preissman's Theorem given in~\cite{Docarmo} is based on
the fact, that if $\alpha,\beta\in\pi_1(M)$ are nontrivial commuting
elements, then there exists a
geodesic in $\bar M$ (the universal covering of $M$) which is mapped to itself
under the action of these elements considered as deck transformations on
$\bar M$. Moreover, these transformations restricted to the geodesic act as
translations. This implies that if $F\neq S^2, T^2, \R P^2, K$ 
is a closed surface, then there exists a unique maximal infinite cyclic
$G'<\pi_1(F)$ containing $G$. This gives the proof of
Proposition~\ref{Preissman} for
closed $F$. 

If $F$ is not closed, then the statement of the Proposition is also
true because in this case $F$ is homotopy equivalent to a bouquet of
circles.
\qed
\end{emf}

\subsection{Proof of Lemma~\ref{zerojump}}
The proof is based on the constructions and the propositions of 
Subsection~\ref{facts}.
We start by making the following observation.
\begin{emf}\label{thoughtful}
In $\Z$ there are no elements of finite order. Thus  
if $m\neq 0$, then $\Delta_{\overline{\St'}}(q)\neq 0
\Leftrightarrow m\Delta_{\overline{\St'}}(q)=
\Delta_{\overline{\St'}}(q^m)\neq 0$. 
Hence, to prove Lemma~\ref{zerojump} 
it suffices to show that $\Delta_{\overline{\St'}}(q^m)=0$, 
for some $m\neq 0$.
\end{emf}

\begin{prop}\label{enoughproject}
Let $q_1,q_2\in \pi_1(\mathcal L, L)$ be loops such that:
\begin{equation} 
\pr_*^1(t_a(q_1))=\pr_*^1(t_a(q_2))\in \pi_1(STF, l(a))
\end{equation} 
Then $\Delta_{\overline{\St'}}(q_2)=\Delta_{\overline{\St'}}(q_1)$.
\end{prop}

\begin{emf}{\em Proof of the Proposition~\ref{enoughproject}.\/}
We know (see~\ref{gammatraject}) that $t_a(\gamma_3)=f_1$ 
(for a proper choice of an orientation of the fiber used to define $f_1$).
The kernel of the homomorphism $\pr^1_*$
is generated by $f_1$, which is in the
center of $\pi_1(CSTF, \vec l(a))$, see~\ref{commute}.
Thus $t_a(q_2)=t_a(q_1)f_1^j$ for some 
$j\in \Z$. 
Proposition~\ref{centralizer} implies that $q_2=q_1\gamma_3^j$. 
Proposition~\ref{zerodeltagamma3} says that
$\Delta_{\overline{\St'}}(\gamma_3)=0$. Hence   
$\Delta_{\overline{\St'}}(q_1)=\Delta_{\overline{\St'}}(q_2)$.
\qed
\end{emf}

{\em We first prove Lemma~\ref{zerojump} 
for $F\neq S^2, \R P^2, T^2, K$, and then separately 
for the cases of
$F=S^2, \R P^2, T^2, K$.\/}

\begin{emf}\label{other}{\em Case $F\neq S^2, T^2, \R P^2, K$.\/} 
Proposition~\ref{centralizer} says that 
$\pi_1(\mathcal L, L)=Z(\vec l)<\pi_1(CSTF, \vec l(a))$. The corresponding
isomorphism (see Section~\ref{H-principle}) 
maps $q\in\pi_1(\mathcal L, L)$ to $t_a(q)\in\pi_1(CSTF, \vec l(a))$.
Thus for any  $q\in\pi_1(\mathcal L, L)$ 
the elements $t_a(q)$ and $\vec l$ commute in 
$\pi_1(CSTF, \vec l(a))$. Hence $L=\pr^2_*\pr^1_*(\vec l)$ 
commutes with $\pr^2_*\pr^1_*(t_a(q))$ in $\pi_1(F, L(a))$.
Proposition~\ref{Preissman} implies that 
there exists an isomorphic to $\Z$ subgroup of 
$\pi_1(F,L(a))$ generated by some $g\in
\pi_1(F,L(a))$ that contains both of these loops.  
Then $L=g^m$ and
$\pr^2_*\pr^1_*(t_a(q))=g^n$, for some $m,n\in\Z$.

Consider a wave front $L_1$ such that $\vec l_1(a)=\vec l(a)$, and  
$g=L_1\in \pi_1(F, L(a))$. 
The kernel of $\pr_*^2$ is generated by $f_2$
which has infinite order in $\pi_1(STF)$ for our surfaces $F$. 
Proposition~\ref{commute} allows us to interchange $f_2$ 
with the other elements of $\pi_1(STF,l(a))$. Thus 
$l=l_1^mf_2^i$, and 
$\pr^1_*(t_a(q))=l_1^nf_2^j$ (in $\pi_1(STF)$) for some $i,j\in \Z$.

We prove Lemma~\ref{zerojump} separately in
the cases of $m\neq 0$   
and $m=0$ in respectively 
Subsubsections~\ref{nontrivial} and~\ref{trivial}.
(Geometrically these two cases correspond to  
$L=1\in\pi_1(F)$ and $L\neq 1\in \pi_1(F)$ respectively.)
\end{emf} 

\begin{emf}\label{nontrivial}{\em Case $m\neq 0$.\/}
To prove Lemma~\ref{zerojump} it suffices to show that 
$\Delta_{\overline{\St'}}(q^m)=0$ (see~\ref{thoughtful}). 
We do it by constructing $\alpha\in \pi_1(\mathcal L, L)$ such that 
$\pr_*^1(t_a(\alpha))=\pr_*^1(t_a(q^m))$ and 
$\Delta_{\overline{\St'}}(\alpha)=0$. After this, 
Proposition~\ref{enoughproject} implies the statement of the lemma.

One can show that: 
\begin{equation}\label{ident}
\pr_*^1(t_a(q^m))=l^nf_2^k \text{ for some } k\in \Z.
\end{equation}
 
For an orientation-preserving $g$ this follows 
from the following calculation (which uses~\ref{commute}):  
\begin{equation}
\pr_*^1(t_a(q^m))=\bigl(\pr_*^1(t_a(q))\bigr)^m=(l_1^n f_2^j)^m=
(l_1^mf_2^i)^nf_2^{jm-ni}=l^nf_2^{jm-ni}.
\end{equation}
(Recall that $\pr^1_*(t_a(q))=l_1^nf_2^j$ for some $j\in \Z$, see~\ref{other}.)

For an orientation-reversing $g$
this follows from the similar calculation (also based on 
Proposition~\ref{commute}).

The fact that $\pr^1_*(t_a(q^m))$ commutes with $l$ 
(since $t_a(q^m)=(t_a(q))^m\in Z(\vec l)$) 
and Proposition~\ref{commute} imply
that $k=0$ in~\eqref{ident}, 
provided that $L$ is an orientation-reversing front.

Consider the case of $L$ being an orientation-preserving front.
Proposition~\ref{gammatraject} says that 
$\pr_*^1(t_a(\gamma_1))=f_2$ (for a proper choice of
the orientation of the fiber used to define $f_2$). 
Hence $\pr_*^1(t_a(\alpha))=\pr_*^1(t_a(q^m))$
for $\alpha\in\pi_1(\mathcal L, L)$ which is:  
$n$ times sliding of $L$
along itself (induced by a rotation of the parameterizing circle) 
composed with $\gamma_1^{j}$.

If $L$ is an orientation-reversing front, 
then as we have shown above $\pr_*^1(t_a(q^m))=l^n$. 
Hence $\pr_*^1(t_a(\alpha))=\pr_*^1(t_a(q^m))$ for 
$\alpha\in\pi_1(\mathcal L, L)$ which is: $n$ times sliding of $L$
along itself (induced by a rotation of the parameterizing circle). 

No discriminant crossings occur under the sliding of $L$ along
itself, and $\Delta_{\overline{\St'}}(\gamma_1)=0$ by the assumption of the 
the lemma. Hence $\Delta_{\overline{\St'}}(\alpha)=0$.
Proposition~\ref{enoughproject} implies that $\Delta_{\overline{\St'}}(q^m)=0$.
Thus, we have proved (see~\ref{thoughtful}) Lemma~\ref{zerojump} 
for $F\neq S^2, \R P^2, T^2, K$ and $m\neq 0$. 
\end{emf}

\begin{emf}\label{trivial}{\em Case $m=0$.\/}
If $m=0$, then $L=1\in\pi_1(F,L(a))$. 
We want to construct $\alpha\in \pi_1(\mathcal L, L)$ such that 
$\pr_*^1(t_a(\alpha))=\pr_*^1(t_a(q^2))$ and 
$\Delta_{\overline{\St'}}(\alpha)=0$.
(After this the statement follows
from~\ref{thoughtful} and Proposition~\ref{enoughproject}.)

For any $q\in\pi_1(\mathcal L, L)$ the projection 
$\pr^2_*\pr_*^1(t_a(q^2))$ is an orientation-preserving loop in $F$.
A straightforward verification shows that $\alpha$ can be obtained 
by a composition of $\gamma_1^{\pm 1}$
(see~\ref{gamma1}) and loops constructed as follows:
  
Push $L$ into a small disc by a generic regular homotopy $r$. 
Slide this small disc along a smooth orientation-preserving curve in $F$ 
and return $L$ to its original shape along $r^{-1}$. 

Clearly, the inputs of $r$ and $r^{-1}$ into 
$\Delta_{\overline{\St'}}$ cancel out, and no discriminant
crossings happen when we slide the  small disc (containing $L$)
along a loop in $F$. 
By the assumption of the lemma $\Delta_{\overline{\St'}}(\gamma_1)=0$.
Thus $\Delta_{\overline{\St'}}(\alpha)=0$.
Proposition~\ref{enoughproject} implies that
$\Delta_{\overline{\St'}}(q^2)=0$, and
we have proved (see~\ref{thoughtful}) Lemma~\ref{zerojump}  
for $F\neq S^2, \R P^2, T^2, K$.
\end{emf}

\begin{emf}\label{S^2}{\em Case $F=S^2$.\/}
One verifies that $\pi_1(STS^2)=\pi_1(\R P^3)=\Z_2$. 
Thus
$\pr_*^1(t_a(q^2))=\pr_*^1(t_a(1))=1\in \pi_1(STF,l(a))$, for any $q\in
\pi_1(\mathcal L, L)$. 
Proposition~\ref{enoughproject} 
implies that $\Delta{\overline{\St'}}(q^2)=0$. 
This finishes (see~\ref{thoughtful}) the proof 
of Lemma~\ref{zerojump} for $F=S^2$. \qed
\end{emf}
 
\begin{emf}\label{T^2}{\em Case $F=T^2$.\/} 
One verifies that
$\pi_1(STT^2)=\Z\oplus\Z\oplus\Z$. As before we fix 
$q\in \pi_1(\mathcal L, L)$ and construct $\alpha\in\pi_1(\mathcal L, L)$ such that 
$\pr_*^1(t_a(\alpha))=\pr_*^1(t_a(q))$.

One verifies that $\alpha$ can be expressed through $\gamma_1$ and loops
$\gamma_4$ and $\gamma_5$ that are  
slidings of $L$ along the 
unit vector fields parallel to the meridian and longitude
of $T^2$ respectively. 

Since $\Delta_{\overline{\St'}}(\gamma_1)=0$ by the assumption of the Lemma,
and no discriminant crossings occur under $\gamma_4$ and $\gamma_5$, we get
that $\Delta_{\overline{\St'}}(\alpha)=0$. Proposition~\ref{enoughproject}
implies that $\Delta_{\overline{\St'}}(q)=0$. 
This finishes the proof of Lemma~\ref{zerojump} for $F=T^2$.\qed
\end{emf}

\begin{emf}\label{RP^2}{\em Case $F=\R P^2$.\/}
One verifies that $\pi_1(ST \R P^2)=\Z_4$. 
Thus $\pr_*^1(t_a(q^4))=\pr_*^1(t_a(1))=1\in \pi_1(STF,l(a))$ 
for any $q\in \pi_1(\mathcal L, L)$. 
Proposition~\ref{enoughproject}
implies that $\Delta{\overline{\St'}}(q^4)=0$. 
This finishes (see~\ref{thoughtful}) the proof
of Lemma~\ref{zerojump} for $F=\R P^2$.\qed
\end{emf}

\begin{emf}\label{K}{\em Case $F=K$.\/}
Proposition~\ref{centralizer} says that 
$\pi_1(\mathcal L, L)$ is isomorphic to $Z(\vec l)<\pi_1(CSTK, \vec l(a))$.
The kernel of the homomorphism $\pr^1_*$ 
is generated by $f_1$ which is in the center of
$\pi_1(CSTF, \vec l(a))$, see~\ref{commute}.
Thus $\pr^1_*(Z(\vec l))$ is isomorphic to $Z(l)$, the 
centralizer of $l\in \pi_1(STF, l(a))$. 
We show that a certain power of any element
of $Z(l)$ can be represented as $\pr_*^1(t_a(\alpha))$, for some $\alpha\in
\pi_1(\mathcal L, L)$ such that $\Delta{\overline{\St'}}(\alpha)=0$. 
This implies the statement of the Lemma for $F=K$.

Consider $K$ as a quotient of a rectangle modulo the identification on its
sides shown in Figure~\ref{klein1.fig}. 
We can assume that $L(a)$ coincides
with the image of a corner of the rectangle. 
Let $L_1$ and $L_2$ be fronts such that 
$\vec l(a)=\vec l_1(a)=\vec l_2(a)$, $L_1=c\in \pi_1(K, L(a))$, $L_2=d\in \pi_1(K, L(a))$. 
(Here $c$ and $d$ are the elements of $\pi_1(K)$
realized by the images of the sides of the rectangle used to construct $K$,
see Figure~\ref{klein1.fig}.)
One can show that: 
\begin{equation}\label{present}
\pi_1(STK, l(a))=\bigl \{l_1, l_2, f_2\big | l_2 l_1^{\pm 1}
=l_1^{\mp 1} l_2, \text{ } l_2f_2^{\pm
1}=f_2^{\mp 1} l_2, \text{ } l_1 f_2=f_2 l_1 \bigr \}.
\end{equation}

The second and the third relations in this presentation follow 
from Proposition~\ref{commute}. 
To get the first relation one notes that  
the identity $d c^{\pm 1}=c^{\mp 1} d\in\pi_1(K, L(a))$ implies 
$l_2 l_1^{\pm 1}=l_1^{\mp 1} l_2 f_2^k$, for some $k\in \Z$.
But $l_2^2$ commutes with $l_1$, since they can be lifted to $STT^2$ the
fundamental group of which is Abelian. Hence $k=0$.

Using relations~\eqref{present} 
one calculates $Z(l)$.
(Note, that these relations allow one to present (in a unique way) an 
element of $\pi_1(STK, l(a))$ as 
$l_1^i l_2^j f_2^k$, for some $i,j,k\in \Z$.) 

This group appears to be:
\begin{description}
\item[a] The whole group $\pi_1(STK,l(a))$ provided that 
$l=l_2^{2i}$, for some $i\in \Z$.

\item[b] A subgroup of $\pi_1(STK, l(a))$ isomorphic to $\Z\oplus\Z\oplus\Z$
provided that $l=l_1^i l_2^{2j} f_2^k$, 
for some $i,j,k\in\Z$ such that $i\neq 0$ or $k\neq 0$.
This subgroup is generated by $\{l_1, l_2^2,f_2\}$.

\item[c] A  subgroup of $\pi_1(STK, l(a))$ isomorphic to $\Z$
provided that $l=l_1^i l_2^{2j+1} f_2^k$ for some $i,j,k\in\Z$.
This subgroup is generated by $\alpha_l=l_1^i l_2^{1} f_2^k$. (Note that
$\alpha_l^2=l_2^2$, and $l=(\alpha_l)^{2j+1}$.)
\end{description}

Using~\eqref{present} one verifies that:
\begin{description}
\item[a] If $L$ is an orientation-preserving front on $K$, then a certain
power of any element of $Z(l)$ can be obtained as 
$\pr_*^1(t_a(\alpha))$, for $\alpha$ being a product of powers of
$\gamma_1$, 
$\gamma_2$ (see~\ref{gamma2}), and 
$\gamma_4$ described below.

\item[b] If $L$ is an orientation-reversing front on $K$, then a certain
power of any element of $Z(l)$ can be obtained as 
$\pr_*^1(t_a(\alpha))$, for $\alpha$ being a power of
$\gamma_5$ described below.
\end{description}

Consider a loop $\beta$ 
in the space of all autodiffeomorphisms of
$K$, which is the sliding of $K$ along the unit vector field parallel 
to the curve $d$ on $K$. 
(Note that $K$ has to slide twice along itself under this loop before every
point of it
comes to the original position.) The loop $\gamma_4$ is the 
sliding of $L$ induced by $\beta$.

The loop $\gamma_5$ is the sliding of $L$ along itself
induced by a rotation of the parameterizing circle.

No discriminant crossings occur under $\gamma_4$ and
$\gamma_5$. By the assumption of the Lemma 
$\Delta_{\overline{\St'}}(\gamma_1)=0$ and 
$\Delta_{\overline{\St'}}(\gamma_2)=0$ (when $\gamma_2$ is well defined). 
Thus $\Delta_{\overline{\St'}}(\alpha)=0$, and because of the reasons explained 
in the beginning of Subsubsection~\ref{K} we have proved
Lemma~\ref{zerojump} for $F=K$. 
\begin{rem}
One can verify that for the front on the Klein bottle shown in
Figure~\ref{klein3.fig} the identity $\Delta_{\overline{\St}}(\gamma_2)=0$
does not follows from
$\Delta_{\overline{\St}}(\gamma_1)=0$. This means that
the condition $\Delta_{\overline{\St}}(\gamma_2)=0$ is needed for 
the integrability of $\psi$.

\begin{figure}[htbp]
 \begin{center}
  \epsfxsize 3cm
  \hepsffile{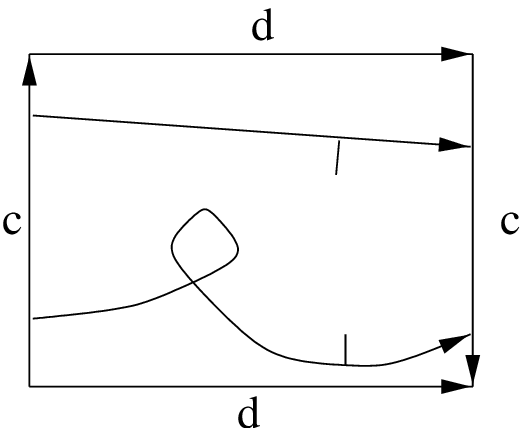}
 \end{center}
\caption{}\label{klein3.fig}
\end{figure}

\end{rem}

\end{emf}

This finishes the proof of Lemma~\ref{zerojump}.
\qed

\section{Proof of Theorem~\ref{barJ+orient}}\label{pfbarJ+orient}
Theorem~\ref{integrability} implies that to prove the
Theorem~\ref{barJ+orient} it suffices to show that
$\Delta_{\overline{J^+}}(\gamma_1)=0$.

Clearly the input into $\Delta_{\overline{J^+}}$ 
of the deformation $r$ of $L$ to a
front with two opposite kinks cancels out with 
the input of $r^{-1}$. 
Thus it suffices to show that $\Delta_{\overline{J^+}}$ 
under the sliding of the kink along $L$ is zero. The only crossings 
of the $K^+$-stratum,
which occur under this sliding happen either when the kink passes
through a small neighborhood of a double point or of a cusp of $L$.

\begin{figure}[htbp]
 \begin{center}
  \epsfxsize 10cm
  \hepsffile{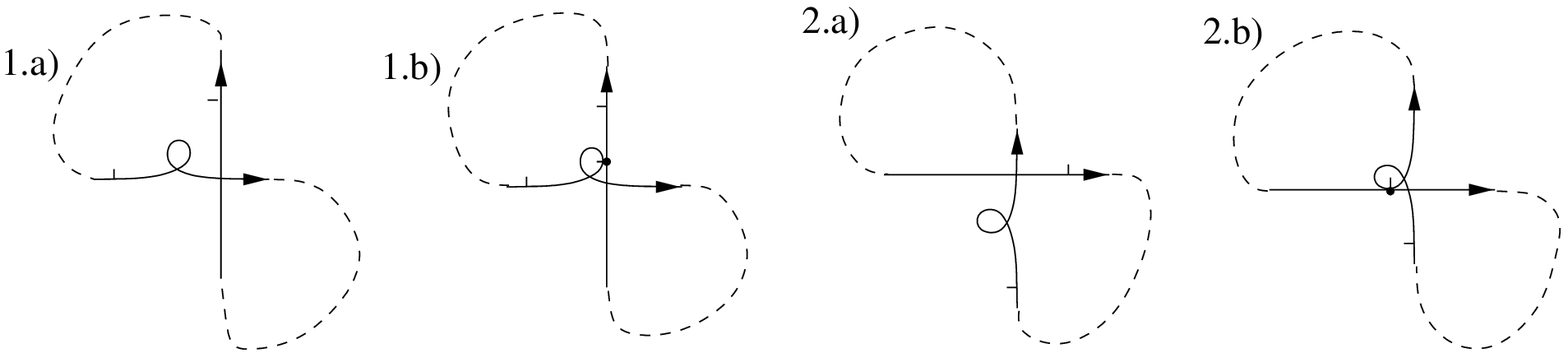}
 \end{center}
\caption{}\label{kinkcross.fig}
\end{figure}

The kink passes twice through each double point of $L$. (Once along each
intersecting branch.) As one can verify (see Figure~\ref{kinkcross.fig}) 
the two $K^+$-equivalence classes corresponding to
these events are equal and the signs of the corresponding 
$K^+$-stratum crossings are opposite. Thus the corresponding two
terms in $\Delta_{\overline{J^+}}$ cancel out. (One can verify that this part
of the proof would not go through, if $F$ is nonorientable and the double
point separates the front into two orientation-reversing loops.)

One can see (using Figure~\ref{kinkcusp.fig}) that 
either two or zero 
crossings of the $K^+$-stratum occur under the passage of the
kink through a neighborhood of a cusp. If the number of crossings 
is zero, then clearly there is no input
into $\Delta_{\overline{J^+}}$. In the case of two crossings one verifies that 
the signs of them are opposite and the corresponding $K^+$-equivalence
classes are equal. Thus the corresponding two terms of
$\Delta_{\overline{J^+}}$ also cancel out, and we have proved that 
$\Delta_{\overline{J^+}}(\gamma_1)=0$.
This finishes the Proof of Theorem~\ref{barJ+orient}. 
\qed

\section{Proof of Theorem~\ref{interpretation}}\label{pfinterpretation}
We prove the statement of the theorem only for the mapping $\psi$.
The proof of the statements about $\psi^+$, $\psi^-$, and $\psi^{\pi}$ 
is obtained in the similar way.

To show that $\psi$ is surjective we take 
$\alpha=(\delta_1, \delta_2, \delta_3, i)\in R_T$ 
and construct $L\in T$ which realizes the $T_i$-equivalence class of
$\alpha$. Consider $L'\in T$ for which 
the element $(\delta_1', \delta_2', \delta_3', i')\in R_T$ corresponding to
it is such that $pr^2_*(\delta_1)=pr^2_*(\delta_1')$,
$pr^2_*(\delta_2)=pr^2_*(\delta_2')$, and
$pr^2_*(\delta_3)=pr^2_*(\delta_3')$ in $\pi_1(F, \pr^2(d))$. 
Then $\delta_1=\delta_1'f_2^k$,  $\delta_2=\delta_2'f_2^m$
$\delta_3=\delta_3'f_2^n$ 
for some $k,m,n\in \Z$.
One verifies that a small extra kink located on one of the three loops of 
$L$ corresponds to the multiplication of 
the element of $\pi_1(STF)$ corresponding 
to this loop by 
$f_2^{\pm 1}$.  (Here the sign depends on which (locally well defined) 
side of the loop the kink points to.) Using this operation we obtain the
front $L$ corresponding to $(\delta_1, \delta_2, \delta_3, j)\in R_T$, for
some $j\in \Z$. Adding an extra pair of cusps of the same sign  we can
change $\mu(L)$ by $\pm 2$.
The three elements of $\pi_1(STF)$ corresponding to $L$ are not
changed by this operation. One can easily show that $i-j$ is even, and hence 
we can change $L$ so that it represent the $T_i$-equivalence class of
$\alpha$. Hence $\psi$ is surjective.

Below we show that $\psi$ is
injective. Let $L_1, L_2\in T$ be $T_i$-equivalent 
fronts. To prove the Theorem we construct a path in the normalization 
of the triple point part of the discriminant that connects the two fronts. 
We deform the fronts, so that under the 
lifting of fronts to Legendrian curves the preimages of triple points are
mapped to the same point $c\in STF$. Let $s_1$ and $s_2$ be the elements
of $\pi_1(STF)\times\pi_1(STF)\times\pi_1(STF)$ 
corresponding to the
deformed fronts. 
Since $L_1$ and $L_2$ are $T_i$-equivalent, one can transfer $s_1$ to $s_2$
by a consequent actions of elements $j\in \Z_3$, $i=(i_1,i_2, i_3)\in \Z^3$,
and $\xi\in \pi_1(STF, c)$. We can assume that
$\xi(i(s_1))=s_2$. 

\begin{figure}[htbp]
 \begin{center}
  \epsfxsize 8cm
  \hepsffile{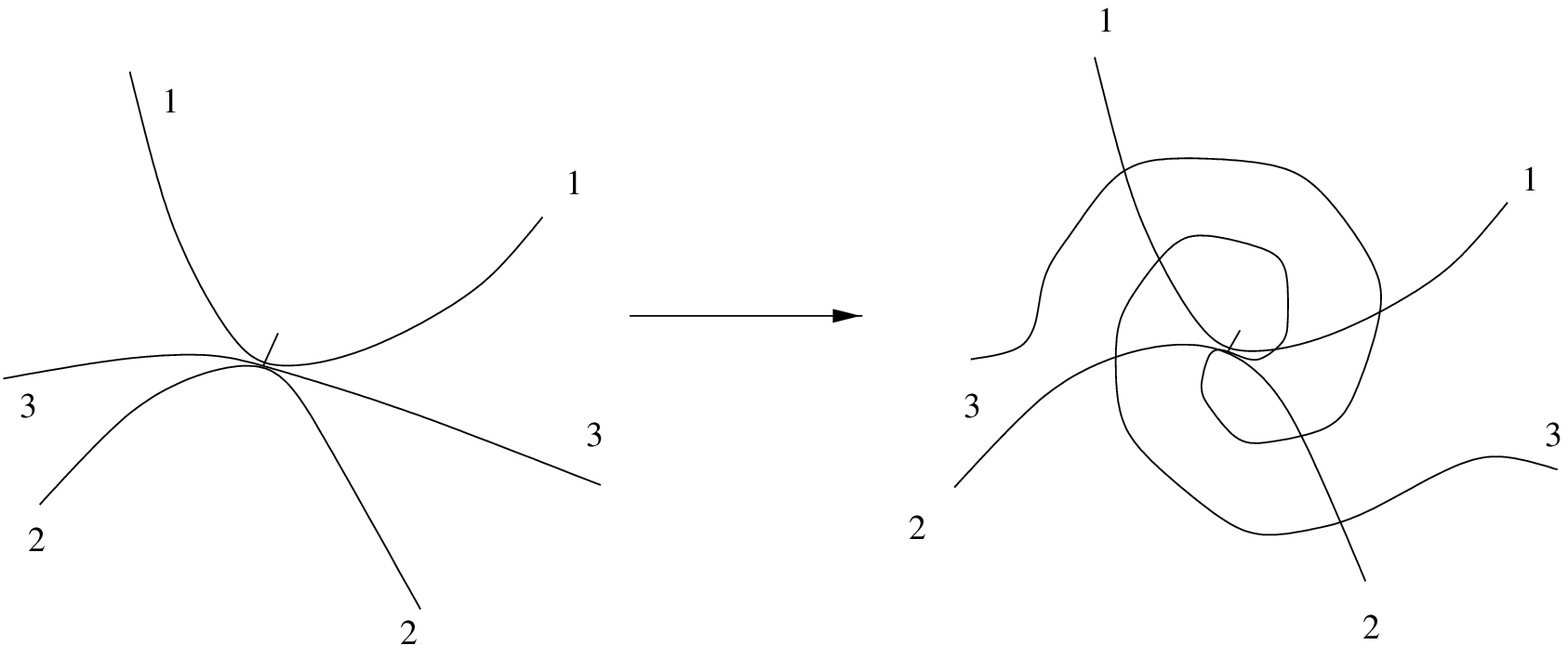}
 \end{center}
\caption{}\label{rotatebranch.fig}
\end{figure}

The local positive rotation by $2\pi$ of one of the three branches of $L_1$ 
passing through the triple point (see Figure~\ref{rotatebranch.fig}) 
induces the 
multiplication on the right by $f_2$ of one of the three loops of $s_1$,
and the multiplication on the left by $f_2^{-1}$ of the next loop of
$s_1$. Clearly this rotation does not change the $\overline T$-equivalence
class corresponding to $L_1$. Applying this rotation sufficiently many times
we deform $L_1$ so that $\xi(s_1)=s_2$.

Let $L$ be a front such that $l=\xi\in \pi_1(STF)$, cf.~\ref{reprfront}.  
Consider a diffeotopy
$\phi_t$, $t\in S^1$, of the small neighborhood of $\pr_2(c)$ in $F$ 
such that $\phi_t(\pr_2(c))=L(t)$ and the differential of $\phi_t$ maps the 
the coorienting normal to $L_1$ at $\pr_2(c)$ to the coorienting normal to 
$L$ at $L(t)$.  
This diffeotopy can be extended to the diffeotopy $\Phi_t$ of the whole $F$.
The diffeotopy $\Phi_t$ induces a deformation of $L_1$ that does not change
its $\overline T$-equivalence class. Clearly $s_1=s_2$ for the deformed
$L_1$.

As it shown in Figure~\ref{rotate.fig}, one can slide a cusp through a 
point of the front in such a way that a point and the coorienting normal at
it do not move under the deformation. For both $L_1$ and $L_2$ slide all the
cusps through a triple point, so that they are located on the first of the
three loops of $s_1$ and $s_2$. Cancel all the pairs of cusps of different sign.
(Clearly $s_1$ and $s_2$ are not changed during the process.)

Compare the directions of the velocity vectors of the corresponding branches of $L_1$ and $L_2$
at the triple point. One verifies that either they are coherent in all the
three pairs or they are opposite in all the three pairs. 
If they are
opposite, then take the last cusp on the saw-like piece of $L_1$
and slide it around $L_1$ (as it described above) till it comes to the
beginning of the saw-like piece. One verifies that after this 
the directions of
velocity vectors of the three branches of $L_1$ 
change sign and become coherent with the directions of the velocity vectors
of the branches of $L_2$. Deform $L_1$ and $L_2$ so that they are identical
in the neighborhood of the triple point.

One can lift an immersed oriented and cooriented interval $\alpha$ to an arc in
$STF$ by mapping a point to the direction of the
velocity vector at it or by mapping it to the direction of the coorienting
vector at it. The two liftings are denoted by $\alpha^v$ and $\alpha^c$
respectively. 
Consider a pair $\alpha_1, \alpha_2$ of immersed oriented cooriented
intervals that are identical in the neighborhood of the end points.
Clearly $\alpha_1^c(0)=\alpha_2^c(0)$, 
$\alpha_1^c(1)=\alpha_2^c(1)$, $\alpha_1^v(0)=\alpha_2^v(0)$, and
$\alpha_1^v(1)=\alpha_2^v(1)$. One verifies that if $\alpha_1^c$ and $\alpha_2^c$
are homotopic as arcs with fixed end points, then $\alpha_1^v$ and
$\alpha_2^v$ are also homotopic as arcs with fixed end points.
A statement
proved by Inshakov~\cite{Inshakov} says that if $\alpha_1^v$ and
$\alpha_2^v$ are homotopic as arcs with fixed points, then $\alpha_1$ is
homotopic to $\alpha_2$ in the class of immersed arcs with fixed end points
and velocity vectors at them. Clearly 
this homotopy $H(t, x): I\times I\rightarrow 
F$ can be chosen so that $H(t, y)=\alpha_1(y)=\alpha_2(y)$, for all $t\in I$
and $y\in [0, \epsilon)\cup (1-\epsilon, 1]$, for some $\epsilon >0$.
(Recall that $\alpha_1$ and
$\alpha_2$ are assumed to be identical in the neighborhood of the end
points.)

The triple point separates a front into three closed arcs. 
Consider a closed arc
$\alpha_1$ of $L_1$ containing all the cusps and the corresponding arc
$\alpha_2$ of $L_2$. Deform the two arcs so that they are identical on
a small piece in the beginning of them that contains all the cusps. Denote by 
$\beta_1$ and $\beta_2$ the subarcs of them where they are still different. 
Since $\alpha_1^c$ and $\alpha_2^c$ are homotopic as arcs with fixed end points 
(they correspond to the
same summand in $s_1=s_2$), we get that $\beta_1^c$ and $\beta_2^c$ are
homotopic as arcs with fixed end points. As it was said
above, this implies that $\beta_1^v$ and $\beta_2^v$ are homotopic as arcs
with fixed end points, and that $\beta_1$ and
$\beta_2$ are homotopic in the class of immersed arcs with fixed end points
and velocity vectors at them. Performing the homotopy we make $\alpha_1$ and
$\alpha_2$ identical. Other pairs of corresponding arcs of the two fronts are
deformed to each other in the similar way. For them the proof is even
simpler since they do not contain cusps. This finishes the proof of
Theorem~\ref{interpretation}.

(The last two steps of the proof can be done easier if one uses the relative
version of the $h$-principle proved for the Legendrian immersions by
T.~Duchamp in his unpublished preprint~\cite{Duchamp}.)

An important observation is that the constructed path in the normalization
of the triple point part of the discriminant can be slightly perturbed so
that its projection to the discriminant crosses only strata of codimension
two and the crossings are transversal. Analogous facts
are true in the case of the other three statements of the Theorem.
\qed

\subsection{Proof of Proposition~\ref{orient}}\label{pforient}
Since we consider only the equivalence classes appearing in $\mathcal C$, it is
clear that all three mappings are surjective. Hence it suffices to show that
these mappings are injective.

Let $L_1, L_2\in \mathcal C$ be $K^+$-equivalent fronts. Clearly
$\mu(L_1)=\mu(L_2)$, which means that $L_1$ and $L_2$ are
$K_i^+$-equivalent. (The free homotopy class of a mapping of $B_2$ is the
same as the element of $\pi_1(STF)\oplus\pi_1(STF)$ modulo the conjugation 
of both summands in it by the element of $\pi_1(STF)$.) Now
Theorem~\ref{interpretation} implies that $\psi^+$ is injective.

To prove that $\psi^-$ is injective we note that 
the free homotopy class of an associated mapping 
of $\phi:B_2\rightarrow PTF$ is the same as the element of 
$\pi_1(PTF)\oplus\pi_1(PTF)$ modulo the conjugation 
of both summands in it by an element of $\pi_1(PTF)$, and that 
the restriction of $\phi$ to a circle of $B_2$ 
represents an element of $\pi_1^-(PTF)$. 
One verifies that the class $f$ of an oriented $S^1$-fiber of 
$PTF\rightarrow F$ is in the center of $\pi_1(PTF)$ (cf.~\ref{commute}), 
and that any 
$\alpha \in \pi_1(PTF)$ is equal to 
$\beta f^k$, for some $\beta\in \pi_1^+(PTF)$ and $k\in\Z$. Hence the results
of the factorization of $\pi_1^-(PTF)$ 
by the actions of $\pi_1(PTF)$ and of $\pi_1^+(PTF)$
via conjugation are the same. 
After this the proof of the fact that $\psi^-$
is injective is the same as for $\psi^+$.

To prove that $\psi$ is injective
it suffices to show that if
$L_1, L_2\in \mathcal C$ are $T$-equivalent, then they are $T_i$-equivalent
(see~\ref{interpretation}). 
The elements 
$(\delta_1, \delta_2, \delta_3,i), (\beta_1, \beta_2, \beta_3, i)\in R_T$
corresponding to them can be chosen, so  
that $\delta_1=\beta_1 f_2^{i_1}, \delta_2=\beta_2
f_2^{i_2}, \delta_3=\beta_3 f_2^{i_3}$. Using the action of $\Z ^3$ we 
can change the element corresponding to $L_1$, so that $\delta_1=\beta_1$, 
$\delta_2=\beta_2$, and $\delta_3=\beta_3 f_2^k$ for some $k\in \Z$.
To prove the proposition it suffices to show that $f_2^k=1$.
 
Since $L_1$ and $L_2$ belong to the same component of $\mathcal L$, the
$h$-principle implies that $\beta=\beta_1 \beta_2 \beta_3$ and 
$\beta'=\beta_1 \beta_2\beta_3f_2^k$ are conjugate in $\pi_1(STF)$. 
Hence 

\begin{equation}\label{relTTi}
\beta \xi=\xi \beta f^k,
\end{equation} 
for some $\xi\in \pi_1(STF)$.

If $F=S^2$, then $\pi_1(STF)=\Z_2$, which implies that $f^k=1$. 

If $F=T^2$, then $\pi_1(STF)=\Z\oplus\Z\oplus\Z$, which implies that $k=0$
and $f^k=1$. 

For $F\neq S^2, T^2$ we note that projections of $\xi$ and $\beta$ to $F$ 
commute in $\pi_1(F)$. From Propositions~\ref{Preissman} and~\ref{commute} 
we get that $\beta=\alpha^if_2^j$ 
and $\xi=\alpha^m f_2^n$, for some $i,j,m,n\in
\Z$ and $\alpha\in \pi_1(STF)$. Substituting these expressions
into~\eqref{relTTi} and using~\ref{commute} we get that $f^k=1$.
(Recall that $F$ is assumed to be orientable.)
This finishes the proof of Proposition~\ref{orient}. 
\qed

\section{Proof of Theorem~\ref{universal}}\label{pfuniversal}
The $\overline{I^+}$ invariant corresponds to some $G$-valued 
function on the set of connected components of the normalization of the
dangerous self-tangency part of the discriminant. (This function  determines
$\overline{I^+}$ in $\mathcal C$ up to the choice of an additive constant.)
The connected components of the normalization of the dangerous self-tangency
part of the discriminant in $\mathcal C$ are in the natural one-to-one
correspondence with the $K^+$ equivalence classes of  
fronts in $K^+\cap\mathcal C$, see~\ref{orient}. Let $R^+$ be the 
factor of the set $\pi_1(STF)\oplus\pi_1(STF)$
modulo the actions of $\Z_2$ acting by permutation of summands and of 
$\pi_1(STF)$ acting by conjugation of both summands. The $h$-principle says
that the component of $\mathcal L$ containing $l$ 
is defined by a conjugacy classes realized 
by $\vec l$ in $\pi_1(CSTF)$ or, which is the same, by the conjugacy class of 
$l\in \pi_1(STF)$ and the Maslov index of $L$. 
Thus the set of $\mathcal K^+$-equivalence classes of fronts in $\mathcal C$ 
is naturally identified 
with the subset $R^+_{\mathcal C}\subset R^+$ whose elements are represented
by $(\alpha_1, \alpha_2)\in\pi_1(STF)\oplus\pi_1(STF)$ such that
$\alpha_1\alpha_2$ is conjugate to $l$ in $\pi_1(STF)$. (Here $L$ is a front
from $\mathcal C$.)
We denote by $\overline{R^+}$ the set which is a factor of 
$\pi_1(F)\oplus\pi_1(F)$ modulo the actions of $\Z_2$ permuting the summands
and of $\pi_1(F)$ acting by conjugation. Since $f_2$ is in the center of 
$\pi_1(STF)$ (see~\ref{commute})
and generates $\ker \pr^2_*$ we get that $\pr^2_*$ induces 
the natural mapping $p:R^+_{\mathcal C}\rightarrow\overline{R^+}$. 
(It is the projection on each summand.)
 
We denote by $\Z[R^+_{\mathcal C}]$ the free 
$\Z$-module of formal finite integer linear
combinations of the elements of $R^+_{\mathcal C}$. We denote by 
$[s_1, s_2]$ the element of $R^+_{\mathcal C}$ realized by $(s_1,
s_2)\in\pi_1(STF)\oplus\pi_1(STF)$.

Let $\gamma$ be a generic path
in $\mathcal C$ connecting $L_1$ to $L_2$. To prove the Theorem it suffices
to show that if
$I^+(L_1)=I^+(L_2)$, then 
the algebraic sum of the signs of crossings of $\gamma$ with the 
part of the $K^+$-stratum corresponding to $k^+\in \mathcal K^+$ is zero for
every $k^+$.

Consider the homomorphism  
$g:\Z[R^+_{\mathcal C}]\rightarrow\Z\mathcal [R^+_{\mathcal C}]$
which maps $[s_1, s_2]$ to $2[s_1, s_2]-[s_1f_2,
s_2f_2^{-1}]-[s_1f_2^{-1}, s_2f_2]$. (This homomorphism is induced 
by the behaviour of $I^+$ under crossings of the $K^+$-stratum.)
To prove the Theorem it suffices to show that $\ker g=0$.
One verifies that $\Z[R^+_{\mathcal C}]$ splits into a direct sum over 
$\Im p(R^+_{\mathcal C})$ of submodules
which are finite linear combinations of elements of $R^+_{\mathcal C}$ 
projecting to the same element of $\overline{R^+}$. Clearly $g$ maps every 
summand to itself. Thus it suffices
to show that the restriction of $g$ to every summand has trivial kernel.

Fix $r^+\in [R^+_{\mathcal C}]$. Below we construct the ordering on
$p^{-1}(p(r^+))$, which makes it isomorphic (as an ordered set) 
to $\N$ or to $\Z$ (depending on $r^+$). One verifies that the matrix 
of the restriction of $g$ to $\Z[p^{-1}(p(r^+))]$
written with respect to the basis which is the
ordered $p^{-1}(p(r^+))$ is tridiagonal with all nonzero entries on the
diagonal below the main one. This implies the statement of the
Theorem.

To construct the ordering on $p^{-1}(p(r^+))$ we need the following
technical proposition.

\begin{prop}\label{technical}
Let $F\neq S^2$ be an orientable surface and $\alpha_1, \alpha_2$ elements of 
$\pi_1(STF)$.

a) $\alpha_1, \alpha_2$ commute in $\pi_1(STF)$ provided  
that $\pr^2_*(\alpha_1)$ and $\pr^2_*(\alpha_2)$ commute in $\pi_1(F)$.

b) If $pr^2_*(\alpha_1)$ and $\pr^2_*(\alpha_2)$ are conjugate in
$\pi_1(F)$, then there exists a unique $i\in\Z$ such that $\alpha_1$ and 
$\alpha_2 f_2^i$ are conjugate in $\pi_1(STF)$.

c) Let $\beta_1, \beta_2\in\pi_1(STF)$ be such that 
$(\delta \alpha_1 \delta^{-1}, \delta \alpha_2 \delta^{-1})=(\beta_1,
\beta_2)\in\pi_1(STF)\oplus\pi_1(STF)$ for some $\delta\in\pi_1(STF)$.
If there exists $\xi\in\pi_1(F)$ such that $\xi \pr^2_*(\alpha_1) \xi
^{-1}=\pr^2_*(\alpha_2)$ and $\xi \pr^2_*(\alpha_2)
\xi^{-1}=\pr^2_*(\alpha_1)$, then 
$\pr^2_*(\alpha_1)=\pr^2_*(\alpha_2)$, $\pr^2_*(\beta_1)=\pr^2_*(\beta_2)$; 
and hence there exist unique $i,j\in\Z$ such that 
$\alpha_1=\alpha_2 f_2^i$, $\beta_1=\beta_2f_2^j$. Moreover $i=j$.
\end{prop}

The proof of the proposition is straightforward. It is based on 
Propositions~\ref{Preissman} and~\ref{commute} and the facts that  
$f_2$ generates $\ker \pr^2_*$ and that $\pi_1(STT^2)=\Z^3$.

For every $r^+\in R^+_{\mathcal C}$ 
the set $p^{-1}(p(r^+))$ has a natural ordering such that 
as an ordered set it is isomorphic to either $\N$ or $\Z$.

The ordering is constructed as follows: 

a) If $r^+$ can be realized as $(\alpha_1, \alpha_2)$ such that 
$\xi \pr^2_*(\alpha_1) \xi^{-1}=\pr^2_*(\alpha_1)$ and 
$\xi \pr^2_*(\alpha_2) \xi^{-1}=\pr^2_*(\alpha_1)$, 
for some $\xi\in\pi_1(F)$, 
then any realization of an element of $p^{-1}(p(r^+))$ has this property.
From~\ref{technical}.c we get that 
every element of $p^{-1}(p(r^+))$ determines
a unique $i\in \N$ such that $k^+$ can be realized as $(\alpha_1,
\alpha_2)$ with $\alpha_1 f_2^i=\alpha_2$. One verifies that these
natural numbers are different for different elements of $p^{-1}(p(r^+))$.
The ordering on $p^{-1}(p(r^+))$ is induced by the magnitude of $i\in\N$ 
and it makes $p^{-1}(p(k^+))$ isomorphic to $\N$.

b) If $r^+$ can not be realized as an element of the type described above,
then none of the elements of  $p^{-1}(p(r^+))$ can. This allows us to
distinguish one loop of $p(r^+)$. We use the $\Z_2$ action on 
$\pi_1(STF)\oplus\pi_1(STF)$ (used to introduce $R^+$) to interchange the
two loops, so that the first loop projects to the distinguished loop of
$p(r^+)$. We get that  
every element of $p^{-1}(p(r^+))$ 
can be realized in a unique way as an element of 
the set $R$ which is the factor of 
$\pi_1(STF)\oplus\pi_1(STF)$ modulo the action of $\pi_1(STF)$ by
conjugation of both summands. 
If $(s_1, s_2)$ and $(s_3, s_4)\in R$ realize two elements of
$p^{-1}(p(r^+))$, then 
there exists a unique $i\in\Z$ such that $s_1 f_2^i$ is conjugate to $s_3$,
see~\ref{technical}.b. As it was said in the beginning of the proof, 
$s_1s_2$ and $s_3s_4$ are conjugate in $\pi_1(STF)$, since they correspond
to nongeneric fronts from the same connected component of $\mathcal L$. 
One uses this to verify that 
if $i=0$, then $(s_1, s_2)$ and $(s_3, s_4)$ realize the same 
element of $p^{-1}(p(r^+))$. The 
ordering on $p^{-1}(p(r^+))$ is induced by the magnitude of $i$ and it makes 
$p^{-1}(p(r^+))$ isomorphic to $\Z$. 

This finishes the proof of Theorem~\ref{universal}. \qed

\section{Proof of theorems describing homotopy groups of $\mathcal L$}
\subsection{Proof of Theorems~\ref{pi1S2},~\ref{pi1T2},~\ref{pi1RP2},
and~\ref{pi1K}}\label{pf4surfaces}
Propositions~\ref{centralizer} and~\ref{pi1CSTF} 
reduce the proof of the Theorems to the calculations of the centralizer
of $l\in \pi_1(STF, l(a))$ and of the subgroup of it consisting of
elements projecting to orientation-preserving loops in $F$.

One verifies that $\pi_1(STS^2)=\Z^2$,
$\pi_1(STT^2)=\Z\oplus\Z\oplus\Z$, $\pi_1(ST\R P^2)=\Z_4$, and that 
the generator of $\pi_1(ST \R P^2)$
projects to an orientation-reversing loop in $\R P^2$.

The centralizers (and the generators of them) in the case where $F$ is the
Klein bottle are described in~\ref{K}.

One verifies that: 

{\bf a:} The subgroup of $\Z\oplus\Z_4$, which is
$(\Z^{ev}\oplus\Z_4^{ev})\cup(\Z^{odd}\oplus\Z_4^{odd})$, is generated by 
$\{(1, 1),(0,2)\}$ and is isomorphic to $\Z\oplus\Z_2$.

{\bf b:} The subgroup of $\Z\oplus\Z$, which is $(\Z^{ev}\oplus\Z^{ev})
\cup(\Z^{odd}\oplus\Z^{odd})$, is generated by $\{(1,1), (0,2)\}$ and is
isomorphic to $\Z\oplus\Z$.
 
{\bf c:} The subgroup of $\Z\oplus\pi_1(STK)$ which is
$(\Z^{ev}\oplus\pi_1^{pres}(STK))\cup(\Z^{odd}\oplus\pi_1^{rev}(STK))$, is
generated by $\{(2, 1), (0, l_1), (1, l_2), (0, f_2)\}$ and is isomorphic to 
$\Z\oplus\pi_1(STK)$. (Here $l_1, l_2, f_2$ are the generators of
$\pi_1(STK)$, see~\eqref{present}.)
 
Combining these results with Propositions~\ref{centralizer}
and~\ref{pi1CSTF} we get the proofs of the theorems.
\qed

\subsection{Proof of Theorem~\ref{pi1orient} and
Theorem~\ref{pi1nonorient}.}\label{pfpi1arbitrary}

We are going to prove, that the statement of Theorem~\ref{pi1nonorient}
is true for any orientable surface $F\neq S^2, T^2$ and any
nonorientable $F\neq \R P^2, K$. (We will see that
$l\in\pi_1(STF, l(a))$ can be presented in the unique way
as $l_g^k f_2^m\in\pi_1(STF, l(a))$ for any $F\neq S^2, \R P^2,
T^2, K$.)

Clearly this gives a proof of Theorem~\ref{pi1nonorient}.
Theorem~\ref{pi1orient} is also an immediate consequence of this fact.

\begin{emf}\label{xineq0}
{\em Proof of Theorems~\ref{pi1orient} 
and~\ref{pi1nonorient} in the case of $L\neq 1\in\pi_1(F, L(a))$.\/}
Proposition~\ref{centralizer} says that $\pi_1(\mathcal L, L)$ is isomorphic to 
$Z(\vec l)$, the centralizer of $\vec l\in \pi_1(CSTF, \vec l(a))$. 
Proposition~\ref{pi1CSTF} allows us to reduce the calculation of 
$Z(\vec l)$ to the calculation of $Z(l)<\pi_1(STF, l(a))$, which is done
below.

Consider a subgroup $G'$ of $\pi_1(F, L(a))$ generated by $L$. It is
an infinite cyclic group (see~\ref{Preissman}).
There is a unique (see~\ref{Preissman}) maximal infinite cyclic
group $G<\pi_1(F, L(a))$ containing $G'$.
Let $g$ be the generator of $G$. Let $L_g$ be a front such that 
$L_g(a)=L(a)$, and $l_g=g\in \pi_1(STF, l(a))$.

Take $\alpha\in Z(l)$.
Since $l$ and $\alpha$ commute in
$\pi_1(STF,l(a))$ we get that their projections to $F$
commute in $\pi_1(F,L(a))$.
Proposition~\ref{Preissman} implies that these projections are in the subgroup
$G$. The kernel of the homomorphism $\pr^2_*$ is generated by $f_2$, which has
infinite order in $\pi_1(STF)$ for our surfaces $F$.
This fact and Proposition~\ref{commute}
imply that there exist unique $i,j,m,n\in \Z$ such that
$g=l_g^i f_2^j$ and $\alpha=l_g^m f_2^n$.

Using Proposition~\ref{commute} we find all 
values of $k,l,m,n$ such that the elements $\alpha$ and $l$ commute.
This allows us to calculate
$Z(l)$. It turns out to be:

{\bf a:} A group isomorphic to $\Z\oplus\Z$ generated by
$\{l_g, f_2\}$, provided that $g$ is an orientation-preserving
loop.

{\bf b:} A group isomorphic to $\Z$ generated by $l_gf_2^j$,
provided that $g$ is an orientation-reversing loop and that $i$ is odd.
(This means that $L$ is an orientation-reversing front.)
Note also that in this case $(l_g f_2^j)^2=l_g^2$.

{\bf c:} A group isomorphic to $\Z\oplus\Z$ generated by $\{l_g^2, f_2\}$,
provided that $g$ is an orientation-reversing loop, $i\neq 0$
is even, and $j\neq 0$.

{\bf d:} A group isomorphic to $\pi_1(K)$ generated by $\{l_g, f_2\}$,
provided that $g$ is an orientation-reversing loop, $i\neq 0$
is
even and $j=0$.

(Note that if $i=0$, then $L=1\in\pi_1(STF)$, which contradicts
to
our assumption.)

One verifies that:

{\bf a:} The subgroup of $\Z\oplus\Z$, which is
$(\Z^{ev}\oplus\Z^{ev})\cup(\Z^{odd}\oplus\Z^{odd})$, is generated by 
$\{(1, 1), (2,0)\}$ and is isomorphic to $\Z\oplus\Z$.

{\bf b:} The subgroup of $\Z\oplus\pi_1(K)$ which is
$(\Z^{ev}\oplus\pi_1^{pres}(K))\cup (\Z^{odd}\oplus\pi_1^{rev}(K))$ is 
generated by $\{(2, 1), (0, b), (1, c)\}$ and is isomorphic to
$\Z\oplus\pi_1(K)$. (Here $b$ and $c$ are the generators of $\pi_1(K)=\{b,
c\big |bc=cb^{-1}\}$.)

Combining these results with Propositions~\ref{centralizer}
and~\ref{pi1CSTF} we obtain the proof of the two theorems for this case.
\end{emf}

\begin{emf}\label{xieq0}
{\em Proof of Theorems~\ref{pi1orient} 
and~\ref{pi1nonorient} in the case of $L=1\in\pi_1(F,
L(a))$.\/}
The kernel of $\pr^2_*$ is generated by $f_2$. 
Since $L=1\in\pi_1(F, L(a))$ we get that $l=f_2^k\in\pi_1(STF,
l(a))$, for some $k\in\Z$. We calculate the centralizer $Z(l)=Z(f_2^k)$ 
of $l\in\pi_1(STF, l(a))$.

For the case of $k\neq 0$ Proposition~\ref{commute} implies that
$Z(l)=Z(f_2^k)$ coincides with 
$\pi_1^{pres}(STF, l(a))$.
If $k=0$, then 
$l=1\in \pi_1(STF, l(a))$ and $Z(l)=\pi_1(STF, l(a))$. 

Combining these results and Propositions~\ref{centralizer}
and~\ref{pi1CSTF} we obtain the proof of the two theorems for this case.
\qed
\end{emf}

\subsection{Proof of Theorem~\ref{pin}.}\label{pfpin}
The space $\mathcal L$ is weak 
homotopy equivalent to the space $\Omega CSTF$ of all free loops in $CSTF$, 
see~\ref{H-principle}.

Consider a fibration of the space $\Omega CSTF$ over $CSTF$. 
The fiber $\Omega_x CSTF$ 
of the fibration over $x\in CSTF$ consists of all loops
$\omega:S^1\rightarrow CSTF$ such that $\omega(a)=x$.

We obtain the following exact sequence:
 
\begin{equation}\label{temp}
\cdots \stackrel{\p}{\rightarrow}
\pi_n(\Omega_{\vec l(a)}CSTF,\vec l)\stackrel{in_*}{\rightarrow}
\pi_n(\Omega CSTF, \vec l)\stackrel{t_*}{\rightarrow}
\pi_n(CSTF, \vec l(a))\stackrel{\p}{\rightarrow}\cdots.
\end{equation}

\begin{lem}\label{splitgroup}
If $F$ is equal to $S^2$ or $\R P^2$ and $n\geq 2$, then
\begin{equation}
\pi_n(\Omega CSTF, \vec l)=\pi_n(\Omega_{\vec l(a)}CSTF, \vec l)\oplus
\pi_n(CSTF, \vec l(a)).
\end{equation}
\end{lem}

\begin{emf}{\em Proof of Lemma~\ref{splitgroup}.\/}
Fix $n>1$. We construct a homomorphism
$g:\pi_n(CSTF, \vec l(a))\rightarrow\pi_n(\Omega CSTF, \vec l)$ such that
$t_*\circ g =\id_{\pi_n(CSTF, \vec l(a))}$. Since the sequence~\eqref{temp}
is exact and the groups are Abelian, the existence of such $g$ implies 
the statement of the lemma.

We describe this construction for $F=\R P^2$. The construction of $g$ for
$F=S^2$ can be easily deduced from this one.
From the exact homotopy sequences of the fibrations 
$CST \R P^2\rightarrow ST\R  P^2$ and $STS^2\rightarrow ST\R P^2$
we get that $\pi_n(CST \R P^2)$, $\pi_n(ST \R P^2)$, and
$\pi_n(STS^2)$, $n\geq 2$, are canonically isomorphic.

Take $s:S^n\rightarrow ST \R P^2$ that corresponds under these
isomorphisms to a given element
of $\pi_n(CST \R P^2, \vec l(a))$.
Let $s':S^n\rightarrow STS^2$ be the mapping which
is a lifting of $s$ under the covering $STS^2\rightarrow ST\R P^2$.
Fix an orientation on $S^2$.
Then for every $x\in S^n$ the local orientation at 
$\pr^2 (s'(x))\in S^2$ induces a local orientation at 
$\pr^2 (s(x))\in \R P^2$.

There is a unique isometric autodiffeomorphism
$I_x$ of $\R P^2$ such that:

a) it maps $\pr^2 (s(*))$ to $\pr^2(s(x))$;

b) the differential of it sends $s(*)$ to $s(x)$;

c) the local orientation at $\pr^2(s(x))$, which is described above,
coincides with the
one induced by the differential of $I_x$ from the local orientation at
$\pr^2(s(*))$.

Let $\bar s:S^n\rightarrow \Omega CST \R P^2$ be the mapping
that sends $x\in S^n$ to $h(I_x(l))$ (the lifting to $CST \R P^2$ of the 
translation of $l$ by $I_x$).

Set the value of $g$ on the element of
$\pi_n(CST\R P^2, \vec l(a))$ corresponding to $s$ to be the element of
$\pi_n(\Omega CST\R P^2, \vec l)$ represented by $\bar s$.
A straightforward verification
shows that $g$ is the desired
homomorphism from $\pi_n(CST \R P^2, \vec l(a))$
to $\pi_n(\Omega CST\R P^2, \vec l)$.
This finishes the proof of Lemma~\ref{splitgroup}.
\qed
\end{emf}

\begin{emf}
One verifies that
$\pi_2(CSTF)=0$ and $\pi_n(CSTF)=\pi_n(S^2)$, $n\geq 3$, for $F$ equal to
$S^2$ or $\R P^2$. 
Now Lemma~\ref{splitgroup}, isomorphism
$\pi_n(\Omega_{\vec l(a)}CSTF,\vec l)=
\pi_{n+1}(CSTF, \vec l(a))$, and the
weak homotopy equivalence given by the
$h$-principle (see~\ref{H-principle}) imply the
first statement of the Theorem. (Note that
$\pi_3(S^2)=\Z$.)

One verifies that $\pi_n(STF)=0$, $n\geq2$, for $F\neq S^2, \R P^2$.
The exactness of sequence~\eqref{temp} and the isomorphism 
$\pi_n(\Omega_{\vec l(a)}CSTF,\vec l)=
\pi_{n+1}(CSTF, \vec l(a))$ imply
that $\pi_n(\Omega CSTF, \vec l)=0$, $n\geq 2$. 
Using the weak homotopy equivalence
given by the $h$-principle we get the second
statement of the Theorem.
This finishes the proof of Theorem~\ref{pin}.
\qed
\end{emf}

\centerline{}
\centerline{\bf Acknowledgments}
Some results of this paper formed a part of my Ph.D.
thesis~\cite{Tchernovthesis} at Uppsala University. Other results were
obtained during my work at the ETH Z\"urich. I am gratefull to the Uppsala
University  and the ETH Z\"urich for the hospitality and support.
I am deeply grateful to Oleg Viro, Tobias Ekholm,
Viktor Goryunov, and 
Maxim Kazarian for many enlightening discussions.

%
%
% 

%\bibitem{Tchernov2}
%V.~Tchernov, {\em Strangeness- and $J^{\pm}$-type invariants of
%immersed oriented curves}, Preprint,
%Dept. of Math., Uppsala Univ., U.U.D.M. report 1997:14.
%\bibitem{Whitney}
%H.~Whitney, {\em On regular closed curves in the plane}, Composito Math.
%{\bf 4} (1936) pp 276--284


\begin{thebibliography}{99999}
\bibitem{Aicardidiscr}
F.~Aicardi: Discriminants and local invariants of planar fronts,
{\em The Arnold-Gelfand math. seminars: Geometry and Singularity theory},
AMS, Birkh\"auser, Boston, (1996), pp. 1--76
\bibitem{Aicardi}
F.~Aicardi: Invariant polynomial of framed knots in the solid torus and
its application to wave fronts and Legendrian knots, 
{\em J. of Knot Th. and Ramif.}, Vol. {\bf 15}, No. 6, (1996), pp. 743--778 
\bibitem{Arnoldtopology} V.I.~Arnold: Topological invariants of plane
curves and caustics, {\em University Lecture Series}, {\bf 5}, (1994), AMS
Providence RI. 
\bibitem{Arnold}
V.I.~Arnold: Plane curves, their invariants, perestroikas and
classifications, Singularities and Bifurcations (V.I.~Arnold, ed.) 
{\em Adv. Sov. Math.}, Vol. {\bf 21}, (1994), pp. 39--91
\bibitem{Arnoldsplit}
V.I.~Arnold: Invarianty i perestroiki ploskih frontov, {\em Trudy Mat.
Inst. Steklov.} {\bf 209} (1995); English
translation: Invariants and perestroikas of wave fronts on the plane,
Singularities of smooth mappings with additional structures,
{\em Proc. Steklov Inst. Math.}, Vol. {\bf 209}, (1995), pp. 11-64.
\bibitem{Docarmo}
M.~do~Carmo: {\em Riemannian Geometry}, Birkh\"auser, Boston (1992)
\bibitem{Duchamp}
T.~Duchamp:{\em The Classification of Legendre Immersions}, preprint (1982),
revised (1996), http://www.math.washington.edu/ duchamp/ 
\bibitem{Gromov}
M.~Gromov: {\em Partial Differential relations}, Springer-Verlag, Berlin
Heidelberg, (1986)
\bibitem{Hansen}
V.L.~Hansen: On the fundamental group of the mapping space,
{\em Compositio Math.}, Vol. {\bf 28} (1974), pp. 33--36
\bibitem{Inshakov} 
A.~V.~Inshakov: Invarianty tipa $j^+, j^-, st$ gladkih krivyh na
dvumernyh mnogoobrazijah, Preprint (1997); to appear in {\em Funct. Anal.
Appl.} 
\bibitem{Inshakov2}
A.~V.~Inshakov: Gomotopicheskie gruppy krivyh na dvumernyh
mnogoobraziyah, {\em Uspehi Mat. Nauk}, vol. 53, N 2 (1998), pp. 147-148
\bibitem{Kalfagianni}
E.~Kalfagianni: Finite type invariants for knots in $3$-manifolds, {\em
Topology}, Vol. 37, No. 3, pp. 673-707, 1998
\bibitem{Polyak}
M.~Polyak: Invariants of curves and fronts via Gauss diagrams;
to appear in {\em Topology}
\bibitem{Polyakpriv}
M.~Polyak: Private communication
\bibitem{Polyakbennequin}
M.~Polyak: On the Bennequin invariant of Legendrian curves and its
quantization, {\em C. R. Acad. Sci. Paris Ser. I Math.} {\bf 322} (1996), No. 1,
pp. 77--82 
\bibitem{Tchernov1}
V.~Tchernov: {\em Arnold-type invariants of curves on surfaces and
homotopy groups of the space of curves},
Preprint, Dept. of Math., Uppsala Univ., U.U.D.M. report 1997:21, to appear 
as: Arnold-type invariants of curves on surfaces, {\em J. of Knot Th. and
Ramif.}, and Homotopy groups of the space of curves on a surface,
{\em Math. Scand.}
\bibitem{Tchernovthesis}
V.~Tchernov: Arnold-type invariants of curves and wave fronts on
surfaces, {\em Uppsala Dissert. in Math.}, {\bf 10}, (1998)
\bibitem{Tchernovamsvolume}
V.~Tchernov: Shadows of wave fronts and
Arnold-Bennequin type invariants of fronts
on surfaces and orbifolds, {\em AMS volume Adv. in Math. Sci.}, to appear
\end{thebibliography}
\end{document}